\documentclass[11pt]{amsart}
\usepackage{etex}
\usepackage[usenames,dvipsnames]{pstricks}
\usepackage[mathscr]{eucal}
\usepackage{amsfonts,amsmath,amssymb,amsthm,amscd, amsxtra}
\usepackage{enumerate,verbatim}
\usepackage[all,2cell,ps]{xy}
\usepackage[notcite, notref]{}
\usepackage[pagebackref]{hyperref}
\usepackage[normalem]{ulem}
\usepackage{calc,color}
\usepackage{wasysym}

\usepackage{mathtools,mathptmx}
\usepackage{enumerate, epsfig,verbatim}
\usepackage{centernot}
\usepackage{todonotes}
\usepackage[T1]{fontenc}
\usepackage[usenames,dvipsnames]{pstricks}
\usepackage[mathscr]{eucal}
\usepackage[notcite,notref]{}
\usepackage{tikz}
\usetikzlibrary{graphs,graphs.standard,matrix, calc, arrows,quotes}

\usepackage{nicefrac}
\usepackage{array}
\usepackage{xcolor}
\usepackage[all]{xypic}

\usepackage{url}
\usepackage{datetime}
\theoremstyle{definition}
\newtheorem{thm}{Theorem}[section]

\newtheorem{prop}[thm]{Proposition}
\newtheorem{cor}[thm]{Corollary}

\theoremstyle{definition}
\newtheorem{lem}[thm]{Lemma}

\newtheorem{eg}[thm]{Example}

\newtheorem{ques}[thm]{Question}

\newtheorem{rmk}[thm]{Remark}

\newtheorem{chunk}[thm]{\hspace*{-1.065ex}\bf}
\numberwithin{equation}{section}
\newcommand{\fm}{\mathfrak{m}}
\newcommand{\fn}{\mathfrak{n}}
\newcommand{\fp}{\mathfrak{p}}
\newcommand{\fq}{\mathfrak{q}}
\newcommand{\fa}{\mathfrak{a}}

\newcommand{\RR}{\mathbb{R}}
\newcommand{\CC}{\mathbb{C}}

\newcommand{\ZZ}{\mathbb{Z}}
\newcommand{\NN}{\mathbb{N}}

\DeclareMathOperator{\Ass}{Ass}

\DeclareMathOperator{\CI}{\textnormal{CI-dim}}
\DeclareMathOperator{\Gdim}{\textnormal{G-dim}}
\DeclareMathOperator{\gr}{\textnormal{grade}}
\DeclareMathOperator{\cx}{cx}
\DeclareMathOperator{\Tr}{Tr}
\DeclareMathOperator{\X}{X}

\DeclareMathOperator{\V}{V}
\DeclareMathOperator{\IH}{E}
\DeclareMathOperator{\mz}{Z}
\DeclareMathOperator{\mw}{W}

\DeclareMathOperator{\coker}{coker}
\DeclareMathOperator{\depth}{depth}

\DeclareMathOperator{\edim}{embdim}
\DeclareMathOperator{\Ext}{Ext}

\DeclareMathOperator{\Ht}{ht}
\DeclareMathOperator{\hh}{H}

\DeclareMathOperator{\E}{E}
\DeclareMathOperator{\Hom}{Hom}

\DeclareMathOperator{\im}{im}

\DeclareMathOperator{\len}{length}

\DeclareMathOperator{\pd}{pd}

\DeclareMathOperator{\rank}{rank}

\DeclareMathOperator{\Supp}{Supp}
\DeclareMathOperator{\Spec}{Spec}

\DeclareMathOperator{\Syz}{\Omega}
\DeclareMathOperator{\Tor}{Tor}

\DeclareMathOperator{\cod}{codim}

\DeclareMathOperator{\NF}{NF}

\author{Mohsen Asgharzadeh}
\address{Mohsen Asgharzadeh\\ Hakimieh, 16599-19556, Tehran, IRAN}
\email{mohsenasgharzadeh@gmail.com}

\author{Olgur Celikbas}
\address{Olgur Celikbas\\School of Mathematical and Data Sciences\\
West Virginia University\\ Morgantown, WV 26506-6310, USA}
\email{olgur.celikbas@math.wvu.edu}

\author[Arash Sadeghi]{Arash Sadeghi}
\address{Arash Sadeghi\\Dokhaniat 49179-66686, Gorgan, IRAN}
\email{sadeghiarash61@gmail.com}

\begin{document}

\title[Tate homology and local cohomology of tensor products of modules]{Studies on a relationship between Tate homology and \\local cohomology of tensor products of modules}



\begin{abstract}
We study the local cohomology of finitely generated modules with respect to specialization-closed subsets over commutative Noetherian rings. Our main result establishes a relationship between the local cohomology of tensor products of modules and Tate homology over Gorenstein rings. We present several consequences of our main result, which are of independent interest. For example, we study various aspects of the cohomological rigidity property of tensor products of modules and generalize results of Huneke--Wiegand in this direction. We determine conditions under which the vanishing of a single local cohomology module of a tensor product implies the vanishing of all lower ones, establish bounds for the depth of tensor products of modules, and obtain criteria for freeness of modules over complete intersection rings. Along the way, we also provide a splitting criterion for vector bundles on smooth complete intersections.
\end{abstract}

\maketitle{}

{\renewcommand\thefootnote{}\footnotetext{This is an expanded and revised version of a portion of the unpublished manuscript~\cite{AOM} by the authors.}}

\setcounter{tocdepth}{1}
\tableofcontents

\section{Introduction} Throughout, $R$ denotes a commutative Noetherian ring. If $R$ is local, we denote its maximal ideal by $\fm$ and its residue field by $k$. All $R$-modules are assumed to be finitely generated.

The primary aim of this paper is to study the local cohomology of tensor products of modules with respect to specialization-closed subsets of $\Spec(R)$; see~\ref{lsc} and~\ref{LCo}. Our main result, Theorem~\ref{t12}, establishes a new relationship between Tate homology and the local cohomology of tensor products; see Section~\ref{mainsection}. A special case of Theorem~\ref{t12}---useful for applications---can be stated as follows:

\begin{thm}\label{main-thm-intro} Let $R$ be a $d$-dimensional Gorenstein local ring. Assume $R$ is an isolated singularity, that is, $R_{\fp}$ is regular for all $\fp \in \Spec(R)-\{\fm\}$. If $M$ and $N$ are maximal Cohen-Macaulay $R$-modules, then $\hh^i_{\fm}(M\otimes_RN)\cong\widehat{\Tor}_{-i}^R(M,N)$ for all $i=0, \ldots, d-1$, and $\widehat{\Tor}_{-d}^R(M,N)$ embeds into $\hh^d_{\fm}(M\otimes_RN)$.
\end{thm}

Theorem~\ref{t12} establishes a more general result than Theorem~\ref{main-thm-intro}. It concerns the local cohomology of tensor products of certain modules with respect to specialization-closed subsets of $\Spec(R)$ (as introduced by Grothendieck~\cite{sga2}) and applies to rings that are not necessarily Gorenstein.

We present various corollaries of Theorem~\ref{t12}. One of these is Proposition~\ref{ctc-intro}, which allows us to compute the depth of tensor products over complete intersections using Tate homology under certain conditions. In what follows, $\NF(M)$ denotes the non-free locus of an $R$-module $M$; we say that $M$ is \emph{locally free on the punctured spectrum of $R$} if $\NF(M) \subseteq \{\fm\}$.

\begin{prop} \label{ctc-intro} Let $R$ be a local complete intersection ring and let $M$ and $N$ be $R$-modules such that $\NF(M)\cap\NF(N)\subseteq\{\fm\}$. Assume $\depth_R(M)+\depth_R(N)-\depth(R) \geq  \cod(R)$.
Then at least one of the following conditions holds:
\begin{enumerate}[\rm(i)]
\item $\depth_R(M)+\depth_R(N)-\depth(R)=\depth_R(M\otimes_RN)$. \smallskip
\item $\depth_R(M)+\depth_R(N)-\depth(R)>\depth_R(M\otimes_RN)=\inf \left\{i\geq0\mid\widehat{\Tor}_{-i}^R(M,N)\neq0\right\}$.
\end{enumerate}
\end{prop}

Proposition~\ref{ctc-intro} provides a method for computing the depth of tensor products of modules without assuming the vanishing of (co)homology, a condition often required in earlier results. It is a consequence of Corollary~\ref{ctc-main}, which in turn follows from Proposition~\ref{ct}. The latter concerns maximal Cohen-Macaulay modules over Gorenstein rings that are not necessarily complete intersections and establishes relationships between the depth of their tensor products and the vanishing of Tate homology; see, for example, Corollaries~\ref{ct-cor} and~\ref{corcombined}.

Auslander~\cite{Au} and Lichtenbaum~\cite{Li} proved that if $M$ and $N$ are modules over a regular ring $R$ such that $M \otimes_R N$ is torsion-free, then $\Tor_i^R(M, N) = 0$ for all $i \geq 1$. This vanishing result fails quickly once one moves away from regular rings. Motivated by this, various researchers have studied conditions under which good depth properties of tensor products force the vanishing of $\Tor$; see, for example,~\cite{MA, Ce, CIPW, Da1, HW1, HW2, Jo1}. A standout result in this direction is a theorem of Dao; see~\ref{DD}. Using Theorem~\ref{t12}, we prove Proposition~\ref{DDmore-intro}, stated as Corollary~\ref{DDmore}, which is a variation of~\ref{DD} and assumes weaker conditions on the local freeness of the modules involved. In fact, Corollary~\ref{DDmore} yields a stronger result than Proposition~\ref{DDmore-intro}, as it incorporates the complexity of the module, a finer invariant than the codimension of the ring; see also~\ref{cx}.

\begin{prop} \label{DDmore-intro} Let $R$ be a local complete intersection ring, and let $M$ and $N$ be $R$-modules. Assume $c$ is an integer such that $c \geq  \cod(R)+1$. Assume further:
\begin{enumerate}[\rm(i)]
\item $\depth_{R_\fp}(M_\fp)+\depth_{R_\fp}(N_\fp)\geq \depth(R_\fp)+(c-1)$ for all $\fp \in \Spec(R)$ with $\Ht_R(\fp)\geq c-1$.
\item $M\otimes_RN$ satisfies $(S_{c})$.
\item $M_{\fp}$ is free over $R_{\fp}$ for all $\fp \in \X^{c-1}(R)$. 
\end{enumerate}
Then $\Tor_{i}^R(M,N)=0$ for all $i\geq 1$.
\end{prop}

Huneke and Wiegand~\cite{HW2}, motivated by the work of Kempf~\cite{Kem} and Luk\textendash Yau~\cite{Luk} on the cohomological rigidity property, studied the vanishing of local cohomology of $M \otimes_R M^{\ast}$ for certain modules $M$ over regular rings $R$, establishing freeness criteria. In Section~\ref{section5}, we revisit one of their results, stated as~\ref{HWL1}, extend it, and prove several variations. For example, Corollary~\ref{c1} extends a result of Huneke and Wiegand by replacing the regularity hypothesis with the assumption that $R$ is Cohen-Macaulay and proves the following; see also Corollaries~\ref{remzidenoturu}, ~\ref{tp3c1}, and Theorem~\ref{mainofsection} for further applications in the same vein.

\begin{prop} Let $R$ be a Cohen-Macaulay local ring and let $M$ be an $R$-module which is locally free on the punctured spectrum of $R$. If $\depth_R(M)\geq 1$ and $\hh^1_{\fm}(M\otimes_RM^*)=0$, then $M$ is free.
\end{prop}

In Sections~\ref{section4},~\ref{section5}, and~\ref{section6}, we give further applications of Theorem~\ref{t12}, which are of independent interest. For example, Corollaries~\ref{cxmm2} and \ref{corcombined} establish useful bounds on the depth of tensor products of modules and corroborate previously established results in the literature, namely~\cite[3.1]{AOMCM} and~\cite[1.4]{CST}.

The proof of our main result, Theorem~\ref{t12}, relies on several auxiliary statements. These are stated in Section~\ref{mainsection} and proved in Section~\ref{PAR}. Among them, Proposition~\ref{l2} provides a decomposition result analogous to those established by \v{C}esnavi\v{c}ius~\cite{KC} for arithmetically Cohen-Macaulay vector bundles on smooth complete intersections; see Propositions~\ref{artCM} and~\ref{t3.1}.

%%%%%%%%%%%%%%%%%%%%%%%%%%%%%%%%%%%%%%%%%%%%%%%%%%%%%%%%%%
%%%%%%%%%%% SECTION %%%%%%%%%%%%%%%%%%%%%%%%%%%%%%%%%%%%%%%%
%%%%%%%%%%%%%%%%%%%%%%%%%%%%%%%%%%%%%%%%%%%%%%%%%%%%%%%%%%
%%%%%%%%%%%%%%%%%%%%%%%%%%%%%%%%%%%%%%%%%%%%%%%%%%%%%%%%%%

\section{Preliminaries: Some definitions and notations}

In this section, we provide the definitions and properties necessary for our results and proofs, particularly those in Section~3. Additional preliminaries will be introduced in the relevant sections as we proceed.

\begin{chunk} \textbf{Setup}. As mentioned in the introduction, throughout, $R$ denotes a commutative Noetherian ring (not necessarily local). If $R$ is local, then $\fm$ and $k$ denote the unique maximal ideal and the residue field of $R$, respectively. Moreover, unless otherwise stated, each $R$-module is assumed to be finitely generated. 

By a local \emph{complete intersection} ring $R$, we mean a ring that is quotient of a power series ring over $k$ modulo a regular sequence, that is, $R=S/(\underline{x})$ for some $(S, \fn)=k[\![x_1, \ldots, x_n]\!]$ and some regular sequence $\{\underline{x}\}$ on $S$ contained in $\fn$. We say $R$ is a local \emph{hypersurface} if $\{\underline{x}\}=\{f\}$ for some $0\neq f\in \fn$.

We set $(-)^{\ast}=\Hom_R(-,R)$. Also, we assume $\depth_R(0)=\infty$, $\dim_R(0)=-1$, and $\pd_R(0)=-\infty$; see~\cite{Au}.
\end{chunk}

\begin{chunk} \textbf{Auslander transpose}. \label{a1} Let $M$ be an $R$-module. If $P$ is a projective $R$-module and $P\overset{\partial}{\twoheadrightarrow} M$ is a surjective $R$-module map, then the kernel of $\partial$ is called the first
\emph{syzygy} of $M$; it is denoted by $\Omega^1_R M$ and is unique up to projectives. Inductively,
we define the $n$th syzygy module of $M$ as $\Omega_R^nM= \Omega^1_R(\Omega_R^{n-1}M)$ for all $n\geq 1$. We simply write $\Omega_R M$ instead of $\Omega^1_R M$ and set $\Omega_R^0M=M$.

Let  $P_1\overset{\partial}{\rightarrow}P_0 \overset{}{\rightarrow} M\to0$ be a finite projective presentation of $M$, that is, an exact sequence of $R$-modules where $P_0$ and $P_1$ are projective $R$-modules. The (Auslander) \emph{transpose} of $M$, denoted by $\Tr_R M$, equals $\coker({\partial}^*)$ and hence is given by the exact sequence
$0\rightarrow M^*\rightarrow P_0^*\overset{\partial^{\ast}}{\rightarrow} P_1^*\rightarrow \Tr_R M\rightarrow 0$. This shows that $M^*\approx \Omega_R^2\Tr_R M$, where $\approx$ denotes stable isomorphism. Note that $\Tr_R M$ is unique up to stable isomorphism, that $\Tr_R\Tr_R M \approx M$, and that minimal projective presentations of $M$ yield isomorphic transposes; see~\cite{AuBr} and~\cite[3.5]{Yo}.

We use the following fact in the sequel: if $\Tor_1^R(\Tr_R M, M)=0$, then $M$ is free; see~\cite[3.7 and 3.9]{Yo}.
\end{chunk}

The notion of Gorenstein dimension, which we recall next, was introduced by Auslander~\cite{Tome1} and subsequently developed by Auslander and Bridger~\cite{AuBr}.

\begin{chunk}\textbf{Gorenstein dimension}. \label{HD} Let $M$ be an $R$-module. Then $M$ is called \emph{totally reflexive} provided that $M\cong M^{**}$  and $\Ext^i_R(M,R)=0=\Ext^i_R(M^*,R)$ for all $i\geq 1$. 

The \emph{Gorenstein dimension} of $M$, denoted $\Gdim_R(M)$, is defined to be the infimum of all integers $n\geq 0$ such that there is an exact sequence $0\rightarrow G_n\rightarrow\cdots\rightarrow G_0\rightarrow  M \rightarrow 0$ in which each $G_i$ is a totally reflexive $R$-module; see~\cite{Tome1, AuBr}. It is known that $R$ is Gorenstein if and only if each $R$-module has finite Gorenstein dimension. Moreover, the following hold:
\begin{enumerate}[\rm(i)]
\item If $\Gdim_R(M)=0$, that is, if $M$ is totally reflexive, then so is $\Tr_R M$.
\item If $R$ is local and $\Gdim_R(M)<\infty$, then $\Gdim_R(M)=\depth(R)-\depth_R(M)$; see~\cite[4.13]{AuBr}.
\item If $\Gdim_R(M)<\infty$, then there is a short exact sequence of $R$-modules $0\to M \to X \to G \to 0$, where $G$ is totally reflexive and 
$\pd_R(X)<\infty$; see~\cite[2.17]{cfh}. 
\end{enumerate}
\end{chunk}
 
\begin{chunk} \textbf{Specialization-closed subsets}.  \label{lsc} A subset $\mz\subseteq\Spec(R)$ is said to be \emph{specialization-closed} provided that the following condition holds: $\mz \owns \fp \subseteq \fq \in \Spec(R) \Longrightarrow \fq\in \mz$; see~\cite{sga2}.

We record a few basic examples of specialization-closed subsets of $\Spec(R)$:

\begin{enumerate}[\rm(i)]
\item Every closed subset of $\Spec(R)$ with respect to Zariski toplogy is specialization-closed.
\item If $R$ is a domain, then $\Spec(R)-\{0\}$ is specialization-closed subset of $\Spec(R)$.
\item If $0\leq n\leq\dim(R)$, then $\{\fp\in\Spec(R): \Ht_R(\fp)\geq n\}$ is a specialization-closed subset of $\Spec(R)$. 
\end{enumerate}
\end{chunk}

\begin{chunk}\label{LCo}\textbf{Local cohomology}. Let $M$ be an $R$-module and let $\mz$ be a specialization-closed subset of $\Spec(R)$. 
\begin{enumerate}[\rm(i)]
\item We set $\Gamma_{\mz}(M)=\{x\in M : \Supp_R(Rx)\subseteq\mz\}$. The $i$-th right derived functor $\hh^i_{\mz}(-)$ of the left exact additive functor $\Gamma_{\mz}(-)$ is called the $i$-th \emph{local cohomology functor with respect to $\mz$}. In other words, if $I^{\bullet}$ is an injective resolution of an $R$-module $M$, then $\hh^i_{\mz}(M) \cong \hh^i(\Gamma_{\mz}(I^{\bullet}))$; see~\cite[3.5.2]{BH} and~\cite{sga2}.

\item If $\mz$ is a closed subset of $\Spec(R)$, that is, if $\mz=\V(\fa)=\{\fp \in \Spec(R): \fa \subseteq \fp\}$ for some ideal $\fa$ of $R$, we denote $\Gamma_{\mz}(-)$ and $\hh^i_{\mz}(-)$ by $\Gamma_{\fa}(-)$ and $\hh^i_{\fa}(-)$, respectively.

\item Assume $R$ is local, and $M$ is nonzero and finitely generated with $\depth_R(M)=t$ and $\dim_R(M)=d$. Then 
$\hh^t_{\fm}(M)\neq 0 \neq \hh^d_{\fm}(M)$, and $\hh^i_{\fm}(M)=0$ for all $i$ with $i<t$ or $i>d$; this fact is known as Grothendieck's (non-)vanishing theorem; see~\cite[6.1.4 and 6.2.8]{Bsh} and~\cite[3.5.7]{BH}.
\end{enumerate}
\end{chunk}

In the following we record some definitions and basic properties of local cohomology modules with respect to specialization-closed subsets; we refer the reader to~\cite{Bsh,sga2} for the proofs of these results.

\begin{chunk} \label{LC} Let $\mz$ and $\mw$ be specialization-closed subsets of $\Spec(R)$. Then the following hold:
\begin{enumerate}[\rm(i)]
\item $\Gamma_{\mz}(M)=\bigcup_{\fa\in\Sigma}\Gamma_{\fa}(M)$, where $\Sigma=\{\fa \text{ is an ideal of } R : \V(\fa)\subseteq\mz\}$. 
\item An exact sequence of $R$-modules $0\to M_1\to M_2\to M_3\to 0$ yields a long exact sequence of $R$-modules: $$\cdots\to \hh^i_{\mz}(M_1)\to \hh^i_{\mz}(M_2)\to \hh^i_{\mz}(M_3)\to \hh^{i+1}_{\mz}(M_1)\to \cdots.$$ 
\item $M$ is said to be torsion with respect to $\mz$ provided that $\Gamma_{\mz}(M)=M$.
\item $M$ is torsion with respect to $\mz$, that is, $\Gamma_{\mz}(M)=M$, if and only if $\Supp_R(M)\subseteq\mz$. 
\item If $M$ is torsion with respect to $\mz$, then $\hh^i_{\mz}(M)=0$ for all $i\geq 1$.
\item $M$ is said to be torsion-free with respect to $\mz$ provided that $\Gamma_{\mz}(M)=0$.
\item $\Gamma_{\mz}\left(\Gamma_{\mw}(M)\right)=\Gamma_{\mz\cap \mw}(M)=\Gamma_{\mw}\left(\Gamma_{\mz}(M)\right)$. 
\end{enumerate}	
\end{chunk}

%%%%%%%%%%%%%%%%%%%%%%%%%%%%%%%%%%%%%%%%%%%%%%%

\begin{chunk}\label{gsc} \textbf{Grade of a module with respect to a specialization-closed subset}. Let $M$ be an $R$-module which is not necessarily finitely generated, and let $\mz$ be a specialization-closed subset of $\Spec(R)$. 

The {\em  grade} of $M$ with respect to $\mz$ is defined as:
$\gr_R(\mz,M)=\inf\{i \geq 0: \hh_{\mz}^{i}(M)\neq0\}$. Note that, if $\gr_R(\mz,M)\geq 1$, then $\Gamma_{\mz}(M)=0$, and hence $M$ is torsion-free with respect to $\mz$; see~\ref{LC}(vi). 	

Assume $M$ is finitely generated and $\mz=\V(\fa)$ for some ideal $\fa$ of $R$. Then we denote $\gr_R(\mz,M)$ by $\gr_R(\fa,M)$, which equals the common length of a maximal $M$-regular sequence contained in $\fa$. Therefore, if $R$ is local, then $\gr_R(\fm,M)$ is nothing but $\depth_R(M)$.
\end{chunk}

%%%%%%%%%%%%%%%%%%%%%%%%%%%%%%%%%%%%%%%%%%%%%%%%%%%%%%%%%%%%%%%%%%%%%%%%

\begin{chunk} {\bf Tate (co)homology}. \label{Tate} Let $T=\cdots\to  T_1 \to T_0 \to T_{-1} \to \cdots$ be a homologically indexed complex of $R$-modules $T_i$ which are not necessarily finitely generated. Then $T$ is called \emph{totally acyclic} provided that $T$ and the dual complex $\Hom_R(T,R)$ have zero homology. 

Let $M$ be an $R$-module. Then, $M$ is totally reflexive if and only if there is a totally acyclic complex $T$ such that $M\cong\coker(T_1\to T_0)$. A \emph{complete resolution} of $M$ is a diagram of complex maps $T\overset{\tau}{\to}P\overset{\simeq}{\to} M$, where $P\overset{\simeq}{\to}M$ is a projective resolution of $M$, $T$ is a totally acyclic complex of free $R$-modules, and $\tau_i$ is an isomorphism for all $i\gg 0$.  It follows that $M$ has a complete resolution if and only if $\Gdim_R(M)<\infty$; see~\cite{AM} for the details.

If $N$ is an $R$-module (not necessarily finitely generated) and $i\in\ZZ$, then the \emph{Tate homology and cohomology} of $M$ and $N$ are defined as $\widehat{\Tor}^R_i(M,N)=\hh_i(T\otimes_RN)$ and $\widehat{\Ext}^i_R(M,N)=\hh_{-i}\big(\Hom_R(T,N)\big)$, respectively. 
\end{chunk}

In the following we catalog some basic, but important, properties of (the vanishing of) Tate (co)homology; for details, see~\cite{AvBu,AM,CJ}.

\begin{chunk}\label{Tate hom} \label{Tate2} Let $R$ be a local ring and let $M$ and $N$ be $R$-modules, where $N$ is not necessarily finitely generated. Assume $\Gdim_R(M)<\infty$. Then the following properties hold:
\begin{enumerate}[\rm(i)]	
\item $\widehat{\Tor}^R_i(M,N)\cong\Tor_i^R(M,N)$ and $\widehat{\Ext}^i_R(M,N)\cong\Ext^i_R(M,N)$ for all $i>\Gdim_R(M)$.
\item If $M$ is totally reflexive, then $\widehat{\Tor}_i^R(M,N)\cong\widehat{\Ext}^{-i-1}_R(M^*,N)$ for all $i\in\ZZ$.
\item If $\pd_R(M)<\infty$, then $\widehat{\Tor}_i^R(M,N)=0=\widehat{\Ext}^i_R(M,N)$ for all $i\in\ZZ$.
\item If $0 \to X \to Y \to Z \to 0$ is an exact sequence of $R$-modules, where $X$, $Y$, and $Z$ have finite Gorenstein dimension, then there is a long exact sequence of $R$-modules for all $i\in \ZZ$:\\
$\cdots \to \widehat{\Tor}^R_{i}(X,N) \to \widehat{\Tor}^R_{i}(Y,N) \to \widehat{\Tor}^R_{i}(Z,N) \to \widehat{\Tor}^R_{i-1}(X,N) \to  \widehat{\Tor}^R_{i-1}(Y,N) \to  \cdots$
\item $\widehat{\Tor}^R_{i+n}(M,N)\cong\widehat{\Tor}^R_i(\Omega_R^n M,N)$ and $\widehat{\Ext}^{i+n}_R(M,N)\cong\widehat{\Ext}^i_R(\Omega_R^n M,N)$ for all $n\geq0$ and all $i\in\ZZ$. 
\end{enumerate}	
\end{chunk}

The depth formula, which we record next, has been proven to be an important tool in homological commutative algebra. 

\begin{chunk} \textbf{Depth formula}. \label{df} Let $R$ be a local ring and $M$ and $N$ be $R$-modules. 

Set $q=\sup\{i\geq 0: \Tor_i^R(M,N)\neq 0\}$. Assume $q<\infty$, that is,  $\Tor_i^R(M,N)=0$ for all $i\gg 0$. We say that the \emph{depth formula} holds for $M$ and $N$, or for the pair $(M,N)$, provided that the following equality holds:
$\depth_R(M)+\depth_R(N)=\depth(R)+\depth_R\big(\Tor_q^R(M,N)\big)-q$.

If $\pd_R(M)<\infty$, and if either $q=0$ or $\depth_R\left(\Tor_q^R(M,N)\right)\leq 1$, then the depth formula holds for $(M,N)$ due to Auslander~\cite[1.2]{Au}. Auslander's result was extended to complete intersection rings by Huneke-Wiegand~\cite[2.5]{HW1}: If $R$ is a local complete intersection ring and $\Tor_i^R(M,N)=0$ for all $i\geq 1$, then the depth formula holds for $M$ and $N$, even if neither $M$ nor $N$ has finite projective dimension; see~\cite{CJ} for more details on the depth formula.
\end{chunk}

%%%%%%%%%%%%%%%%%%%%%%%%%%%%%%%%%%%%%%%%%%%%%%%%%%%%%%%%%%%%%%%%
%%%%%%%%%%%%%%%%%%%%%%%%%%%%%%%%%%%%%%%%%%%%%%%%%%%%%%%%%%%%%%%%%
%%%%%%%%%%%%%%%%%%%SECTION%%%%%%%%%%%%%%%%%%%%%%%%%%%%%%%%%%%%%%%
%%%%%%%%%%%%%%%%%%%%%%%%%%%%%%%%%%%%%%%%%%%%%%%%%%%%%%%%%%%%%%%%
%%%%%%%%%%%%%%%%%%%%%%%%%%%%%%%%%%%%%%%%%%%%%%%%%%%%%%%%%%%%%%%%%

\section{Main theorem: A relationship between local cohomology and Tate homology}\label{mainsection}

In this section we prove our main theorem, Theorem~\ref{t12}, which establishes a new relationship between local cohomology of tensor products of modules and Tate homology (recall that, unless otherwise stated, all modules are assumed to be finitely generated).

\begin{thm}\label{t12} Let $M$ and $N$ be $R$-modules, $\mz$ be a specialization-closed subset of $\Spec(R)$, and let $n\geq 0$ be an integer. Assume the following conditions hold:
\begin{enumerate}[\rm(i)]
\item $\Gdim_R(M)<\infty$.
\item $\depth_{R_\fp}(M_\fp)+\depth_{R_\fp}(N_\fp)\geq\depth(R_{\fp})+n$ for all $\fp\in\mz$.
\item $\Supp_R\big(\Tor_i^R(M,N)\big) \bigcup \Supp_R\Big(\widehat{\Tor}_j^R(M,N)\Big) \subseteq\mz$ for all $i\geq 1$ and for all $j\in\ZZ$.
\end{enumerate}
Then $\hh^i_{\mz}(M\otimes_RN)\cong\widehat{\Tor}_{-i}^R(M,N)$ for all $i=0, \ldots, n-1$. Moreover, there is an injective $R$-module homomorphism $\widehat{\Tor}_{-n}^R(M,N)\hookrightarrow\hh^n_{\mz}(M\otimes_RN)$.
\end{thm}

Note that the hypothesis of part~(iii) in Theorem~\ref{t12} holds if the nonfree locus of $M$ is contained in $\mz$, that is, if $\NF(M)=\{\fp \in \Spec(R) : M_{\fp} \text{ is not free over } R_{\fp} \} \subseteq \mz$. Hence, letting $n=\dim(R)$ and $\mz=\V(\fm)$, Theorem~\ref{main-thm-intro} immediately follows from Theorem~\ref{t12}.

To the best of our knowledge, Theorem~\ref{t12} is new, even if the specialization-closed subset $\mz$ in the theorem is a closed subset of $\Spec(R)$. There are several consequences of Theorem~\ref{t12}, proven in the next sections, which contribute to various areas of the commutative algebra literature; see, for example, Corollaries~\ref{DDmore} and~\ref{cxmm2}, and Proposition~\ref{ct}.

Next, we state a few preparatory results, each necessary for proving Theorem~\ref{t12}. To avoid disturbing the flow of the paper, we give a proof of Theorem~\ref{t12} in this section, and defer the proofs of these auxiliary results to Section~\ref{PAR}.

\begin{lem} \label{lem1} Let $M$ and $N$ be $R$-modules, and let $\mz$ be a specialization-closed subset of $\Spec(R)$. Set $\mz'=\mz\cap\Supp_R(M)$. Then $\hh^i_{\mz}(M\otimes_RN)\cong\hh^i_{\mz'}(M\otimes_RN)$ for all $i\geq 0$.
\end{lem}

\begin{lem}\label{pr} Let $M$ and $N$ be $R$-modules and let $\X\subseteq \Spec(R)$ be a subset. Assume $M_\fp$ is totally reflexive over $R_\fp$ for all $\fp\in\X$. Then  the following conditions are equivalent:
\begin{enumerate}[\rm(i)]
\item $\widehat{\Tor}_i^{R_\fq}(M_\fq,N_\fq)=0$ for all $i\in\ZZ$ and all $\fq\in\X$.
\item $\Supp_R\left(\Ext^i_R(\Tr_R M,N)\right)\bigcup\Supp_R\left(\Tor_i^R(M,N)\right)\subseteq \Spec(R)-\X $ for all $i\geq 1$. 
\end{enumerate}
\end{lem}

Lemma~\ref{grz} plays an important role in this paper and can be viewed as a generalization of~\cite[4.1]{TTY}. The conclusion of Lemma~\ref{grz}, for the special case where $\mz$ is a closed subset of $\Spec(R)$, was proved by Grothendieck; see~\cite[III.2.9]{sga2}.  

\begin{lem}\label{grz} Let $M$ be an $R$-module and let $\mz$ be a specialization-closed subset of $\Spec(R)$. Then,
$$\gr_R(\mz,M)=\inf\big\{i\geq 0 : \hh_{\mz}^{i}(M)\neq0\}=\inf \left\{ \depth_{R_{\fp}}(M_{\fp}): \fp\in\mz \right\}.$$
\end{lem}

\begin{prop}\label{l2} Let $M$ and $N$ be $R$-modules, where $N$ is not necessarily finitely generated, and $\mz$ be a specialization-closed subset of $\Spec(R)$. If $n$ is an integer such that $0\leq n\leq\gr_R(\mz,N)$ and $\Supp_R\left(\Ext^i_R(\Tr_R M,N)\right)\subseteq\mz$ for all $i=1, \ldots, n+1$, then the following hold:
\begin{enumerate}[\rm(i)]
\item $\hh^{i}_{\mz}(M\otimes_RN)\cong\Ext^{i+1}_R(\Tr_R M,N)$ for all $i=0, \ldots, n-1$.
\item There is an injection $\Ext^{n+1}_R(\Tr_R M,N)\hookrightarrow\hh^{n}_{\mz}(M\otimes_RN)$. 
\end{enumerate}
\end{prop}

Next is the proof of the main theorem of this paper:

\begin{proof}[Proof of Theorem~\ref{t12}] We may replace $\mz$ with $\mz\cap\Supp_R(M)$ since this does not affect the hypotheses and the conclusion of the theorem; see Lemma~\ref{lem1}. Therefore, we assume $\mz\subseteq\Supp_R(M)$. 

As $\Gdim_R(M)<\infty$, we know there is an exact sequence of $R$-modules
\begin{equation}\label{t121}\tag{\ref*{t12}.1}
0\rightarrow M\rightarrow X\rightarrow G\rightarrow0,
\end{equation}
where $G$ is totally reflexive and $\pd_R(X)<\infty$; see~\ref{HD}(iii). The exact sequence~\eqref{t121}, in view of~\ref{Tate2}(iv), yields the following long exact sequence of stable homology modules:
\begin{equation}\label{t122}\tag{\ref*{t12}.2}
\cdots \rightarrow \widehat{\Tor}_{j+1}^R(X,N) \rightarrow\widehat{\Tor}_{j+1}^R(G,N)\rightarrow\widehat{\Tor}_{j}^R(M,N)\rightarrow\widehat{\Tor}_{j}
^R(X,N)\rightarrow
\widehat{\Tor}_{j}^R(G,N)\rightarrow\cdots
\end{equation}
Note, by~\ref{Tate hom}(iii), we have that $\widehat{\Tor}_{j}^R(X,N)=0$ for all $j\in\ZZ$. Hence,~\eqref{t122} gives the isomorphism:
\begin{equation}\label{t123}\tag{\ref*{t12}.3}
\widehat{\Tor}_{j+1}^R(G,N)\cong\widehat{\Tor}_{j}^R(M,N) \text{ for all } j\in\ZZ.
\end{equation}
As $G$ is totally reflexive, we have that $\Tor_i^R(G,N) \cong \widehat{\Tor}_i^R(G,N)$ for all $i\geq 1$; see~\ref{Tate hom}(i). Therefore,~\eqref{t123} implies that:
\begin{equation}\label{t124}\tag{\ref*{t12}.4}
\Supp_R\big(\Tor_i^R(G,N)\big)=\Supp_R\big(\widehat{\Tor}_i^R(G,N)\big)=\Supp_R\big(\widehat{\Tor}_{i-1}^R(M,N)\big)\subseteq\mz \text{ for all } i\geq 1.
\end{equation}
Applying the functor $-\otimes_RN$ to the exact sequence~\eqref{t121}, we get the following exact sequence:
\begin{equation}\label{t125}\tag{\ref*{t12}.5}
\cdots\rightarrow\Tor_i^R(M,N)\to\Tor_i^R(X,N)\to\Tor_i^R(G,N)\rightarrow\cdots.
\end{equation}
In view of hypothesis (iii), we see from~(\ref{t12}.4) and~(\ref{t12}.5) that: 
\begin{equation}\label{t126}\tag{\ref*{t12}.6}
\Supp_R\left(\Tor_i^R(X,N)\right)\subseteq\Supp_R(\Tor_i^R(M,N))\cup\Supp_R\left(\Tor_i^R(G,N)\right)\subseteq\mz \text{ for all } i\geq 1.
\end{equation}

We deduce $\depth_{R_{\fp}}(M_\fp)=\depth_{R_{\fp}}(X_\fp)$ for all $\fp\in\mz\subseteq\Supp_R(M)\subseteq\Supp_R(X)$ by using the exact sequence (\ref{t12}.1). Therefore, hypothesis (ii) shows that:		
\begin{equation} \label{t127}\tag{\ref*{t12}.7}
\depth_{R_{\fp}}(X_\fp)+\depth_{R_{\fp}}(N_\fp)\geq\depth(R_{\fp})+ n \text{ for all } \fp\in\mz.
\end{equation}

Next we claim: $\Tor_i^R(X,N)=0$ for all $i\geq 1$.\\
Proof of the claim. Set $q=\sup\{j\geq 0: \Tor_{j}^R(X,N)\neq0\}$. Suppose $q\neq 0$ and pick $\fp\in\Ass_R(\Tor_q^R(X,N))$. Note that $\fp \in \mz$; see~(\ref{t12}.6). As $\depth_{R_\fp}(\Tor_q^R(X,N)_\fp)=0$, the depth formula~\ref{df} and~(\ref{t12}.7) imply:
$$0\leq n\leq \depth_{R_\fp}(X_\fp)+\depth_{R_\fp}(N_\fp)-\depth(R_\fp)=-q<0.$$
This contradiction completes the proof of the claim. 

As the claim shows $\Tor_1^R(X,N)=0$, by tensoring the exact sequence~(\ref{t12}.1) with $N$, we obtain the following exact sequence:
\begin{equation}\label{t128}\tag{\ref*{t12}.8}
0\rightarrow\Tor_1^R(G,N)\rightarrow M\otimes_RN\rightarrow X\otimes_RN\rightarrow G\otimes_RN\rightarrow0.
 \end{equation}
 Thus, by applying the functor $\Gamma_{\mz}(-)$ to the exact sequence~(\ref{t12}.8), we obtain the following injection:
 \begin{equation} \label{t129}\tag{\ref*{t12}.9}
 \hh^0_{\mz}(\Tor_1^R(G,N))\hookrightarrow\hh^0_{\mz}(M\otimes_RN).
\end{equation}

Recall that $\Supp_R\left(\Tor_1^R(G,N)\right)\subseteq\mz$; see~(\ref{t12}.6). Therefore, $\Tor_1^R(G,N)$ is torsion with respect to $\mz$; see~\ref{LC}(iv). Hence, $\Tor_{1}^R(G,N)\cong\hh^0_{\mz}(\Tor_1^R(G,N))$, and the following hold: 
\begin{equation}\label{t1210}\tag{\ref*{t12}.10}
\widehat{\Tor}_0^R(M,N)\cong\widehat{\Tor}_1^R(G,N)\cong\Tor_{1}^R(G,N)\cong\hh^0_{\mz}\left(\Tor_1^R(G,N)\right)\hookrightarrow\hh^0_{\mz}(M\otimes_RN).
\end{equation}
Here, in~(\ref{t12}.10) the first and the second isomorphisms hold by~(\ref{t12}.3) and~\ref{Tate hom}(i), respectively, and the injection is due to~(\ref{t12}.9). This completes the proof of the theorem if $n=0$. Consequently, from now on, we assume $n\geq 1$.

We know, as we established in the claim, that $\Tor_i^R(X,N)=0$ for all $i\geq 1$. As $\pd_R(X)<\infty$, the depth formula~\ref{df} gives:
\begin{equation}\tag{\ref{t12}.11}
\depth_{R_\fp}(X_\fp)+\depth_{R_\fp}(N_\fp)=\depth(R_{\fp})+\depth_{R_\fp}((X\otimes_RN)_\fp) \text{ for all } \fp\in\Supp_R(X\otimes_RN).
\end{equation}	

Let $\fq \in \mz$. Then $\fq \in \mz \subseteq  \Supp_R(M)$ so that $\fq \in \Supp_R(X)$; see~(\ref{t12}.1). If $\fq \in \Supp_R(N)$, then
$\depth_{R_\fq}((X\otimes_RN)_\fq)\geq n$ due to the (in)equalities of~(\ref{t12}.7) and~(\ref{t12}.11). Also, if $\fq \notin \Supp_R(N)$, then $\depth_{R_\fq}((X\otimes_RN)_\fq)=\infty$. These show that $\depth_{R_\fp}\left((X\otimes_RN)_\fp\right)\geq n$ for all $\fp \in \mz$. Hence, Lemma~\ref{grz} implies:
\begin{equation}\tag{\ref{t12}.12}
\hh^i_{\mz}(X\otimes_RN)=0 \text{  for  all  } i=0, \ldots, n-1.
\end{equation}

We obtain the following two exact sequences of $R$-modules from~(\ref{t12}.8):
\begin{equation}\tag{\ref{t12}.13}
0\rightarrow\Tor_1^R(G,N)\rightarrow M\otimes_RN\to Y\to0\ \ \text{ and }\ \ 0\to Y\rightarrow X\otimes_RN\rightarrow G\otimes_RN\rightarrow0.
\end{equation}
By applying the functor $\Gamma_{\mz}(-)$ to the rightmost exact sequence in~(\ref{t12}.13), in view of~(\ref{t12}.12), we get: $$\hh^{0}_{\mz}(Y)\hookrightarrow \hh^{0}_{\mz}(X\otimes_RN)=0,\;\; \hh^i_{\mz}(Y)\cong\hh^{i-1}_{\mz}(G\otimes_RN) \text{  for  all  } i=1, \ldots, n-1, \text{ and } \hh^{n-1}_{\mz}(G\otimes_RN)\hookrightarrow \hh^n_{\mz}(Y).$$

Recall that $\Tor_1^R(G,N)$ is torsion with respect to $\mz$. Hence $\hh^i_{\mz}(\Tor_1^R(G,N))=0$ for all $i\geq 1$; see~\ref{LC}(v). So, by applying the functor $\Gamma_{\mz}(-)$ to the leftmost exact sequence in~(\ref{t12}.13), we obtain:
$$ \hh^i_{\mz}(M\otimes_RN)\cong\hh^i_{\mz}(Y) \text{ for all }i\geq 1$$
Therefore, by the foregoing discussion, we conclude the following; see, also~(\ref{t12}.10):
\begin{align}
\notag{} & \hh^0_{\mz}(M\otimes_RN)\cong\hh^0_{\mz}(\Tor_1^R(G,N))\cong\Tor_1^R(G,N)\cong\widehat{\Tor}_1^R(G,N)\cong\widehat{\Tor}_0^R(M,N), \\ & 
\tag{\ref{t12}.14} \hh^i_{\mz}(M\otimes_RN)\cong\hh^i_{\mz}(Y)\cong\hh^{i-1}_{\mz}(G\otimes_RN) \text{  for  all  } i=1, \ldots, n-1 \text{, and } &\\
\notag{}&\hh^{n-1}_{\mz}(G\otimes_RN)\hookrightarrow\hh^n_{\mz}(Y)\cong\hh^n_{\mz}(M\otimes_RN).
\end{align}

Recall that $\Supp_R\big(\widehat{\Tor}^R_i(G,N)\big)\subseteq\mz$ for all $i\in\ZZ$; see~(\ref{t12}.4). Hence, if $\X=\Spec(R)-\mz$, then $\widehat{\Tor}^R_i(G,N)_{\fp}=0$ for all $\fp\in\X$ and for all $i\in\ZZ$. Now, as $G$ is totally reflexive, Lemma~\ref{pr} shows:
\begin{equation}\tag{\ref{t12}.15}
\Supp_R\big(\Ext^i_R(\Tr_R G,N)\big)\subseteq\mz \text{ for all } i\geq 1. 
\end{equation}

Recall that, we assume $\Gdim_R(M)<\infty$ and $\depth_{R_\fp}(M_\fp)+\depth_{R_\fp}(N_\fp) \geq \depth(R_{\fp})+n$ for all $\fp \in \mz$, as hypotheses of the theorem. Hence, by using the Auslander-Bridger formula, we have
$$\depth_{R_\fp}(N_\fp)\geq \depth(R_{\fp})-\depth_{R_\fp}(M_\fp)+n = \Gdim_{R_\fp}(M_\fp)+n\geq n \text{ for all } \fp\in\mz\subseteq\Supp_R(M).$$
Therefore, Lemma~\ref{grz} implies that:
\begin{equation}\tag{\ref{t12}.16}
\gr_R(\mz,N)\geq n. 
\end{equation}
We use of~(\ref{t12}.15) and~(\ref{t12}.16), and conclude by Proposition~\ref{l2} that $\hh^{j-1}_{\mz}(G\otimes_RN)\cong\Ext^j_R(\Tr G,N)$ for all $j=1, \ldots, n$. This yields the following isomorphisms for all $j=1, \ldots, n$:
\[\begin{array}{rl}\tag{\ref{t12}.17}
\hh^{j-1}_{\mz}(G\otimes_RN)&\cong\Ext^j_R(\Tr_R G,N)\\
&\cong\widehat{\Ext}^j_R(\Tr_R G,N)\\
&\cong\widehat{\Ext}^{j-2}_R(G^{\ast},N)\\
&\cong\widehat{\Tor}_{-j+1}^R(G,N).
\end{array}\]
Here, in~(\ref{t12}.17), the second isomorphism holds since $\Tr_R G$ is totally reflexive; see~\ref{HD}(i). The third isomorphism is due to the fact that $G^{\ast} \cong \Omega^2_R \Tr_R G$; see~\ref{a1} and~\ref{Tate2}(v). Also, the last isomorphism follows from~\ref{Tate2}(ii). Now we use~(\ref{t12}.3),~(\ref{t12}.14), and~(\ref{t12}.17), and obtain:
$$\hh^{i+1}_{\mz}(M\otimes_RN)\cong \hh^{i}_{\mz}(G\otimes_RN) \cong \widehat{\Tor}_{-i}^R(G,N)\cong \widehat{\Tor}_{-i-1}^R(M,N) \text{ for all } i=0, \ldots, n-2, \text{ and } $$
$$\widehat{\Tor}_{-n}^R(M,N) \cong \widehat{\Tor}_{-n+1}^R(G,N) \cong \hh^{n-1}_{\mz}(G\otimes_RN)\hookrightarrow \hh^{n}_{\mz}(M\otimes_RN).$$
This completes the proof of the theorem.
\end{proof}

%%%%%%%%%%%%%%%%%%%%%%%%%%%%%%%%%%%%%%%%%%%%%%%%%%%%%%%%%%%%%%%%%
%%%%%%%%%%%%%%%%%%% SECTION  %%%%%%%%%%%%%%%%%%%%%%%%%%%%%%%%%%%%%%%%
%%%%%%%%%%%%%%%%%%% SECTION  %%%%%%%%%%%%%%%%%%%%%%%%%%%%%%%%%%%%%%%%
%%%%%%%%%%%%%%%%%%% SECTION  %%%%%%%%%%%%%%%%%%%%%%%%%%%%%%%%%%%%%%%%
%%%%%%%%%%%%%%%%%%%%%%%%%%%%%%%%%%%%%%%%%%%%%%%%%%%%%%%%%%%%%%%%%

\section{On the depth of tensor products of modules}\label{section4}

This section, which is divided into two subsections, is primarily devoted to applications of Theorem~\ref{t12}. In the first subsection, we use Theorem~\ref{t12} to prove a variation of a result of Auslander, previously studied by Huneke-Wiegand, Dao, and others; see~\ref{DD} and Corollary~\ref{DDmore}. In the second subsection, we present further applications of Theorem~\ref{t12} and investigate the depth of tensor products of modules. For example, Corollary~\ref{ctc-main} provides a proof of Proposition~\ref{ctc-intro}, as advertised in the introduction; see also Corollary~\ref{corcombined} for a related application.

%%%%%%%%%%%%%%%%%SUBSECTION%%%%%%%%%%%%%%%%%%%%%%%%%%%%%%%%

\subsection*{On a result of Auslander} In this subsection, we consider results that stem from a fundamental theorem of Auslander~\cite{Au}. To facilitate the discussion, we recall a few definitions and briefly review the relevant literature.

\begin{chunk}{\bf Serre's condition}. \label{SL} An $R$-module $M$ is said to satisfy \emph{Serre's condition} $(S_n)$ for some integer $n\geq 0$ provided that $\depth_{R_{\fp}}(M_\fp)\geq\min\{n,\Ht_R(\fp)\}$ for all $\fp\in\Spec(R)$ (recall that, by convention, $\depth_R(0)=\infty$.) The following properties hold; see~\cite[3.6]{EG}.
\begin{enumerate}[\rm(i)]
\item If $R$ is unmixed (that is, $R$ has no embedded primes), $M$ satisfies $(S_1)$ if and only if $M$ is torsion-free. 
\item If $R$ is Gorenstein, then $M$ satisfies $(S_2)$ if and only if $M$ is reflexive.
\item If $R$ is Gorenstein and $n\geq 0$, then $M$ satisfies $(S_n)$ if and only if $M$ is an $n$th syzygy module.
\end{enumerate}
\end{chunk}

It is not difficult to see that the torsion-freeness of the tensor product does not, in general, imply the vanishing of $\Tor$ over local rings that are not regular. On the other hand, Auslander~\cite{Au} (for unramified regular rings) and Lichtenbaum~\cite{Li} (for ramified regular rings) proved the following:

\begin{chunk} [{\cite[3.1]{Au}; see also~\cite[Cor. 2]{Li}}] \label{AuLi} Let $R$ be a regular local ring and let $M$ and $N$ be $R$-modules. If $M\otimes_RN$ satisfies $(S_1)$, that is, torsion-free, then $\Tor_i^R(M,N)=0$ for all $i\geq 1$.
\end{chunk}

Regular rings are complete intersections of codimension zero. In view of Theorem~\ref{AuLi}, it is natural to seek a similar vanishing result over complete intersections of codimension one, that is, over hypersurfaces. In general, over a hypersurface ring, a torsion-free tensor product---even a reflexive one---does not imply the vanishing of $\Tor$. Huneke and Wiegand~\cite{HW1} established an analogue of Theorem~\ref{AuLi} in the setting of hypersurfaces. We recall their result after a few definitions.
 
Given an integer $n \geq 0$, define $\X^{n}(R) = \{ \fp \in \Spec(R) : \depth(R_{\fp}) \leq n \}$. An $R$-module $M$ is said to have \emph{rank on} $\X^n(R)$ if there exists an integer $s \geq 0$ such that $M_{\fp} \cong R_{\fp}^{\oplus s}$ for all $\fp \in \X^n(R)$. In particular, $M$ is said to have \emph{rank} if it has rank on $\X^0(R)$.

The following result slightly extends a theorem of Huneke and Wiegand, who proved the case $c=1$ in~\cite[2.7]{HW1}.

\begin{chunk}[{\cite[4.20]{Ce}}] \label{Ce-HW1} Let $R$ be a local complete intersection of codimension $c$, and let $M$ and $N$ be $R$-modules, at least one of which has rank. Assume $M \otimes_{R} N$ satisfies $(S_2)$, that is, reflexive. If $c \geq 2$, assume further:
\begin{enumerate}[\rm(i)]
\item Both $M$ and $N$ satisfy $(S_1)$.
\item $\Tor^{R}_{1}(M,N) = \Tor^{R}_{2}(M,N) = \cdots = \Tor^{R}_{c-1}(M,N) = 0$.
\end{enumerate}
Then $\Tor^{R}_{i}(M,N) = 0$ for all $i \geq 1$.
\end{chunk}

In view of the aforementioned results of Auslander, Lichtenbaum, and Huneke and Wiegand, one naturally raises the following question, which still remains open.

\begin{ques} \label{que-HWC} Let $R$ be a local complete intersection of codimension $c$, and let $M$ and $N$ be $R$-modules, at least one of which has rank. If $M \otimes_{R} N$ satisfies $(S_{c+1})$, then must $\Tor^{R}_{i}(M,N) = 0$ for all $i \geq 1$? What if $R$ is a domain?
\end{ques} 

Huneke-Wiegand~\cite{HW2} established several results analogous to~\ref{Ce-HW1}, primarily for hypersurface rings. Huneke-Jorgensen-Wiegand~\cite{HJW} continued the work initiated in~\cite{HW1, HW2}, examining the relationship between the vanishing of $\Tor$ and the depth of tensor products of modules over complete intersection rings. Their techniques were effective for rings of codimension two and three but failed in higher codimensions. We present one of their results, which extends~\cite[2.7]{HW1} to complete intersections of codimension two.

\begin{chunk}[{\cite[2.4]{HJW}}] Let $R$ be a local complete intersection of codimension $2$, and let $M$ and $N$ be $R$-modules. Assume the following conditions hold:
\begin{enumerate}[\rm(i)]
\item $M$ has rank on $\X^1(R)$ and $N$ has rank.
\item $M$ and $N$ satisfy $(S_{2})$.
\end{enumerate}
If $M\otimes_RN$ satisfies $(S_{3})$, then $\Tor^{R}_{i}(M,N) = 0$ for all $i \geq 1$.
\end{chunk}

Our results are motivated by the following, which gives a partial answer to Question~\ref{que-HWC} by replacing the rank assumption with a local freeness condition on $M$ at primes of height at most $c$, and by requiring certain depth conditions on $M$ and $N$.

\begin{chunk} [{\cite[7.7]{Da2} and \cite[3.14]{CIPW}}] \label{DD} Let $R$ be a local complete intersection ring, and let $M$ and $N$ be $R$-modules. Assume $c$ is an integer such that $c \geq \cod(R)$. Assume further:
\begin{enumerate}[\rm(i)]
\item $M$ and $N$ satisfy $(S_{c})$.
\item $M\otimes_RN$ satisfies $(S_{c+1})$.
\item $M_{\fp}$ is free over $R_{\fp}$ for all $\fp \in \X^{c}(R)$. 
\end{enumerate}
Then $\Tor_{i}^R(M,N)=0$ for all $i\geq 1$.
\end{chunk}

The result stated as~\ref{DD} was initially proved by Dao~\cite{Da2} for complete intersection rings that are quotients of unramified regular local rings. The requirement that the ambient regular ring be unramified was later removed by Celikbas-Iyengar-Piepmeyer-Wiegand~\cite{CIPW}. Our goal is to establish Proposition~\ref{t11}, a general corollary of Theorem~\ref{t12}, and then use it to prove Corollary~\ref{DDmore}, which is a variation of~\ref{DD}. We now introduce additional terminology needed for our arguments.

\begin{chunk} \textbf{Cosyzygies}.\label{cs} Let $M$ be an $R$-module and let $M\overset{\partial_{-1}}{\longrightarrow} P_{-1}$ be a left projective approximation of $M$. For example, if $P$ is a projective $R$-module and $\phi:P\twoheadrightarrow M^*$ is a surjective $R$-module map, then $e_M\circ\phi^*$ is a left projective approximation of $M$, where $e_M$ is the natural map $M\to M^{**}$. 

The cokernel of $\partial_{-1}$ is called the first \emph{cosyzygy} of $M$, and it is denoted by $\Omega_R^{-1}M$. We define define the $n$th cosyzygy 
of $M$ inductively: $\Omega_R^{-n}M= \Omega_R^{-1} \Omega_R^{-(n-1)}M $ for $n\geq 1$; see, for example,~\cite[2.3]{ST2}. 
\end{chunk}

\begin{chunk}\label{cx}\textbf{Complexity}. Let $R$ be a local ring and let $M$ and $N$ be $R$-modules. The \emph{complexity} of the pair $(M,N)$, defined by Avramov and Buchweitz~\cite{AvBu}, is
$$\cx_R(M,N)=\inf\{ r\in\NN \mid \exists\; a \in \RR \text{ such that } \mu_R\big(\Ext^n_R(M,N)\big)\leq a n^{r-1} \text{ for all } n\gg 0\},$$
where  $\mu_R(\Ext^n_R(M,N))$ is the minimal number of generators of $\Ext^n_R(M,N)$. Accordingly, the complexity $\cx_R(M)$ of $M$, initially introduced in local algebra by Avramov~\cite{Av1}, equals $\cx_R(M,k)$; see also~\cite{AGP}. 

It follows from definition that $\cx_R(M,N)=0$ if and only if $\Ext^i_R(M,N)=0$ for all $i\gg0$; $M$ has bounded Betti numbers if and only if $\cx_R(M)\leq1$; and $\pd_R(M)<\infty$ if and only if $\cx_R(M)=0$. 

If $R$ is a complete intersection, then by~\cite[5.7]{AvBu} we have:
\begin{equation}\tag{\ref*{cx}.1} \label{cx.1}
0\leq \cx_R(M,N)=\cx_R(N,M)\leq\min\{\cx_R(M),\cx_R(N)\}\leq\cod(R)=\edim(R)-\dim(R).
\end{equation}
\end{chunk}

The following facts are used in the proof of Proposition~\ref{t11} and later in the sequel.

\begin{chunk}\label{fp} Let $R$ be a local complete intersection ring and let $M$ and $N$ be $R$-modules.
\begin{enumerate}[\rm(i)]
\item It follows that $\Ext^i_R(M,N)=0$ for all $i\gg0$ if and only if  $\Tor_i^R(M,N)=0$ for all  $i\gg 0$  if and only if $\widehat{\Tor}_j^R(M,N)=0$ for all $j\in\ZZ$; see~\cite[4.9 and 6.1]{AvBu}. Thus, if $\mz$ is a specialization-closed subset of $\Spec(R)$ and $\Supp_R\big(\Tor_{i}^R(M,N)\big)\subseteq\mz$ for all $i\gg 0$, then $ \Supp_R\Big(\widehat{\Tor}_j^R(M,N)\Big) \subseteq\mz$ for all  $j\in\ZZ$.
\item Assume $\cx_R(M,N)\leq r-1$ for some integer $r\geq 1$. If $r\geq 2$, assume further $M$ or $N$ is maximal Cohen-Macaulay and $\Tor_i^R(M,N)=0$ for all $i=1, \ldots, r$. Then $\widehat{\Tor}_{j}^R(M,N)=0$ for all $j\in \ZZ$; see~\cite[3.5]{CST}.
\item If $\Tor_{i}^R(M,N)=0$ for all $i\gg0$, then~\cite[2.2]{J2} shows that: 
\begin{equation}\notag{}
\sup\{i\geq 0 \mid\Tor_{i}^R(M,N)\neq0\}=\sup\big\{\depth(R_\fp)-\depth_{R_{\fp}}(M_\fp)-\depth_{R_{\fp}}(N_\fp): \fp\in\Spec(R)\big \}.	
\end{equation}
\end{enumerate}
\end{chunk}

Next is the main result of this section:

\begin{prop}\label{t11} Let $R$ be a local complete intersection ring, $\mz$ be a specialization-closed subset of $\Spec(R)$, and let $M$ and $N$ be $R$-modules. Assume there are integers $c$ and $n$ such that $n\geq c-1 \geq \cx_R(M,N)$. Assume further the following hold:
\begin{enumerate}[\rm(i)]
\item $\depth_{R_\fp}(M_\fp)+\depth_{R_\fp}(N_\fp)\geq \depth(R_\fp)+n$ for all $\fp\in\mz$.
\item $\hh^i_{\mz}(M\otimes_RN)=0$ for all $i=n-c+1, \ldots, n$.
\item $\Supp_R\big(\Tor_{i}^R(M,N)\big)\subseteq\mz$ for all $i\geq 1$ (for example, $\NF(M)\cap\NF(N)\subseteq\mz$).
\end{enumerate}
Then $\gr_R(\mz,M\otimes_RN)\geq n+1$, $\Tor_i^R(M,N)=0$ for all $i\geq 1$, and $\widehat{\Tor}_j^R(M,N)=0$ for all $j\in \ZZ$.
\end{prop}

\begin{proof} It follows that $\Supp_R\big(\Tor_i^R(M,N)\big) \bigcup \Supp_R\Big(\widehat{\Tor}_j^R(M,N)\Big) \subseteq\mz$ for all $i\geq 1$ and for all $j\in\ZZ$ since we assume $\Supp_R\big(\Tor_{i}^R(M,N)\big)\subseteq\mz$ for all $i\geq 1$; see~\ref{fp}(i). Thus, Theorem~\ref{t12} yields:
\begin{equation}\tag{\ref{t11}.1}
\hh^i_{\mz}(M\otimes_RN)\cong\widehat{\Tor}_{-i}^R(M,N) \text{ for all } i=0, \ldots, n-1 \text{ and } \widehat{\Tor}_{-n}^R(M,N)\hookrightarrow\hh^n_{\mz}(M\otimes_RN).
\end{equation}
Hence, since $n-c+1 \geq 0$ and $\hh^n_{\mz}(M\otimes_RN)=0$, we have by (\ref{t11}.1) that:
\begin{equation}\tag{\ref{t11}.2}
\widehat{\Tor}_{-i}^R(M,N)=0 \text{ for all } i=n-c+1, \ldots, n.
\end{equation}

Consider a maximal Cohen-Macaulay approximation of $M$, that is, an exact sequence of $R$-modules of the form $0\to L\to X\to M\to0$, where $X$ is totally reflexive and $\pd_R(L)<\infty$. As $\pd_R(L)<\infty$, we have by~\ref{Tate hom}(iii) that $\widehat{\Tor}_i^R(L,N)=0$ for all $i\in\ZZ$. So, in view of~\ref{Tate hom}(v), the sequence $0\to L\to X\to M \to 0$ induces the following isomorphisms for all $i \in \ZZ$:
\begin{equation}\tag{\ref{t11}.3}
\widehat{\Tor}_i^R(X,N) \cong \widehat{\Tor}_i^R(M,N).
\end{equation}

We set $Y=\Omega_R^{-(n+1)}X$; see~\ref{cs}. Then $Y$ is totally reflexive, $X \cong \Omega_R^{n+1}Y$, and the following hold:
\begin{equation}\tag{\ref{t11}.4}
\Tor_i^R(Y,N)\cong\widehat{\Tor}_i^R(Y,N)\cong\widehat{\Tor}_{i-n-1}^R(X,N)\cong\widehat{\Tor}_{i-n-1}^R(M,N) \text{ for all } i\geq 1.
\end{equation}
Here, in (\ref{t11}.4), the first isomorphism is true by~\ref{Tate2}(i) because $Y$ is totally reflexive; the second one holds due to~\ref{Tate2}(v); and the third isomorphism follows from (\ref{t11}.3).

It follows from \ref{fp}(ii) and (\ref{t11}.1) that $\widehat{\Tor}_i^R(Y,N)=0$ for all $i\in\ZZ$. This implies, since $X \cong \Omega_R^{n+1}Y$, that $\widehat{\Tor}_i^R(X,N)=0$ for all $i\in\ZZ$; see \ref{Tate2}(iv). Hence, we conclude from (\ref{t11}.3) that 
\begin{equation}\tag{\ref{t11}.5}
\widehat{\Tor}_i^R(M,N)=0 \text{ for all } i\in \ZZ.
\end{equation}
Now (\ref{t11}.1)  and (\ref{t11}.5) yield $\hh^i_{\mz}(M\otimes_RN)\cong\widehat{\Tor}_{-i}^R(M,N)=0$ for all $i=0, \ldots, n-1$. As we already assume $\hh^n_{\mz}(M\otimes_RN)=0$, we see that $\hh^i_{\mz}(M\otimes_RN)$ for all $i=0, \ldots, n$. Therefore, it follows that $\gr_R(\mz,M\otimes_RN)\geq n+1$; see \ref{gsc}. Note, by~\ref{Tate hom}(i) and (\ref{t11}.5), we have that $\Tor_{i}^R(M,N)=0$ for all $i\gg0$.  

Let $\fp \in \Spec(R)$. First consider the case where $\fp \notin \mz$. Then $\fp \notin \Supp_R\big(\Tor_{i}^R(M,N)\big)$ for all $i\geq 1$. Therefore, it follows from \ref{df} that $\depth(R_{\fp})-\depth_{R_{\fp}}(M_{\fp})-\depth_{R_{\fp}}(N_{\fp})=-\depth_{R_{\fp}}(M\otimes_RN)_{\fp}\leq 0$. Next assume $\fp \in \mz$. Then, by our assumption, we have $\depth(R_{\fp})-\depth_{R_{\fp}}(M_{\fp})-\depth_{R_{\fp}}(N_{\fp})\leq -n\leq 0$. 

Hence, we observed that $\depth(R_{\fq})-\depth_{R_{\fq}}(M_{\fq})-\depth_{R_{\fq}}(N_{\fq})\leq 0$ for all $\fq \in \Spec(R)$. So, we deduce from~\ref{fp}(iii) that $\Tor_i^R(M,N)=0$ for all $i\geq 1$.
\end{proof}	
 
The following is an immediate consequence of Proposition~\ref{t11}, for the case where $\mz=\{\fm\}$, and the depth formula~\ref{df}.

\begin{cor} \label{cxmm1} Let $R$ be a local complete intersection ring and let $M$ and $N$ be $R$-modules. Assume there are integers $c$ and $n$ such that $n\geq c-1 \geq \cx_R(M,N)$. Assume further the following hold:
\begin{enumerate}[\rm(i)]
\item $\depth_{R}(M)+\depth_{R}(N)-\depth(R) \geq n$.
\item $\hh^i_{\fm}(M\otimes_RN)=0$ for all $i=n-c+1, \ldots, n$.
\item $\NF(M)\cap\NF(N)\subseteq\{\fm\}$.
\end{enumerate}
Then $\Tor_i^R(M,N)=0$ for all $i\geq 1$, and $\depth_R(M\otimes_RN)=\depth_{R}(M)+\depth_{R}(N)-\depth(R)\geq n+1$.
\end{cor}

We are now ready to prove a variation of~\ref{DD}, as previously mentioned. Recall that, given an integer $n\geq 0$, we set $\X^{n}(R)=\{\fp \in \Spec(R): \depth(R_{\fp})\leq n \}$.

\begin{cor} \label{DDmore} Let $R$ be a local complete intersection ring, and let $M$ and $N$ be $R$-modules. Assume $c$ is an integer such that $c \geq  \cx_R(M,N)+1$. Assume further:
\begin{enumerate}[\rm(i)]
\item $\depth_{R_\fp}(M_\fp)+\depth_{R_\fp}(N_\fp)\geq \depth(R_\fp)+(c-1)$ for all $\fp \in \Spec(R)$ with $\Ht_R(\fp)\geq c-1$.
\item $M\otimes_RN$ satisfies $(S_{c})$.
\item $M_{\fp}$ is free over $R_{\fp}$ for all $\fp \in \X^{c-1}(R)$. 
\end{enumerate}
Then $\Tor_{i}^R(M,N)=0$ for all $i\geq 1$.
\end{cor}

Before proving Corollary~\ref{DDmore}, we highlight a special case in which condition~(i) holds:

\begin{rmk} If $R$ is a local complete intersection ring, $M$ is a maximal Cohen-Macaulay $R$-module, and $N$ is an $R$-module which satisfies $(S_{c-1})$, then $\depth_{R_\fp}(M_\fp)+\depth_{R_\fp}(N_\fp)\geq \depth(R_\fp)+(c-1)$ for all $\fp \in \Spec(R)$ with $\Ht_R(\fp)\geq c-1$.
\end{rmk}

\begin{proof}[Proof of Corollary~\ref{DDmore}] Set $\mz=\{\fp \in \Spec(R) : \Ht_R(\fp)\geq c\}$. Then $\mz$ is a specialization-closed subset of $\Spec(R)$; see~\ref{lsc}(iii). 

Let $\fq \in \Supp_R\big(\Tor_j^R(M,N)\big)$ for some $j\geq 1$. Then $\Tor_j^R(M,N)_{\fq}\neq 0$ so that $M_{\fq}$ is not free over $R_{\fq}$. Since we assume that $M_{\fp}$ is free over $R_{\fp}$ for all $\fp \in \X^{c-1}(R)$, it follows that $\fq \notin  \X^{c-1}(R)$. Thus, $\Ht_R(\fq)\geq c$ and so $\fq \in \mz$. This proves that $\Supp_R\big(\Tor_i^R(M,N)\big)\subseteq\mz$ for all $i\geq 1$.

Note that, by Lemma~\ref{grz}, we have:
$$\gr_R(\mz,M\otimes_RN)=\inf\big\{i\geq 0 : \hh_{\mz}^{i}(M\otimes_RN)\neq0\}=\inf\{\depth_{R_{\fp}}\left(M\otimes_RN\right)_{\fp}: \fp\in\mz\big\}.$$ 
Also, if $\fp \in \mz$, then $\Ht_R(\fp)\geq c$ and so $\depth_{R_{\fp}}(M\otimes_RN)_{\fp}\geq \min\{c, \Ht_R(\fp)\}=c$ since $M\otimes_RN$ satisfies $(S_{c})$.
Thus, $\gr_R(\mz,M\otimes_RN)\geq c$ and hence $\hh^i_{\mz}(M\otimes_RN)=0$ for all $i=0, \ldots, c-1$. Therefore, as $c-1\geq  \cx_R(M,N)$, we use Proposition~\ref{t11} for the case where $n=c-1$, and conclude that $\Tor_{i}^R(M,N)=0$ for all $i\geq 1$.
\end{proof}

Note that if the dimension of the ring in Theorem~\ref{DD} equals its codimension, then the theorem is vacuous, as the module $M$ in question would be free by hypothesis~(iii). On the other hand, Corollary~\ref{DDmore} provides meaningful information for modules $M$ and $N$ with $\cx_R(M,N)<\cod(R)$, even when $\cod(R) = \dim(R)$; see also~\cite[3.4]{Ce} for a related result concerning~\ref{DD}. We illustrate the existence of such modules with the following example:

\begin{eg} [{\cite[4.14]{CST}}] \label{excx} Let $R=k[\![x,y,z]\!]/(x^2, y^2, z^2)$, $M=R/(x,y)$, and let $N$ be the ideal of $R$ generated by $yz$. Then $R$ is a complete intersection and $\cx_R(M,N)=1<\min\{\cx_R(M),\cx_R(N)\}=2<\cod(R)=3$. 
\end{eg}

Since complexity is preserved under faithfully flat extensions, one can use Example~\ref{excx} to construct similar examples of modules over complete intersections of higher codimension and dimension; see~\cite[4.2.5(2)]{Av2}.

Next, we present a consequence of Theorem~\ref{t12}:

\begin{cor} \label{cxmm2} Let $R$ be a local complete intersection ring which has an isolated singularity, and let $M$ and $N$ be $R$-modules. Assume $M$ is maximal Cohen-Macaulay and $\cx_R(M,N)\leq \depth_R(N)$. Then we have:
\begin{enumerate}[\rm(i)]
\item $\depth_R(M\otimes_RN)\leq\depth_R(N)$.
\item If $M\otimes_RN$ is maximal Cohen-Macaulay, then so is $N$.
\end{enumerate}
\end{cor}

\begin{proof} It is enough to prove part (i). For that we set $c=\cx_R(M,N)+1$ and $n=\depth_R(N)$. 

Suppose $\hh^i_{\fm}(M\otimes_RN)=0$ for all $i=n-c+1, \ldots, n$. Then, since $n\geq c-1=\cx_R(M,N)$, we see from Corollary \ref{cxmm1} that $n=\depth_{R}(N)=\depth_{R}(M)+\depth_{R}(N)-\depth(R)\geq n+1$. So $\hh^i_{\fm}(M\otimes_RN)\neq0$ for some integer $i$ with $n-c+1 \leq i \leq n$. This implies that $\depth_R(M\otimes_RN)\leq\depth_R(N)$; see \ref{gsc}.
\end{proof}

We can easily see that the assumption $\cx_R(M,N)\leq \depth_R(N)$ is necessary for Corollary~\ref{cxmm2}.

\begin{eg} [{\cite[page 164]{HW2}}] Let $R=k[\![x,y]\!]/(xy)$, $M=R/(x)$, and let $N=R/(x^2)$. Then $R$ is a hypersurface with an isolated singularity. Also, $M\otimes_RN \cong M$ so that $M$ and $M\otimes_RN$ are both maximal Cohen-Macaulay. On the other hand, $N$ is not maximal Cohen-Macaulay and $\cx_R(M,N)=1>\depth_R(N)=0$.
\end{eg}

%%%%%%%%%%%%%%%%%SUBSECTION%%%%%%%%%%%%%%%%%%%%%%%%%%%%%%%%

\subsection*{Further applications of Theorem~\ref{t12}} In this subsection, we provide further applications of Theorem~\ref{t12}. We discuss torsion and determine new ways to compute depths of tensor products of modules. In particular, we prove  Proposition~\ref{ctc-intro}; see also Proposition~\ref{ct} and Corollary~\ref{corcombined}. %We start by recalling the following result from \cite{CST}:

\begin{chunk} \label{cst} Let \( R \) be a local complete intersection ring, and let $M$ and $N$ be nonzero $R$-modules such that $\len_R\left(\Tor_i^R(M, N)\right)<\infty$ for all $i \geq 1$. Set $b=\depth_R(M) + \depth_R(N) - \depth(R)$.
Then we have:
\begin{enumerate}[\rm(i)]
\item (\cite[3.1]{CST}) If $\depth_R(M \otimes_R N) \geq \cx_R(M, N)+1$, then the following conditions are equivalent:\smallskip
\begin{equation*}
\begin{aligned}
\cx_R(M, N) \leq b &\Longleftrightarrow  \depth_R(M \otimes_R N) \leq b \\
&\Longleftrightarrow \depth_R(M \otimes_R N) =b\\
&\Longleftrightarrow \Tor_i^R(M, N) = 0 \text{ for all } i \geq 1.
\end{aligned}
\end{equation*}
\item (\cite[4.1]{CST}) At least one of the following inequalities is satisfied:
\[
\depth_R(M \otimes_R N) \geq b \quad \text{or} \quad \depth_R(M \otimes_R N) \leq \cx_R(M, N).
\]
\end{enumerate}
\end{chunk}

\begin{rmk} \label{boundedbounded} Let $R$ be a local complete intersection ring, and let $M$ and $N$ be nonzero $R$-modules. Assume $\len_R\left(\Tor_i^R(M, N)\right)<\infty$ for all $i \geq 1$. The implication $\cx_R(M, N) \leq b \Longrightarrow \depth_R(M \otimes_R N) \leq b$ in~\ref{cst}(i) holds even without assuming $\depth_R(M \otimes_R N) \geq \cx_R(M, N)+1$. As we discuss next, this can be seen either by using Corollary~\ref{cxmm1} or applying~\ref{cst}(i).
\begin{enumerate}[\rm(i)]
\item Suppose $\cx_R(M,N) \leq b<\depth_R(M\otimes_RN)$, and seek a contradiction. It follows from our assumption that $\depth_R(M\otimes_RN)\geq \cx_R(M,N)+1$. Therefore, $b=\depth_R(M\otimes_RN)$; see \ref{cst}(i). This contradiction shows that $\depth_R(M\otimes_RN)\leq b$.
\smallskip
\item Assume $\cx_R(M, N) \leq b$. Set $c= \cx_R(M,N)+1$ and $n=\depth_R(M)+\depth_R(N)-\depth(R)$. Then we have $n\geq c-1$. Suppose contrarily that $\depth_R(M\otimes_RN)\geq n+1$. This implies that $\hh^i_{\fm}(M\otimes_RN)=0$ for all $i=n-c+1, \ldots, n$; see \ref{LCo}(iii). Hence, by Corollary~\ref{cxmm1},  it follows that $\Tor_i^R(M,N)=0$ for all $i\geq 1$. Now, the depth formula yields $\depth_R(M\otimes_RN)=n$, which is a contradiction; see~\ref{df}. This proves that $\depth_R(M\otimes_RN)\leq b$.
\end{enumerate}
\end{rmk}

Our next result, namely Proposition~\ref{ct}, is one of the key applications of Theorem~\ref{t12}. Note that parts~(i) and~(ii) of the proposition were also established in~\cite{Fukano}---after the appearance of~\cite{AOM}---via a different method.

\begin{prop} \label{ct} Let $M$ and $N$ be nonzero $R$-modules such that $\NF(M)\cap\NF(N)\subseteq\{\fm\}$. Consider the following conditions:
\begin{enumerate}[\rm(i)]
\item $R$ is Gorenstein, $M$ and $N$ are maximal Cohen-Macaulay, and $M\otimes_RN$ is not maximal Cohen-Macaulay.
\item $R$ is Gorenstein, and $M$, $N$, and $M\otimes_RN$ are maximal Cohen-Macaulay.
\item $R$ is a complete intersection and $$\cx_R(M,N)\leq  \depth_R(M)+\depth_R(N)-\depth(R) \neq \depth_R(M\otimes_RN).$$
\end{enumerate}
If condition (i) or condition (iii) holds, then $\depth_R(M\otimes_RN)=\inf \left\{i\geq0\mid\widehat{\Tor}_{-i}^R(M,N)\neq0\right\}$. Also, if condition (ii) holds, then $\widehat{\Tor}_{i}^R(M,N)=0$ for all $i=-d+1, \ldots, 0$. 
\end{prop}

\begin{proof} Recall that $\depth_R(M\otimes_RN)=\inf\{i\geq0\mid\hh^i_{\fm}(M\otimes_RN)\neq0\}$; see \ref{LCo}(iii). Hence, to establish both part~(i) and~(iii), it suffices to find an integer $n\geq 1$ such that $\depth_R(M\otimes_RN)\leq n-1$ and 
\begin{equation}\tag{\ref{ct}.1}
\hh^i_{\fm}(M\otimes_RN)\cong\widehat{\Tor}_{-i}^R(M,N) \text{ for  all }  i=0, \ldots, n-1.
\end{equation}

Assume condition~(i) holds. Set $n=\depth_R(M\otimes_RN)+1$. Then $\depth_R(M)+\depth_R(N)\geq\depth(R)+n$ since $M$ and $N$ are maximal Cohen-Macaulay, but $M\otimes_RN$ is not maximal Cohen-Macaulay. We now use Theorem~\ref{t12} for the case where $\mz=\{\fm\}$ and see that the isomorphisms in~(\ref{ct}.1) hold.

Assume condition~(iii) holds. Set $n=\depth_R(M)+\depth_R(N)-\depth(R)$ and $r=\cx_R(M,N)$. If $\depth_R(M\otimes_RN)\leq r$, then $\depth_R(M\otimes_RN) \leq n$ since $r \leq n$. On the other hand, if $\depth_R(M\otimes_RN)\geq r+1$, then it follows from~\ref{cst} that $\depth_R(M\otimes_RN)\leq n$. Thus, in either case, $\depth_R(M\otimes_RN)\leq n$. As we assume $n \neq \depth_R(M\otimes_RN)$, we deduce that $\depth_R(M\otimes_RN)\leq n-1$. Hence, using Theorem~\ref{t12} for the case where $\mz=\{\fm\}$, we see that the isomorphisms in~(\ref{ct}.1) hold.

Assume condition~(ii) holds. As $M\otimes_RN$ are maximal Cohen-Macaulay, it follows that $\hh^i_{\fm}(M\otimes_RN)=0$ for all $i=0, \ldots, d-1$, where $d=\dim(R)$; see~\ref{LCo}. Therefore, by Theorem~\ref{t12}, the claim holds, that is, $\widehat{\Tor}_{-i}^R(M,N)=0$ for all $i=0, \ldots, d-1$. 
\end{proof}

Recall that, if $R$ is a complete intersection and $M$ and $N$ are $R$-modules, then $\cx_R(M,N) \leq \cod(R)$; see~\eqref{cx.1}. Therefore, Proposition~\ref{ctc-intro}, advertised in the introduction, is subsumed by the next result.

\begin{cor} \label{ctc-main} Let $R$ be a local complete intersection ring and let $M$ and $N$ be nonzero $R$-modules such that $\NF(M)\cap\NF(N)\subseteq\{\fm\}$. Set $b=\depth_R(M)+\depth_R(N)-\depth(R)$. Assume $\cx_R(M,N) \leq b$. Assume further $\depth_R(M\otimes_RN)\neq b$. Then the following hold:
\[
 \inf \left\{i\geq0\mid\widehat{\Tor}_{-i}^R(M,N)\neq0\right\}=\depth_R(M\otimes_RN) \leq \cx_R(M,N) \leq b. 
 \]
Therefore, $\depth_R(M\otimes_RN) < \cx_R(M,N)$ or $\cx_R(M,N)<b$.
\end{cor}

\begin{proof} Note that, since $\cx_R(M,N) \leq b$, it follows from Remark~\ref{boundedbounded} that $\depth_R(M\otimes_RN)\leq b$. Thus, $\depth_R(M\otimes_RN)< b$ since we assume $\depth_R(M\otimes_RN)\neq b$. Therefore, by~\ref{cst}(ii), we have that $\depth_R(M\otimes_{R} N)\leq \cx_R(M,N)$. Consequently, Proposition~\ref{ct}(iii) completes the proof.
\end{proof}

In view of part (iii) of Proposition \ref{ct}, it is worth noting that the inequality $\depth_R(M \otimes_R N)\neq b$ does not necessarily imply that $\cx_R(M, N)\leq b$.

\begin{eg} Let $R=k[\![x,y,z]\!]/(xz-y^2, xy-z^2)$, $M=R/(x,y)$, and let $N=R/(x,z)$. Then $R$ is a complete intersection ring, $\dim_R(M\otimes_RN)=\dim_R(M)=\dim_R(N)=0$, and $\NF(M)\cap\NF(N)\subseteq\{\fm\}$. Moreover, $\Tor_1^R(M,N)\neq 0=\Tor_i^R(M,N)$ for all $i\geq 2$, and $\depth_R(M \otimes_R N)=0 \neq -1=b<0=\cx_R(M,N)$; see~\ref{fp}(i) and~\cite[4.2]{Jo1}. 
\end{eg}

An important point in part~(iii) of Proposition~\ref{ct} is that we do \emph{not} assume the vanishing of higher $\Tor_i^R(M, N)$. In fact, the condition $\cx_R(M, N) \leq \depth_R(M) + \depth_R(N) - \depth(R) \neq \depth_R(M \otimes_R N)$ considered in the proposition, as well as its conclusion, fails if all higher Tor modules vanish. These are discussed in the next remark.

\begin{rmk} \label{rmk-ct} Let $R$ be a local complete intersection ring and let $M$ and $N$ be nonzero $R$-modules. Assume $\Tor_i^R(M,N)=0$ for all $i\gg 0$.
\begin{enumerate}[\rm(i)]
\item It follows that $\widehat{\Tor}_{i}^R(M,N)=0$ for all $i\in \ZZ$; see~\ref{fp}(i). Since $M$ and $N$ are not zero, $\depth_R(M\otimes_RN)$ is a nonnegative integer and hence is not equal to $\inf \big\{i\geq0\mid\widehat{\Tor}_{-i}^R(M,N)\neq0\big\}=\infty$. 

\item Note that $\cx_R(M,N)=0$; see~\ref{fp}(i). Set $q=\sup\{i\geq 0 \mid\Tor_i^R(M,N)\neq0\}$. We have, by~\ref{fp}(iii), that:
\[\tag{\ref{rmk-ct}.1}
q=\sup\big\{\depth(R_\fp)-\depth_{R_{\fp}}(M_\fp)-\depth_{R_{\fp}}(N_\fp): \fp\in\Spec(R)\big \}.
\]
Set $b=\depth_R(M)+\depth_R(N)-\depth(R)$ and assume $\NF(M)\cap\NF(N)\subseteq\{\fm\}$. Note that, $b\geq 0$ if $b\geq \cx_R(M,N)$. It follows from~(\ref{rmk-ct}.1) that $q=0$. This implies that $\Tor_i^R(M,N)=0$ for all $i\geq 1$. So, the depth formula holds and this gives $b=\depth_R(M\otimes_RN)$; see~\ref{df}. 
\end{enumerate}
\end{rmk}

The following example demonstrates that the hypothesis $\NF(M) \cap \NF(N) \subseteq \{\fm\}$ is necessary for part~(iii) of Proposition~\ref{ct}.

\begin{eg} Let $R=k[\![x,y,z,u]\!]/(xz-yu)$ and $M=N=R/(x)$. Then $\pd_R(M)=1$ and $\depth_R(M)=2$. Thus $\cx_R(M,N)=0<1= 2+2-3=b \neq 2=\depth_R(M\otimes_RN)$. On the other hand, $(x,y) \in \NF(M) \not\subseteq \{\fm\}$. Note that $\depth_R(M\otimes_RN)$ is not equal to $\inf \big\{i\geq0\mid\widehat{\Tor}_{-i}^R(M,N)\neq0\big\}=\infty$; see Remark~\ref{rmk-ct}.
\end{eg}

The following facts are used in the proof of Corollaries~\ref{ct-cor} and~\ref{cq}:

\begin{chunk} \label{TVFPD} Let $R$ be a local complete intersection and let $M$ and $N$ be $R$-modules. Assume $R$ is hypersurface, or $M=N$. If $\Tor_i^R(M,N)=0$ for all $i\gg 0$, then $\pd(M)<\infty$ or $\pd(N)<\infty$; see~\cite[1.9]{HW2} and~\cite[1.2]{Jo2}. 
\end{chunk}

\begin{cor} \label{ct-cor} Let $R$ be a $d$-dimensional local complete intersection ring of codimension $c$ such that $d\geq c+1$, and let $M$ and $N$ be nonfree maximal Cohen-Macaulay $R$-modules. Assume $M$ and $N$ are locally free on the punctured spectrum of $R$ (for example, $R$ is an isolated singularity). Then $$\depth_R(M\otimes_RN)=\inf\left\{i\geq 0: \widehat{\Tor}_{-i}^R(M,N)\neq 0 \right\}$$ provided that at least one of the following condition holds:
\begin{enumerate}[\rm(i)]
\item $M=N$.
\item $\cx_R(M)+\cx_R(N)>c$.
\end{enumerate}
\end{cor}

\begin{proof} If $\cx_R(M,N) \leq \depth_R(M)+\depth_R(N)-\depth(R) \neq \depth_R(M\otimes_RN)$, then the claim follows from Proposition~\ref{ct}(iii). As $\cx_R(M,N) \leq c<c+1\leq d=\depth_R(M)+\depth_R(N)-\depth(R)$, it suffices to show that $M\otimes_RN$ is not maximal Cohen-Macaulay; see \eqref{cx.1}. Suppose $M\otimes_RN$ is maximal Cohen-Macaulay, and seek a contradiction.

Let $\fp \in \Spec(R)$ such that $\dim(R_{\fp})\leq c$. Then $M_{\fp}$ is free over $R_{\fp}$ since $\dim(R_{\fp})\leq d-1$. Thus, $M_{\fp}$ is free over $R_{\fp}$ for all $\fp \in \X^{c}(R)$. Since $M$, $N$, and $M\otimes_RN$ are maximal Cohen-Macaulay, they satisfy $(S_n)$ for all $n\geq 0$; see~\ref{SL}. Therefore,~\ref{DD} shows that $\Tor_i^R(M,N)=0$ for all $i\geq 1$. 

For part (i), if $M=N$, then the vanishing of $\Tor_i^R(M,M)$ implies that $M$ is free and yields a contradiction due to~\ref{TVFPD}. For part (ii), the assumption that $\cx_R(M)+\cx_R(N)>c$ yields a contradiction because the vanishing of $\Tor_i^R(M,N)$ for all $i\geq 1$ implies that $\cx_R(M)+\cx_R(N)=\cx_R(M\otimes_RN)\leq \cod(R)=c$; see~\cite[2.1]{Mi}. Consequently, $M\otimes_RN$ cannot be maximal Cohen-Macaulay, and the claims follow.
\end{proof}

\begin{eg} Let $R=\CC[\![x,y,z,u,v]\!]/(x^2+y^2+z^2+u^2+v^2, x^2+2y^2+3z^2+4u^2+5v^2)$. Then $R$ is a three-dimensional local complete intersection ring of codimension two. So, by Corollary~\ref{ct-cor}, we have that $\depth_R(M\otimes_RM)=\inf\{i\geq 0: \widehat{\Tor}_{-i}^R(M,M)\neq0\}$ for each non-free maximal Cohen-Macaulay $R$-module $M$.
\end{eg}

Next we proceed to investigate torsion submodules of tensor products of nonfree maximal Cohen-Macaulay modules over hypersurface singularities; we make use of Proposition~\ref{ct} and the following:

\begin{chunk} \label{SRT} If $R$ is a local hypersurface ring with an isolated singularity, and $M$ and $N$ are nonfree maximal Cohen-Macaulay $R$-modules, then $\depth_R(M\otimes_RN)\leq1$; see~\cite[4.7]{CST}. 
\end{chunk}

\begin{cor} \label{corcombined} Let $R$ be a $d$-dimensional local hypersurface ring, with $d\geq 2$, and let $M$ and $N$ be nonfree maximal Cohen-Macaulay $R$-modules. If $R$ is an isolated singularity, then one of the following holds:
\begin{enumerate}[\rm(i)]
\item $\depth_R(M\otimes_RN)=1$, $M\otimes_RN$ is torsion-free, and $\Tor_{2i-1}^R(M,N) \neq 0 = \Tor_{2i}^R(M,N)$ for all $i\geq 1$.
\item $\depth_R(M\otimes_RN)=0$, $M\otimes_RN$ has (nonzero) torsion, and $\Tor_{2i}^R(M,N) \neq 0 $ for all $i\geq 1$.
\end{enumerate}
\end{cor}
\begin{proof} 

Note that $\depth_R(M\otimes_RN)\leq 1$ since $M$ and $N$ are nonfree maximal Cohen-Macaulay $R$-modules; see~\ref{SRT}. Therefore, since $\dim(R)\geq 2$, $M\otimes_R N$ is not maximal Cohen-Macaulay. Hence, Proposition~\ref{ct}(i) shows that $\depth_R(M\otimes_RN)=\inf\{i\geq0 : \widehat{\Tor}_{-i}^R(M,N)\neq0\}$. As $M$ is locally free on $\Spec(R)-\{\fm\}$ and $\depth_R(M\otimes_RN)\leq 1$, it follows that $M\otimes_RN$ is torsion-free if and only if $\depth_R(M\otimes_RN)=1$. 

If $M\otimes_RN$ is torsion-free, then $\widehat{\Tor}_{0}^R(M,N)=0 \neq \widehat{\Tor}_{-1}^R(M,N)$. On the other hand, if $M\otimes_RN$ is not torsion-free, that is, if $M\otimes_RN$ has torsion, then $\depth_R(M\otimes_RN)=0$ so that $\widehat{\Tor}_{0}^R(M,N)\neq 0$. Without loss of generality, we may assume $M$ has no nonzero free direct summand. Thus, $M\cong \Omega^2_R(M)$; see~\cite[6.1(ii)]{Ei}. Hence, the claims about the vanishing of Tors hold; see~\ref{Tate hom}(v). 
\end{proof}

We point out by the next example that the Tor vanishing conclusion in Corollary~\ref{corcombined} can fail if the modules in question are not maximal Cohen-Macaulay.

\begin{eg} Let $R=k[\![x,y,w,z]\!]/(xw-yz)$ and let $M$ be the cokernel of the following map
\[
\xymatrixcolsep{13pt}
\xymatrix{
R \ar[rrr]^-{\left( \begin{array}{c} w\; y\; x\; z \end{array} \right)^T} &&& R^{\oplus 4}
}
\]

Then $R$ is a three-dimensional hypersurface ring which has an isolated singularity, $M$ is a torsion-free $R$-module which is not maximal Cohen-Macaulay ($\depth_R(M)=2$ and $\pd_R(M)=1$), and $\Tor^{R}_{i}(M,M)=0$ for all $i\geq 1$; see~\cite[4.5]{Ce}. It follows from~\ref{df} that $\depth_R(M\otimes_RM)=1$.
\end{eg}

%%%%%%%%%%%%%%%%%%%%%%%%%%%%%%%%%%%%%%%%%%%%%%%%%%%%%%%%%%%%%%%%%
%%%%%%%%%%%%%%%%%%% SECTION  %%%%%%%%%%%%%%%%%%%%%%%%%%%%%%%%%%%%%%%%
%%%%%%%%%%%%%%%%%%% SECTION  %%%%%%%%%%%%%%%%%%%%%%%%%%%%%%%%%%%%%%%%
%%%%%%%%%%%%%%%%%%% SECTION  %%%%%%%%%%%%%%%%%%%%%%%%%%%%%%%%%%%%%%%%
%%%%%%%%%%%%%%%%%%%%%%%%%%%%%%%%%%%%%%%%%%%%%%%%%%%%%%%%%%%%%%%%%

\section{Generalizations of certain results of Huneke and Wiegand} \label{section5}

Let $\mathcal{E}$ be a vector bundle on a projective scheme $X$ equipped with a very ample line bundle $\mathcal{O}(1)$. Denote the dual bundle $\mathcal{H}om(\mathcal{E}, \mathcal{O}_X)$ by $\mathcal{E}^*$. Define the graded cohomology module
\[
\hh^i_*(X, \mathcal{E} \otimes \mathcal{E}^*) = \bigoplus_{n \in \mathbb{Z}} \hh^i\big(X, \mathcal{E} \otimes \mathcal{E}^*(n)\big).
\]
We say that $\mathcal{E}$ has the \emph{cohomological rigidity property} provided that the following implication holds: If $\hh^i_*(X, \mathcal{E} \otimes \mathcal{E}^*) = 0$ for some $i \geq 1$, then $\mathcal{E}$ is trivial---that is, $\mathcal{E}$ is isomorphic to a direct sum of line bundles.

The cohomological rigidity property was first studied by Kempf~\cite{Kem} and Luk-Yau~\cite{Luk} using different approaches. A commutative algebra interpretation of this property was explored by Huneke and Wiegand~ \cite{HW2}. For a reflexive module $M$ over a hypersurface ring $R$, Huneke and Wiegand established new criteria for the freeness of $M$ in terms of the vanishing of $\hh_{\fm}(M \otimes_R M^\ast)$. They also generalized some of the results of Luk-Yau~\cite{Luk} using commutative algebra techniques.

This section is motivated by the work of Huneke and Wiegand~\cite{HW2}. Our main goal is to present new criteria for detecting the freeness of certain modules $M$ via the vanishing of local cohomology of $M \otimes_R M^\ast$. In particular, we aim to extend the results found in~\ref{HWL1}. The key tools we use are Lemma~\ref{grz} and Proposition~\ref{l2}. Additionally, we study arithmetically Cohen-Macaulay vector bundles on complete intersections, as described in Proposition~\ref{artCM}.

Recall that an $R$-module $M$ is said to be \emph{locally free on the punctured spectrum of $R$} if $M_{\fp}$ is free over $R_{\fp}$ for all $\fp \in \Spec(R) \setminus \{\fm\}$. For such an $R$-module $M$, note that if $\mz = \V(\fm)$, then $\NF(M) \subseteq \mz$. We define $(-)^{\vee} = \Hom_R\left(-, \E_R(k)\right)$, where $\E_R(k)$ denotes the injective hull of the residue field $k$ of the local ring $R$.

\begin{chunk} [{Huneke-Wiegand \cite[4.1 and 4.2]{HW2}}] \label{HWL1} \label{HWL3} Let $R$ be a $d$-dimensional regular local ring, with $d\geq 2$, and let $M$ be an $R$-module which is locally free on the punctured spectrum of $R$.
\begin{enumerate}[\rm(i)]
\item If $\hh^1_{\fm}(M\otimes_RM^*)=0$, then $M \cong F\oplus T$ for some free $R$-module $F$ and a finite length module $T$.
\item If $\hh^2_{\fm}(M\otimes_RM^*)=0$, then $M^{\ast}$ is free.
\item If $d\geq 3$ and $N$ is an $R$-module which is locally free on the punctured spectrum of $R$, then,  for each integer $i$ with $2 \leq i \leq d-1$, it follows that: $$\hh_{\fm}^i(M^\ast\otimes_RN^\ast)\cong \big(\hh_{\fm}^{d-i+1}(M \otimes_RN )\big)^{\vee}.$$
\end{enumerate}
\end{chunk}

Our aim is to generalize the first and third parts of~\ref{HWL1}; see Theorems~\ref{t7} and~\ref{tomax}. We also study the vanishing of $\hh^i_{\fm}(M \otimes_R M^*)$ and establish results analogous to the second part of~\ref{HWL1}, even over non-regular rings; see Corollaries~\ref{remzidenoturu} and~\ref{tp3c1}, as well as Proposition~\ref{artCM} for related applications of geometric nature. These results are developed across three separate subsections.

%%%%%%%%%%%%%%%%%%%%%%%%%%%%%
\subsection*{A generalization of the first part of \ref{HWL1}} We proceed to generalize~\ref{HWL1}(i). For that, we first prepare a lemma which extends~\cite[1.4.19]{BH}.  Recall that, given an $R$-module $M$ and a specialization-closed subset $\mz$ of $\Spec(R)$, we have $\gr_R(\mz,M)=\inf\{i \geq 0 : \hh_{\mz}^{i}(M)\neq0\}$; see~\ref{gsc}.

\begin{lem}\label{gr} Let $\mz$ be a specialization-closed subset of $\Spec(R)$, and let $M$ and $N$ be $R$-modules, where $N$ is not necessarily finitely generated. Then $\gr_R\big(\mz,\Hom_R(M,N)\big)\geq\min\{2,\gr_R(\mz,N)\}$
\end{lem}

\begin{proof} Let $F_1\to F_0\to M\to0$ be a free presentation of $M$. This yields exact sequences of $R$-modules: 
$$0\to\Hom_R(M,N)\to\Hom_R(F_0,N)\to X\to0\ \ \text{and}\ \ 0\to X\to\Hom_R(F_1,N).$$
We apply~\ref{LC}(ii) to obtain the following exact sequences of $R$-modules:
\[
0 \rightarrow \Gamma_{\mz}\big(\Hom_R(M,N)\big)\rightarrow\Gamma_{\mz}\big(\Hom_R(F_0,N)\big)\to\Gamma_{\mz}(X)\to\hh^1_{\mz}\big(\Hom_R(M,N)\big)\to\hh^1_{\mz}\big(\Hom_R(F_0,N)\big),
\]
and
\[
0\to\Gamma_{\mz}(X)\to\Gamma_{\mz}(\Hom_R(F_1,N)).
\]
Now the assertion follows from the above exact sequences.
\end{proof}

Next we record some facts that are used in the proof of Theorem~\ref{t7}; see also~\ref{a1} and~\ref{cs}.

\begin{chunk} \label{lan} Let $M$ be an $R$-module. Then there is an (four-term) exact sequence of $R$-modules 
\begin{equation}\tag{\ref{lan}.1}
0\rightarrow\Ext^1_R(\Tr_R M,N)\rightarrow M\otimes_RN\rightarrow\Hom_R(M^*,N)\to\Ext^2_R(\Tr_R M,N)\to0,
\end{equation}
where the middle map is the evaluation map; see~\cite[2.6]{AuBr}. Setting $N=R$, we obtain an exact sequence
\begin{equation}\tag{\ref{lan}.2}
0\rightarrow\Ext^1_R(\Tr_R M,R)\rightarrow M \rightarrow M^{**}\to\Ext^2_R(\Tr_R M,R)\to0, 
\end{equation}
where the middle map is the natural one.

Given an integer $n\geq 1$, there is an exact sequence~\cite[2.8]{AuBr}:
\begin{equation}\tag{\ref{lan}.3}
\Tor_2^R\left(\Tr_R\Omega_R^nM,N\right)\rightarrow\Ext^n_R(M,R)\otimes_RN\rightarrow\Ext^n_R(M,N)\rightarrow\Tor_1^R\left(\Tr_R\Omega_R^nM,N\right)\rightarrow0.
\end{equation}

Map minimally a free $R$-module $F$ onto $M$. Then the cokernel of the composition $M \to M^{\ast\ast} \hookrightarrow F^{\ast} \cong F$ is called the first \emph{cosyzygy} of $M$ (or pushforward of $M$~\cite{HJW}), and is denoted by $\Omega_R^{-1} M$. As $\Tr_R \Tr_R M\approx M$, it follows from~\cite[4.1]{ST2} that: 
\begin{equation}\tag{\ref{lan}.4}
\Omega_R^{-i}M \approx \Tr_R\Omega_R^i\Tr_R M \text{ for all } i\geq 1. 
\end{equation}
Hence, by replacing $M$ with $\Tr_R M$ in (\ref{lan}.3) and using (\ref{lan}.4), we obtain the following exact sequence:
\begin{equation}\tag{\ref{lan}.5}
\Tor_2^R(\Omega_R^{-n}M,N) \rightarrow \Ext^n_R(\Tr_R M,R)\otimes_RN\rightarrow\Ext^n_R(\Tr_R M,N)\rightarrow\Tor_1^R(\Omega_R^{-n}M,N)\rightarrow0.
\end{equation}
\end{chunk}

In the following we generalize~\ref{HWL1}(i):

\begin{thm}\label{t7} Let $R$ be a local ring, $M$ be an $R$-module, and let $\mz$ be a specialization-closed subset of $\Spec(R)$. Assume the following conditions hold:
\begin{enumerate}[\rm(i)]
\item $\NF(M)\subseteq\mz$. 
\item $\gr_R(\mz,R)\geq1$.
\item $\hh^1_{\mz}(M\otimes_RM^*)=0$.  
\end{enumerate}
Then $M^{\ast}$ is free and $M\cong M^{\ast\ast}\oplus \Ext^1_R(\Tr_R M,R)$. So, setting $F=M^{\ast\ast}$ and $T= \Ext^1_R(\Tr_R M,R)$, we have that $M\cong F \oplus T$ where $F$ is free, and $T$ is torsion with respect to $\mz$. Therefore, if $M$ is torsion-free with respect to $\mz$, then $M$ is free.
\end{thm}

\begin{proof} We have by Lemma~\ref{gr} that $\gr_R(\mz,M^*)\geq\min\{2,\gr_R(\mz,R)\}\geq1$. So, Proposition~\ref{l2}(ii) implies that $\Ext^2_R(\Tr M,M^*) \hookrightarrow \hh^1_{\mz}(M\otimes_RM^*)$. Therefore, since $\hh^1_{\mz}(M\otimes_RM^*)=0$, we see that $\Ext^2_R(\Tr M,M^*)=0$. Hence,~(\ref{lan}.5) gives $\Tor_1^R(\Omega_R^{-2}M,M^{\ast})=0$.

Recall that $\Omega_R^2\Tr_R M \cong M^{\ast}$; see~\ref{a1}. Therefore, by~(\ref{lan}.4), we see that $\Omega_R^{-2}M \cong \Tr_R \Omega_R^2 \Tr_R M \cong \Tr_R M^*$. Therefore, since $\Tor_1^R(\Omega_R^{-2}M,M^{\ast})=0$ and $\Omega_R^{-2}M \cong \Tr_R M^*$, we obtain $\Tor_1^R(\Tr_R M^*,M^*)=0$. Hence $M^{\ast}$ is free; see~\ref{a1}. Consequently, the fact that $\Ext^2_R(\Tr_R M,M^*)=0$ implies $\Ext^2_R(\Tr_R M,R)=0$. So, the sequence~(\ref{lan}.2) yields the exact sequence $0\rightarrow\Ext^1_R(\Tr_R M,R)\rightarrow M\rightarrow M^{**}\rightarrow0$. 
As $M^{**}$ is free, the exact sequence splits and gives an isomorphism $M\cong M^{\ast\ast}\oplus \Ext^1_R(\Tr_R M,R)$. 

Set $T=\Ext^1_R(\Tr_R M,R)$. Then $\Supp_R(T)\subseteq \mz$ since $\NF(M)\subseteq\mz$. So, $T$ is torsion with respect to $\mz$; see~\ref{LC}(iv). If $M$ is torsion-free with respect to $\mz$, since $T \hookrightarrow M$, we see that $T=\Gamma_{\mz}(T)\hookrightarrow \Gamma_{\mz}(M)=0$, and hence $M \cong M^{\ast\ast}$, that is, $M$ is free.
\end{proof}

In passing, let us note that we can prove the result of Huneke-Wiegand stated as~\ref{HWL1}(i) by using Theorem~\ref{t7}. Note also that our proof of Theorem~\ref{t7} does not appeal to the proof of \ref{HWL1}(i) given in~\cite{HW2}. 

\begin{cor}\label{c1} Let $R$ be a Cohen-Macaulay local ring and let $M$ be an $R$-module which is locally free on the punctured spectrum of $R$. If $\depth_R(M)\geq 1$ and $\hh^1_{\fm}(M\otimes_RM^*)=0$, then $M$ is free.
\end{cor}

\begin{proof} Set $\mz=\V(\fm)$. Then $\mz$ is a specialization-closed subset of $\Spec(R)$ and $\NF(M)\subseteq\mz$. Moreover, $\gr_R(\mz,R)=\depth(R)=\dim(R)\geq \depth_R(M)=\gr_R(\mz,M)\geq 1$. Hence, if $\hh^1_{\fm}(M\otimes_RM^*)=0$, then Theorem~\ref{t7} shows that $M \cong F\oplus T$ for some free $R$-module $F$ and some finite length $R$-module $T$; in that case, $T=0$ since $\depth_R(M)\geq 1$. This proves the claim.
\end{proof}

We see by the next example that the locally freeness assumption is necessary for Corollary~\ref{c1}.

\begin{eg}\label{directsum} Let $R=k[\![x,y,z]\!]/(x^2)$ and let $M=R/xR$. Then $M$ is not free and $\depth_R(M)=\dim_R(M)=2$. Moreover, $M \cong M^\ast$ and $\hh^{ 1}_{\fm} (M\otimes_RM^\ast)=\hh^{ 1}_{\fm} (M)=0$. Note that $M$ is not locally free on the punctured spectrum of $R$; for example, $M_{\fp}$ is not free when $\fp=(x)\in \Spec(R)$.
\end{eg}

\begin{rmk}\label{directsum2} The vanishing of $\hh^1_{\fm}(M\otimes_RM^*)$ is needed for the conclusion of Corollary~\ref{c1} to hold, even if $\hh^n_{\fm}(M\otimes_RM^*)$ vanishes for other values of $n$. We discuss this next:
\begin{enumerate}[\rm(i)]
\item Let $R$ be a $d$-dimensional Cohen-Macaulay local ring, with $d\geq 3$, admitting a canonical module $\omega_R$. Assume $R$ is not regular and set $M=\fm$. Then $M$ is not free, is locally free on the punctured spectrum of $R$, $\depth_R(M)=1$, and $\dim_R(M)=d$. Moreover, as $\Ext^1_R(k,R)=0$, it follows that $M^{\ast} \cong R$. Thus, $M\otimes_RM^{\ast} \cong M$. Therefore, $$\hh^{0}_{\fm}(M\otimes _R M^\ast)=0=\hh^{j}_{\fm}(M\otimes _R M^\ast) \text{ for all }j\geq d+1 \; \text{ and }\; \hh^{1}_{\fm}(M\otimes _R M^\ast)\neq 0\neq \hh^{d}_{\fm}(M\otimes _R M^\ast).$$
The short exact sequence $0 \to \fm \to R \to k \to 0$ yields an exact sequence $\hh^i_{\fm}(k) \to \hh^{i+1}_{\fm}(\fm) \to \hh^{i+1}_{\fm}(R)$ for each $i\geq 1$. Since $\hh^i_{\fm}(k)=0$ for all $i\geq 1$, and $\hh^{j}_{\fm}(R)=0$ if and only if $j\neq d$, we conclude that $\hh^{i}_{\fm}(M\otimes _R M^\ast)=0$ for all $i=2, \ldots, d-1$. 
\smallskip

\item One can also construct examples for non-regular rings analogous to the one in part (i). For example, if $R$ is a $d$-dimensional non-Gorenstein Cohen-Macaulay local ring, with $d\geq 3$, admitting a canonical module $\omega_R$ such that $R_{\fp}$ is Gorenstein for all $\fp \in \Spec(R)-\{\fm\}$ and $M=\omega_R$, then $M$ is not free and it follows from~\cite[5.8]{MA} that $\hh^{1}_{\fm}(M\otimes _R M^\ast)\neq 0 \neq \hh^{d}_{\fm}(M\otimes _R M^\ast)$ and $\hh^{i}_{\fm}(M\otimes _R M^\ast)=0$ for all $i=2, \ldots, d-1$.
\end{enumerate}
\end{rmk}

%%%%%%%%%%%%%%%%%%%%%%%%%%%%%
\subsection*{New results similar to the second part of~\ref{HWL1}} We give applications of Proposition~\ref{l2} and prove results similar to~\ref{HWL1}(ii) over rings which are not necessarily regular. Our results show that the regularity hypothesis from~\ref{HWL1}(ii) can be removed for Tor-rigid modules, or over complete intersection rings, under certain conditions; see~\ref{TorRigid}, Proposition~\ref{tp3}, and Corollary~\ref{tp3c1}.
%%%%%%%%%%%%%%%%%%%%%%%%%%%%

\begin{lem} \label{remzi} Let $R$ be a local Cohen-Macaulay ring, $\mz$ be a specialization-closed subset of $\Spec(R)$, and let $M$ be an $R$-module. If $\NF(M)\subseteq\mz$ and $\gr_R(\mz,M)\geq 2$, then $M$ is reflexive. 
\end{lem}

\begin{proof} Note that, by Lemma~\ref{grz}, we have $\gr_R(\mz,M)=\inf\{\depth_{R_\fp}(M_\fp) : \fp\in\mz\}$. This implies that $\depth_{R_\fp}(M_\fp) \geq 2$ for all $\fp \in \mz$. So, if $\fp\in\X^{1}(R)$, then $\fp \notin \mz$ so that $\fp \notin \NF(M)$, that is, $M_\fp$ is free.

Next let $\fq\in \Spec(R)$ such that $\depth(R_\fq)\geq 2$. If $\fq \notin \mz$, then $\fq \notin \NF(M)$ and so $\depth_{R_\fq}(M_\fq)\geq 2$. On the other hand, if $\fq \in \mz$, then 
$\depth_{R_\fq}(M_\fq) \geq \gr_R(\mz,M)\geq 2$. Thus, by~\cite[1.4.1(b)]{BH}, $M$ is reflexive.
\end{proof}

\begin{lem} \label{remzii} Let $R$ be a Cohen-Macaulay local ring, $\mz$ be a specialization-closed subset of $\Spec(R)$, and let $M$ be an $R$-module. Assume:
\begin{enumerate}[\rm(i)]
\item $\NF(M)\subseteq\mz$.
\item $2\leq \gr_R(\mz,M)$ and $\hh^n_{\mz}(M\otimes_RM^*)=0$ for some integer $n$ where $2\leq n\leq\gr_R(\mz,M)$.
\end{enumerate} 
Then $\Ext^{n-1}_R(M,M)=0$. 
\end{lem}

\begin{proof} We have, by Proposition~\ref{l2}, an injective map $\Ext^{n+1}_R(\Tr_R M^{\ast}, M)\hookrightarrow\hh^{n}_{\mz}(M\otimes_RM^{\ast})$. Hence  $\Ext^{n+1}_R(\Tr_R M^{\ast}, M)=0$. Note that $M$ is reflexive due to Lemma~\ref{remzi}. Thus, $M\cong M^{**}\cong \Omega_R^2\Tr_R M^*$; see \ref{a1}. Consequently, $\Ext^{n-1}_R(M,M) \cong \Ext^{n-1}_R(\Omega_R^2 \Tr_R M^{\ast},M) \cong \Ext^{n+1}_R(\Tr_R M^{\ast}, M)=0$, as claimed.
\end{proof}

Here is an example showing that the condition $n\leq\gr_R(\mz,M)$ is necessary for Lemma~\ref{remzii}.

\begin{eg} Let $R=k[\![x,y]\!]$, $M=\Syz_R^1k=\fm$, and let $\fa=xR$ for some $0\neq x \in \fm$. Set $\mz=\V(\fa)$. Then $\gr_R(\mz,M)=1$. Moreover, $\hh^i_{\fa}(M\otimes_RM^*)=0$ for all $i\geq 2$ since $\fa$ is principal; see~\cite[3.3.1]{Bsh}. However, $\pd_R(M)=1$ and so $\Ext^{1}_R(M,M)\neq 0$; see~\cite[Lemma 1(iii), page 154]{Mat}.
\end{eg}

Our main aim is to use Lemma \ref{remzii} to prove Corollary~\ref{remzidenoturu} which is necessary for Proposition~\ref{cortomax}. However, Lemma~\ref{remzii} has other useful applications; see Proposition~\ref{tp3}. 

\begin{chunk} \textbf{Tor-rigidity}. \label{TorRigid} An $R$-module $M$ is called \emph{Tor-rigid} provided that the vanishing of a single $\Tor^R_n(M,N)$ for each $R$-module $N$ and each $n\geq 1$ implies the vanishing of $\Tor^R_i(M,N)$ for all $i\geq n$. 

Clearly, over a local ring, every syzygy module of the residue field is Tor-rigid. Results of Auslander~\cite{Au} and Lichtenbaum~\cite{Li} show that, if $R$ is a regular local ring, then each $R$-module is Tor-rigid. On the other hand, there are non-regular rings over which each $R$-module is Tor-rigid; see~\cite[2.8 and 3.16]{Da1}. 

Tor-rigidity is a rather subtle topic in commutative algebra and remains one of the themes yet to be fully understood. As some of our results rely on Tor-rigidity, we collect a few facts and examples about Tor-rigid modules, and refer the reader to~\cite{CeRo, DaoSurvey} for further details and examples. 

In the following remarks, we assume $R$ is a local ring.

\begin{enumerate}[\rm(i)] 
\item If $M$ is a Tor-rigid $R$-module and $\Ext^n_R(M,M)=0$ for some $n\geq 1$, then the work of Auslander ~\cite{Au} and Jothilingam \cite{Jot} implies that $\pd_R(M)<n$; see, for example,~\cite[3.1.2]{DaoSurvey}.
\item Let $R$ be a hypersurface ring which is quotient of an unramified regular ring (for example, a regular ring containing a field). Then each $R$-module which has finite projective dimension, or has finite length is Tor-rigid; see~\cite[2.4]{HW1} and~\cite[Thm. 2]{Li}.  
\item Let $I$ be an ideal of $R$ such that $\Tor^R_n(R/I,M)=0$ for some $R$-module $M$ and integer $n\geq 0$. Then $\pd(M)\leq n-1$ and so $R/I$ is a Tor-rigid $R$-module provided that one of the following conditions holds:
\begin{enumerate}[\rm(a)] 
\item $I$ is integrally closed and $\fm$-primary; see~\cite[3.3]{CHKV}.
\item $\depth(R)\geq 1$ and $I=(A:_R\fm)$ for some $\fm$-primary ideal $A$ of $R$ such that  $A\neq \fm$; see~\cite[2.2]{CSTT} and~\cite[2.9 and 2.10]{CeTos}.
\end{enumerate}
\end{enumerate}
\end{chunk}

We present a consequence of Lemma~\ref{remzii}:

\begin{cor} \label{remzidenoturu} Let $R$ be a Cohen-Macaulay local ring and let $M$ be an $R$-module. Assume $\depth_R(M)\geq 2$, $M$ is locally free on the punctured spectrum of $R$, and $n$ is an integer such that $2\leq n\leq\depth_R(M)$.
\begin{enumerate}[\rm(i)]
\item If $\hh^n_{\fm}(M\otimes_RM^{\ast})=0$, then $\Ext^{n-1}_R(M,M)=0$.
\item If $\hh^n_{\fm}(M\otimes_RM^{\ast})=0$ and $M$ is Tor-rigid, then $\pd_R(M)<n-1$. 
\end{enumerate}
\end{cor}

\begin{proof} We set $\mz=\V(\fm)$ and obtain by Lemma \ref{remzii} that $\Ext^{n-1}_R(M,M)=0$. If $M$ is Tor-rigid, this yields $\pd_R(M)<n-1$; see~\ref{TorRigid}(i).
\end{proof}

As syzygies of finite-length modules are Tor-rigid over certain hypersurfaces, the following holds due to~\ref{TorRigid}(i) and Corollary~\ref{remzidenoturu}.

\begin{cor} Let $R=k[\![x_1, \ldots, x_{d+1}]\!]/(f)$ be a hypersurface ring, with $d\geq 2$, and let $M=\Omega_R^d(L)$ for some finite-length $R$-module $L$. If $\hh^n_{\fm}(M\otimes_RM^{\ast})=0$ for some $n$ with $2\leq n \leq d$, then $M$ is free. 
\end{cor}

Tor-rigidity is a necessary hypothesis in Corollary~\ref{remzidenoturu}, as the following example shows.

\begin{eg} \label{TRE} Let $R=k[\![x,y,u,v]\!]/(xy-uv)$ and let $M$ be the ideal of $R$ generated by $x$ and $u$. Then $M$ is a non-free maximal Cohen-Macaulay $R$-module such that $M$ is locally free on the punctured spectrum of $R$ and $\hh_{\fm}^2(M\otimes_RM^{\ast}) \cong \hh_{\fm}^2(\fm)=0$; see~\cite[1.8]{HW2}. Moreover, $M$ is not a Tor-rigid $R$-module since $\Tor^R_2(M,M^{\ast})=0\neq \Tor^R_3(M,M^{\ast})$.
\end{eg}

Next, we proceed to investigate the vanishing of $\hh^i_{\mz}(M\otimes_RM^*)$ over complete intersection rings. Since the vanishing of local cohomology is connected to the vanishing of $\Ext$ via Proposition~\ref{l2}, the following facts are of interest:

\begin{chunk} \label{JorJot} Let $R$ be a local complete intersection and let $M$ be an $R$-module. If $\Ext^n_R(M,M)=0$ for some even integer $n\geq 2$, then $\pd_R(M)<\infty$; see~\cite[4.2]{AvBu}. Moreover, if $\Ext^2_R(M,M)=0$, then $\pd_R(M)\leq 1$; see~\cite[2.5]{JorExt}.
\end{chunk}

\begin{prop}\label{tp3} Let $R$ be a $d$-dimensional local complete intersection ring, $\mz$ be a specialization-closed subset of $\Spec(R)$, and let $M$ be an $R$-module such that $\NF(M)\subseteq\mz$. Set $g=\gr_R(\mz,M)$. Then the following hold:
\begin{enumerate}[\rm(1)]
\item If $n$ is an odd integer such that $1\leq n\leq g$, and $\hh^{n}_{\mz}(M\otimes_RM^*)=0$, then $\pd_R(M)<\infty$. 
\item Assume $d\geq 4$ and $g\geq3$. 
\begin{enumerate}[\rm(i)]
\item If $\hh^3_{\mz}(M\otimes_RM^*)=0$, then $\pd_R(M)\leq1$. 
\item If $\hh^2_{\mz}(M\otimes_RM^*)=\hh^3_{\mz}(M\otimes_RM^*)=0$, then $M$ is free.
\end{enumerate}
\end{enumerate}
\end{prop}

\begin{proof} (1) Assume $n=1$ so that $\hh^{1}_{\mz}(M\otimes_RM^*)=0$. Note that $\gr_R(\mz,R)\geq g \geq 1$ and $M$ is torsion-free with respect to $\mz$; see~\ref{LC}(iii) and Lemma~\ref{grz}. Therefore, $M$ is free by Theorem~\ref{t7}. This shows we may assume $n\geq 2$. In that case, Lemma~\ref{remzii} shows that $\Ext^{n-1}_R(M,M)=0$. As $R$ is a complete intersection and $n-1$ is even, we conclude that $\pd_R(M)<\infty$; see~\ref{JorJot}.

(2) Note that, since $g \geq 3$, Proposition \ref{l2} (the case where $n=3$) yields:
\begin{equation}\tag{\ref{tp3}.1}
\hh^{2}_{\mz}(M^{\ast}\otimes_RM)\cong\Ext^{3}_R(\Tr_R M^{\ast}, M) \text{ and } \Ext^{4}_R(\Tr_R M^{\ast}, M)\hookrightarrow\hh^{3}_{\mz}(M^{\ast}\otimes_RM).
\end{equation}
Recall that $M$ is reflexive and $M\cong M^{\ast\ast} \cong \Omega_R^2 \Tr_R M^{\ast}$; see \ref{a1} and Lemma \ref{remzi}.
Thus, if $\hh^3_{\mz}(M\otimes_RM^*)=0$, then we see from (\ref{tp3}.1) that $\Ext^2_R(M^{\ast\ast},M^{\ast\ast})=\Ext^2_R(M,M)=0$. In that case, $\pd_R(M)\leq1$ due to \ref{JorJot}. 

If $\hh^2_{\mz}(M\otimes_RM^*)=\hh^3_{\mz}(M\otimes_RM^*)=0$, then we deduce the vanishing of $\Ext^1_R(M,M)$ and $\Ext^2_R(M,M)$ from~(\ref{tp3}.1).
As $\pd_R(M)\leq1$, we conclude that $M$ is free. 
\end{proof}

The following corollary corresponds to a special case of Proposition~\ref{tp3}. It provides conditions under which the regularity assumption on $R$ in~\ref{HWL1}(ii) can be replaced by the complete intersection property.

\begin{cor} \label{tp3c1} Let $R$ be a $d$-dimensional local complete intersection and let $M$ be an $R$-module. Assume $M$ is locally free on the punctured spectrum of $R$.
\begin{enumerate}[\rm(i)]
\item If $\hh^i_{\fm}(M\otimes M^{\ast})=0$ for some odd integer $i$ such that $1\leq i \leq \depth_R(M)$, then $\pd_R(M)<\infty$. 
\item If $d\geq 4$, $\depth_R(M)\geq 3$, and $\hh^3_{\fm}(M\otimes_RM^*)=0$, then $\pd_R(M)\leq 1$. 
\item If $d\geq 4$, $\depth_R(M)\geq 3$, and $\hh^2_{\fm}(M\otimes_RM^*)=0=\hh^3_{\fm}(M\otimes_RM^*)=0$, then $M$ is free.
\end{enumerate}
\end{cor}

\begin{rmk} We briefly discuss that the conclusion of Corollary \ref{tp3c1}(i) is sharp.
\begin{enumerate}[\rm(i)]
\item It is necessary that the ring in Corollary~\ref{tp3c1}(i) is, at least, Gorenstein. For example, if we pick $R=\mathbb{C}[\![x, y, z, u, v]\!]/(yv-zu, yu-xv, xz-y^2)$, then 
$R$ is a three-dimensional Cohen-Macaulay ring such that $R$ is not Gorenstein but $R_{\fp}$ is regular for all $\fp \in \Spec(R)-\{\fm\}$. Therefore, if $M=\omega_R$, then $\pd_R(M)=\infty$, $M$ is locally free on the punctured spectrum of $R$, and $\hh^{1}_{\fm}(M\otimes _R M^\ast)\neq 0 \neq \hh^{3}_{\fm}(M\otimes _R M^\ast)$ and $\hh^{2}_{\fm}(M\otimes _R M^\ast)=0$; see Remark~\ref{directsum2}(ii). 
\item The locally freeness assumption on the module $M$ is needed for Corollary \ref{tp3c1}(i): If $R$ and $M$ are as in Example~\ref{directsum}, then $R$ is a complete intersection ring and $M$ is a non-free $R$-module such that $\hh^1_{\fm}(M\otimes M^{\ast})=0$ and $M$ is not locally free on the punctured spectrum of $R$.
\item The vanishing of $\hh^n_{\fm}(M\otimes M^{\ast})=0$ for some odd integer $n$ is necessary in Corollary~\ref{tp3c1}(i), even if $1\leq n \leq \depth_R(M)$:
If $R$ and $M$ are as in Example~\ref{TRE}, then $R$ is a complete intersection ring and $M$ is an $R$-module such that $\depth_R(M)=3$, $M$ is locally free on the punctured spectrum of $R$, $\hh^2_{\fm}(M\otimes M^{\ast})=0$, and $\pd_R(M)=\infty$.
\end{enumerate}
\end{rmk}

Let $X$ be a projective variety over an algebraically closed field $k$ and let $\mathcal{E}$ be a vector bundle. Set $\Gamma_\ast(\mathcal{E})=\bigoplus_{i\in\mathbb{Z}}\Gamma \big(X, \mathcal{E} (i)\big)$. We say $ \mathcal{E}$ is \emph{arithmetically Cohen-Macaulay} if $\hh^i\left(X, \mathcal{E}(j)\right)=0$ for all $j\in \mathbb{Z}$ and all $i = 1, \ldots , \dim(X)-1$. In passing we use of Corollary~\ref{tp3c1}(i) and prove the next result for arithmetically Cohen-Macaulay vector bundles over smooth complete intersections; cf.~\cite[1.2]{KC} and~\cite[1.5]{Da3}.
  
\begin{prop}\label{artCM} Let $k$ be an algebraically closed field, $X\subseteq \mathbb{P}^m_k$ be a globally complete intersection of dimension $d\geq3$, and let $\mathcal{E}$ be a vector bundle. Assume $\hh^{n}\big(X, \mathcal{E}\otimes\mathcal{E}^*(j)\big)=0$ for all $j\in\ZZ$ and for some integer $n$ such that $1\leq n\leq d-1$. 
\begin{enumerate}[\rm(i)]
\item If $n$ is even and $\depth\left(\Gamma_\ast(\mathcal{E})\right)\geq n+1$, then $\Gamma_\ast(\mathcal{E})$ has finite projective dimension over its affine cone. If, in addition, $\mathcal{E}$ is arithmetically Cohen-Macaulay, then it is a direct sum of powers of $\mathcal{O}(1)$.
\item If $d$ is odd and $\mathcal{E}$ is arithmetically Cohen-Macaulay, then $\mathcal{E}$ is a direct sum of powers of $\mathcal{O}(1)$.
\end{enumerate} 
\end{prop}

\begin{proof} Note that there is a homogeneous ideal $I$ of $k[x_0,\ldots, x_m]$ generated by a regular sequence such that $R=k[x_0, \ldots, x_m]/I$ and $X= \textmd{Proj}(R)$. Note that $\bigoplus_{n\in \mathbb{Z}} \hh^0\big(X,\mathcal{O}_X(n)\big) \cong R$, $\dim(R)=d+1$, and there exists an $R$-module $M$ such that $\mathcal{E} =\widetilde{M}$. In fact, $M=\Gamma_\ast(\mathcal{E})$, which is graded and reflexive. 

Arithmetically Cohen-Macaulay vector bundles correspond to maximal Cohen-Macaulay modules over the associated graded rings. Therefore, $M_{\fm}$ is maximal Cohen-Macaulay if and only if $\mathcal{E}$ is arithmetically Cohen-Macaulay, where $\fm$ denotes the graded maximal ideal of $R$.
Also, we have the following:
\begin{equation}\label{art} \tag{\ref*{artCM}.1}
\hh^{n+1}_{\fm}(M\otimes_RM^*)\cong\bigoplus_{i\in\mathbb{Z}}\hh^n \big(X,\mathcal{\mathcal{E}}\otimes \mathcal{\mathcal{\mathcal{E}}}^\ast(i)\big)=0.
\end{equation}

(i): We conclude from Corollary~\ref{tp3c1} and \eqref{art} that $\pd_{R_{\fm}}(M_{\fm})<\infty$. So $\pd_R(M)<\infty$; see~\cite[1.5.15(e)]{BH}. 
 If, in addition, $\mathcal{E}$ is arithmetically Cohen-Macaulay,
 by Auslander-Buchsbaum formula, it is free. As $M$ is graded, there is a finite set $L\subseteq \mathbb{Z}$ such that $M=\bigoplus_{ \ell\in L}R(\ell)$. Consequently, the observation $\mathcal{E} =\widetilde{M}=\bigoplus_{ \ell\in L}\widetilde{R(\ell)}=\bigoplus_{ \ell\in L}\mathcal{O}(1)^{\ell}$ completes the proof.

(ii): Since $d:=\dim(X)$ is odd, it follows $R$ is of even dimension.  This allows us to apply ~\ref{CelDa}. Also, note that $M=\Gamma_\ast(\mathcal{E})$ is maximal Cohen-Macaulay and $\mathcal{E} =\widetilde{M}$. In view of \eqref{art} and~\cite[3.9]{CeD} (see also~\ref{CelDa}), we have that $\pd_{R_{\fm}}(M_{\fm})=0$. Hence, $\pd_R(M)=0$ by~\cite[1.5.15(e)]{BH}. It remains to apply the argument of part~(i).
\end{proof}

We give an example and show that the conclusion of Proposition~\ref{artCM} is sharp.

\begin{eg}$\phantom{}$
\begin{enumerate}[\rm(i)]
\item Let $R=k[x,y,u,v]/(xy-uv)$ and let $I=(x,u)$. Let $X=\textmd{Proj}(R)$ and $\mathcal{E}=\widetilde{I}$. By Example~\ref{TRE}, we have
$$\hh^1_\ast\big(X,\mathcal{E}nd_{\mathcal{O}_X}(\mathcal{\mathcal{E}})\big)=\bigoplus_{i\in\mathbb{Z}}\hh^1\big(X,\mathcal{E}nd_{\mathcal{O}_X}(\mathcal{\mathcal{E}})(i)\big)=\bigoplus_{i\in\mathbb{Z}}\hh^1\big(X,\mathcal{\mathcal{E}}\otimes \mathcal{\mathcal{\mathcal{E}}}^\ast(i)\big) =\hh^2_{\fm}(I\otimes_RI^*)=0.$$ However, $\mathcal{E}$ is not trivial. The point is that $\dim(X)=\dim(R)-1=2$, which is strictly less than $3$.\smallskip
\item Let $X\subseteq \mathbb{P}^5_k$ be a smooth hypersurface of degree $d \geq 2$ and let $\mathcal{E}$ be an indecomposable arithmetically Cohen-Macaulay vector bundle of rank two. Then, by~\cite[1.1.(2)]{Kum} (see also~\cite[1.1.2]{Kem}), it follows that $\hh^1_{\ast}(X,\mathcal{E}nd_{\mathcal{O}_X}(\mathcal{\mathcal{E}}))=0$. As $\mathcal{E}$ is indecomposable, this example shows that $n$ needs to be an even integer in Proposition~\ref{artCM}(i).
\end{enumerate}
\end{eg}

%%%%%%%%%%%%%%%%%%%%%%%%%%%%%
\subsection*{A generalization of the third part of \ref{HWL1}}
%%%%%%%%%%%%%%%%%%%%%%%%%%%%%

Our next aim is to utilize Lemma~\ref{grz} and Proposition~\ref{l2} to extend the result of Huneke-Wiegand as stated in~\ref{HWL3}(iii) to a more general setting, leading to Corollary~\ref{HWL3cor}. First, we need to make some preparations:

\begin{lem} \label{star} Let $\mz$ be a specialization-closed subset of $\Spec(R)$, and let $M$ and $N$ be $R$-modules. If $\NF(N)\cup\NF(M)\subseteq\mz$, then $\hh_{\mz}^i(M^\ast\otimes_RN^\ast) \cong \hh_{\mz}^i\big((M\otimes_RN)^\ast\big)$ for all $i\geq 2$.
\end{lem}

\begin{proof} We consider the natural map $\psi:M^{\ast}\otimes_RN^{\ast}\to(M\otimes_RN)^{\ast}$ given by $\psi(f\otimes g)(x\otimes y)=f(x)g(y)$ for $f\in M^{\ast}$, $g\in N^{\ast}$, $x\in M$ and $y\in N$. Note that, given $\fp \in \Spec(R)$, if $M_{\fp}$ or $N_{\fp}$ is a free $R_{\fp}$-module, then $\psi_{\fp}$ is an isomorphism.

Set $K=\ker(\psi)$ and $C=\coker(\psi)$. Then $\Supp_R(C) \cup \Supp_R(K) \subseteq \NF(N)\cup\NF(M)\subseteq\mz$. So, $C$ and $K$ are both torsion with respect to $\mz$; see \ref{LC}(iv). Consequently, $\hh_{\mz}^i(C)=\hh_{\mz}^i(K)=0$ for all $i\geq 1$; see~\ref{LC}(v). We now split the exact sequence $0 \to C \to M^{\ast} \otimes_RN^{\ast} \overset{\psi}{\longrightarrow} (M\otimes_RN)^{\ast} \to K \to 0$ into two short exact sequences and deduce by using~\ref{LC}(ii) that $\hh_{\mz}^i(M^\ast\otimes_RN^\ast)\cong \hh_{\mz}^i\big((M\otimes_RN)^\ast\big)$ for all $i\geq 2$. 	
\end{proof}

Local duality is a classical result of Grothendieck, but we state it here for the convenience of the reader; see, for example,~\cite[3.5.8]{BH}. Recall that $(-)^{\vee}=\Hom_R(-,E_R(k))$, where $E_R(k)$ denotes the injective hull of the residue field $k$ of the local ring $R$.

\begin{chunk} \textbf{Local duality}. \label{Ld} If $R$ is a $d$-dimensional Cohen-Macaulay local ring with a canonical module $\omega_R$ and $M$ is an $R$-module, then $\Ext^{i}_R(M,\omega_R) \cong \hh^{d-i}_{\fm}(M)^{\vee}$ for all $i\in \ZZ$.
\end{chunk}

\begin{chunk}\textbf{Grade of a pair of modules}.  \label{chunk:bildikk}
Let $R$ be a local ring and let $M$ and $N$ be $R$-modules. The \emph{grade} of the pair $(M,N)$ is defined as 
\begin{equation}\label{eq:bildikk1}\tag{\ref*{chunk:bildikk}.1}
\gr_R(M,N)=\inf\{i: \Ext^i_R(M,N)\neq 0\}.
\end{equation}
It follows from~\cite[4.5]{AuBr} that 
\begin{equation}\label{eq:bildikk2}\tag{\ref*{chunk:bildikk}.2}
\gr_R(M,N)=\inf\big \{\depth_{R_{\fp}}(N_\fp) : \fp\in\Supp_R(M)\big\}.
\end{equation}
For simplicity, we denote  the grade of the pair $(M,R)$ by $\gr_R(M)$.
\end{chunk}

\begin{thm}\label{tomax} Let $R$ be a $d$-dimensional Gorenstein local ring, $\mz$ be a specialization-closed subset of $\Spec(R)$, and let $M$ and $N$ be $R$-modules. Assume $\NF(N)\cup\NF(M)\subseteq\mz$ and $\gr(\mz,R)\geq 3$. If $i$ is an integer such that $2\leq i\leq \gr(\mz,R)-1$, then
$\hh_{\mz}^i(M^\ast\otimes_RN^\ast) \cong \hh^{d-i+1}_{\fm} (M\otimes_RN)^{\vee}$. 
\end{thm}

\begin{proof} We set $L=M\otimes_RN$ and consider the natural map $\gamma:L\to L^{**}$. We also set $K=\ker(\gamma)$, $C=\coker(\gamma)$, and $g=\gr_R(\mz,R)$.

As $\NF(N)\cup\NF(M)\subseteq\mz$, it follows that $\NF(L)\subseteq\mz$. Hence, if $\fp \in \Spec(R)-\mz$, then $L_{\fp}$ is free so that $L_{\fp} \cong  L_{\fp}^{**}$, or equivalently, $\gamma_{\fp}$ is an isomorphism; see~\cite[1.1.9(b)]{Gdimbook}. Thus $\Supp_R(K)\subseteq \mz$ and $\Supp_R(C)\subseteq \mz$. This implies, in view of Lemma \ref{grz} and \eqref{eq:bildikk2}, that $g\leq \gr_R(K)$ and $g\leq \gr_R(C)$. We conclude that $\Ext^i_R(K,R)=0=\Ext^i_R(C,R)$ for all $i=0, \ldots, g-1$. Hence, the exact sequence $0 \to K \to L \overset{\gamma}{\rightarrow} L^{\ast\ast} \to C \to 0$ induces the following isomorphism
\begin{equation}\label{to}\tag{\ref*{tomax}.1}
\Ext^i_R(L,R)\cong\Ext^i_R(L^{**},R) \text{ for  all } i \text{ with } i\leq g-2.
\end{equation}

Then, if $i$ is an integer with $2\leq i\leq g-1$, it follows that:
\begin{align} \notag{}
\hh_{\mz}^i(M^\ast\otimes_RN^\ast) \cong \hh_{\mz}^i\big((M\otimes_RN)^\ast\big)=\hh_{\mz}^i(L^\ast)  \cong \Ext^{i+1}_R(\Tr_R L^*,R) & \cong \Ext^{i-1}_R(L^{**},R)  \notag{} \\ & \cong \Ext^{i-1}_R(L,R) \notag{} \\ & \cong \hh^{d-i+1}_{\fm}(L)^{\vee} \tag{\ref*{tomax}.2} \\ & = \hh^{d-i+1}_{\fm} \notag{} (M\otimes_RN)^{\vee}.
\end{align}
Here, in (\ref*{tomax}.2), the first isomorphism holds by Lemma \ref{star}, the third one holds by the fact $L^{**}\cong \Omega_R^2\Tr_R L^*$, the fourth one is due to \eqref{to}, the fifth one holds by \ref{Ld}. Moreover, we use Proposition \ref{l2}(i) and observe $\hh^{i}_{\mz}(L^{\ast}\otimes_RR)\cong\Ext^{i+1}_R(\Tr_R L^{\ast}, R)$, that is, the second isomorphism in~(\ref*{tomax}.2) holds.
\end{proof}

Next is a consequence of Theorem~\ref{tomax} which generalizes~\ref{HWL1}(iii):

\begin{cor} \label{HWL3cor} Let $R$ be a $d$-dimensional Gorenstein local ring, with $d\geq 3$, and let $M$ and $N$ be $R$-modules which are both locally free on the punctured spectrum of $R$. If $i$ is an integer such that $2 \leq i \leq d-1$, then $\hh_{\fm}^i(M^\ast\otimes_RN^\ast)\cong \big(\hh_{\fm}^{d-i+1}(M \otimes_RN )\big)^{\vee}$.
\end{cor}

\begin{cor} \label{depthcor} Let $R$ be a $d$-dimensional local complete intersection, with $d\geq 4$, and let $M$ be an $R$-module such that $\depth_R(M)\geq 3$ and $M$ is locally free on the punctured spectrum of $R$. 
\begin{enumerate}[\rm(i)]
\item If $\hh^{d-1}_{\fm}(M\otimes_RM^{\ast})=0$, then $\hh^{2}_{\fm}(M\otimes_RM^{\ast})=0$.
\item If $\hh^{d-2}_{\fm}(M\otimes_RM^{\ast})=0$, then $\hh^{3}_{\fm}(M\otimes_RM^{\ast})=0$.
\item If $\hh^{d-1}_{\fm}(M\otimes_RM^{\ast})=\hh^{d-2}_{\fm}(M\otimes_RM^{\ast})=0$, then $M$ is free.
\end{enumerate}
\end{cor}

\begin{proof} Part (iii) follows from parts~(i) and~(ii) due to Corollary~\ref{tp3c1}(iii). As $\depth_R(M)\geq 3$ and $M$ is locally free on the punctured spectrum of $R$, it follows that $M\cong M^{\ast\ast}$; see Lemma~\ref{remzi}. We use Corollary~\ref{HWL3cor} with $N=M^{\ast}$ and obtain the isomorphism $\hh_{\fm}^{2}(M^\ast\otimes_RM) \cong \hh_{\fm}^{2}(M^\ast\otimes_RM^{\ast\ast})  \cong \big(\hh_{\fm}^{d-1}(M \otimes_RM^{\ast} )\big)^{\vee}$. Similarly, we see that $ \hh_{\fm}^{3}(M^\ast\otimes_RM) \cong \hh_{\fm}^{3}(M^\ast\otimes_RM^{\ast\ast})  \cong \big(\hh_{\fm}^{d-2}(M \otimes_RM^{\ast} )\big)^{\vee}$. This proves parts~(i) and~(ii). 
\end{proof}

We do not know whether or not the depth assumption on $M$ in Corollary~\ref{depthcor} is necessary. However, in some cases, a similar result works for modules of depth at least two. The following proposition is a translation of a result of {C}esnavi\v{c}ius~\cite{KC} into the setting of commutative algebra:

\begin{prop}\label{t3.1} Let $R$ be a $d$-dimensional graded normal complete intersection over a field, with $d\geq 4$, and let $M$ be a graded $R$-module such that $\depth_R(M)\geq 2$ and $M$ is locally free on the punctured spectrum of $R$. Then $M$ is free if at least one of the following conditions holds:
\begin{enumerate}[\rm(i)]
\item $\hh^2_{\fm}(M\otimes_RM^{\ast})=\hh^3_{\fm}(M\otimes_RM^{\ast})=0$.
\item $\hh^{d-2}_{\fm}(M\otimes_RM^{\ast})=\hh^{d-1}_{\fm}(M\otimes_RM^{\ast})=0$.
\end{enumerate}
\end{prop}

\begin{proof} Note that $\displaystyle{R=\frac{k[x_0,\ldots, x_n]}{I}}$ for some polynomial ring $k[x_0,\ldots, x_n]$ and homogeneous ideal $I$. Set $X= \textmd{Proj}(R)$ and $\mathcal{E}=\widetilde{M}$. Then $X\subseteq \mathbb{P}^n$ is globally a complete intersection, $\mathcal{E}$ is  a vector bundle, and $\bigoplus_{n\in \mathbb{Z}} \hh^0\big(X,\mathcal{O}_X(n)\big)\cong R$. For part (i), we have:
$$\hh^1_\ast\big(X,\mathcal{E}nd_{\mathcal{O}_X}(\mathcal{\mathcal{E}})\big)=\bigoplus_{i\in\mathbb{Z}}\hh^1\big(X,\mathcal{E}nd_{\mathcal{O}_X}(\mathcal{\mathcal{E}})(i)\big)=\bigoplus_{i\in\mathbb{Z}}\hh^1\big(X,\mathcal{\mathcal{E}}\otimes \mathcal{\mathcal{\mathcal{E}}}^\ast(i)\big) =\hh^2_{\fm}(M\otimes_RM^*)=0.$$ 

Similarly, it follows that $\hh^2_\ast\big(X,\mathcal{E}nd_{\mathcal{O}_X}(\mathcal{\mathcal{E}})\big)=0$. Thus, by~\cite[1.2]{KC}, $\mathcal{E}$ is a direct sum of powers of $\mathcal{O}(1)$. So, $M= \hh^0_\ast(X, \mathcal{E})$ is free. In view of Theorem~\ref{tomax}, part~(ii) follows from part~(i).
\end{proof}

We end this section with another application of Theorem~\ref{tomax}. To the best of our knowledge, the result we obtain next is new, even for regular local rings.

\begin{prop}\label{cortomax} Let $R$ be a $d$-dimensional Gorenstein local ring, with $d\geq3$, and let $M$ be an $R$-module. Assume the following conditions hold:
\begin{enumerate}[\rm(i)]
\item $d$ is odd.
\item $\depth_R(M)=(d+1)/2$.
\item $M$ is Tor-rigid (for example, $R$ is regular).
\item $M$ is locally free on the punctured spectrum of $R$.
\end{enumerate}
Then $\hh^i_{\fm}(M\otimes_RM^*)\neq 0$ for all $i=1, \ldots, d-1$.  
\end{prop}

\begin{proof} Suppose $\hh^i_{\fm}(M\otimes_RM^*)=0$ for some $i$ with $1\leq i \leq d-1$ and seek a contradiction. Set $t=(d+1)/2$. If $i=1$, then Corollary~\ref{c1} shows that $M$ is free so that $\depth_R(M)=t=d=1$. So, we may assume $2\leq i \leq d-1$. Note that $t\leq d-1$. 

We claim that, if $j$ is an integer such that $2\leq j\leq t  \leq d-1$, then $\hh^j_{\fm}(M\otimes_RM^*)\neq 0$. Indeed, if $j$ is such an integer and $\hh^j_{\fm}(M\otimes_RM^*)=0$, then 
Corollary~\ref{remzidenoturu}(ii) implies that $\pd_R(M)<j-1$. This yields the contradiction $d<2t-1=d$ since $d-t=\depth(R)-\depth_R(M)=\pd_R(M)<j-1\leq t-1$. So the claim holds. Now, by the claim, we may assume that $i\geq t+1$. 

Note that $M^{\ast}$ is locally free on the punctured spectrum of $R$. Moreover, $M$ satisfies $(S_2)$ since $\depth_R(M)=t\geq 2$ and $M$ is locally free on the punctured spectrum of $R$. Therefore, $M$ is reflexive; see~\ref{SL}. Now we use Corollary~\ref{HWL3cor} with $N=M^{\ast}$ and conclude that $0=\hh_{\fm}^i(M^\ast\otimes_RM)\cong \big(\hh_{\fm}^{d-i+1}(M \otimes_RM^{\ast})\big)^{\vee}$. Thus $\hh_{\fm}^{d-i+1}(M \otimes_RM^{\ast})=0$. As $i\geq t+1$, it follows that $d-i+1\leq d-t-1+1=(d-1)/2<(d+1)/2=t$. Consequently, by setting, $j=d-i+1$, we obtain $\hh_{\fm}^{j}(M \otimes_RM^{\ast})=0$ where $2\leq j\leq t  \leq d-1$. However, such an integer $j$ cannot exist due to the claim we proved. This gives the contradiction we seek.
\end{proof}

Recall that modules of finite projective dimension are Tor-rigid over certain hypersurface rings; see~\ref{TorRigid}(ii). Hence a special case of Proposition~\ref{cortomax} is the following:

\begin{cor} Let $R=k[\![x_1, \ldots, x_{d+1}]\!]/(f)$ be a hypersurface ring, with $d\geq 3$ is odd, and let $M$ be an $R$-module which is locally free on the punctured spectrum of $R$. If $\pd_R(M)=(d-1)/2$, then $\hh^i_{\fm}(M\otimes_RM^*)\neq 0$ for all $i=1, \ldots, d-1$.  
\end{cor}

We conclude this section with an example that highlights the necessity of the locally free assumption in Proposition~\ref{cortomax}.

\begin{eg}\label{notvb2} Let $R$ be a $3$-dimensional Gorenstein local ring and let $0\neq x\in \fm$ be a non zero-divisor on $R$. Set $M= R/xR$. Then $\depth_R(M)=2=(3+1)/2$ and $\pd_R(M)=1$. So, $M$ is Tor-rigid. Moreover, $M^{\ast}=0$ and hence $\hh^i_{\fm}(M\otimes_RM^*)=0$ for all $i$. Pick $\fp \in \Spec(R)-\{\fm\}$ such that $x \in \fp$. Then $\pd_{R_{\fp}}(M_{\fp})=1$ and hence $M$ is not locally free on the punctured spectrum of $R$.
\end{eg}

%%%%%%%%%%%%%%%%%%%%%%%%%%%%%%%%%%%%%%%%%%%%%%%%%%%%%%%%%%%%%%%%%
%%%%%%%%%%%%%%%%%%% SECTION  %%%%%%%%%%%%%%%%%%%%%%%%%%%%%%%%%%%%%%%%
%%%%%%%%%%%%%%%%%%% SECTION  %%%%%%%%%%%%%%%%%%%%%%%%%%%%%%%%%%%%%%%%
%%%%%%%%%%%%%%%%%%% SECTION  %%%%%%%%%%%%%%%%%%%%%%%%%%%%%%%%%%%%%%%%
%%%%%%%%%%%%%%%%%%%%%%%%%%%%%%%%%%%%%%%%%%%%%%%%%%%%%%%%%%%%%%%%%

\section{Further applications and modules of finite complete intersection dimension}\label{section6}

 The aim of this section is twofold. First, we prove Theorem~\ref{mainofsection} which corroborates the results of the previous section and establishes an additional freeness criteria for modules $M$ over Gorenstein rings in terms of the vanishing of $\hh^{i}_{\fm}(M\otimes_RM^*)$. In the second part, we give applications of Theorem~\ref{t12} concerning modules of finite complete intersection dimension and contribute to the literature; see~\ref{CI} and Proposition~\ref{cx1}. Our first result, Theorem~\ref{mainofsection}, is motivated by the following theorem of Huneke and Wiegand:

\begin{chunk} [{Huneke-Wiegand \cite[4.3]{HW1}}] \label{bildik} If $R$ is a $d$-dimensional regular local ring and $M$ is an $R$-module which is locally free on the punctured spectrum of $R$ such that $\hh^{d-1}_{\fm}(M\otimes_RM^*)=0$, then $M^{\ast}$ is free.
\end{chunk}

We give a proof of~\ref{bildik} for completeness by using our results from Section~\ref{section5}.

\begin{proof} Note that $\depth_R(M^{\ast})\geq \min\{2, d\}$; see~\cite[1.4.19]{BH}. Hence, we may assume $d\geq 3$. As $M^{\ast}$ is locally free on the punctured spectrum of $R$ and $\depth_R(M^{\ast})\geq 2$, we see that $M^{\ast}$ is reflexive; see~\ref{SL}(ii). We use Corollary~\ref{HWL3cor} with $N=M^{\ast}$ and $i=2$, and conclude that $\hh_{\fm}^{2}(M^\ast\otimes_RM^{\ast\ast})\cong \big(\hh_{\fm}^{d-1}(M \otimes_RM^{\ast} )\big)^{\vee}=0$. Recall that each $R$-module is Tor-rigid; see~\ref{TorRigid}. Thus, $M^{\ast}$ is free by Corollary~\ref{remzidenoturu}(ii).
\end{proof}
 
\begin{thm} \label{mainofsection} Let $R$ be a $d$-dimensional Gorenstein local ring, with $d\geq 1$, and let $M$ be a reflexive $R$-module which is locally free on the punctured spectrum of $R$. If $M$ is Tor-rigid and $\hh^{d-1}_{\fm}(M\otimes_RM^{\ast})=0$, then $M$ is free.
\end{thm}

We prove Theorem \ref{mainofsection} after some preparations. Let us mention here that the reflexivity assumption is necessary for the theorem. For example, if $R$ is a three-dimensional Cohen-Macaulay local ring, then $\fm$ is not a reflexive $R$-module, $\fm^{\ast}\cong R$, and $\hh^{2}_{\fm}(\fm\otimes_R\fm^{\ast})=0$; see, for example,~\cite[5.4]{MA}.

Theorem \ref{mainofsection} has several consequences. For example, it generalizes \ref{bildik} (recall that, if $M$ is as in \ref{bildik}, then $M^{\ast}$ is reflexive for the nontrivial case where $d\geq 2$.) Theorem~\ref{mainofsection} also gives useful information over rings which are not necessarily regular. For example, in view of~\ref{TorRigid}(iii), Theorem~\ref{mainofsection} yields the following:

\begin{cor} Let $R$ be a Gorenstein local ring of positive dimension and let $M=\Omega^i_R(R/I)$, where $i\geq 2$, $I=(A:_R\fm)$, and $A$ is an $\fm$-primary ideal of $R$ such that $A \neq \fm$. If $\hh^{d-1}_{\fm}(M\otimes_RM^*)=0$, then $R$ is regular.
\end{cor}

The next lemma is key to establishing Theorem \ref{mainofsection} and should be of independent interest. Recall that the grade of the pair $(M,N)$ is $\gr_R(M,N)=\inf\{i\geq 0: \Ext^i_R(M,N)\neq 0\}$; see~\ref{chunk:bildikk}.

\begin{lem}\label{duality} Let $M$, $N$, and $K$ be $R$-modules (not necessarily finitely generated) and let $\ell \geq 1$ be an integer. Assume $\Ext^j_R(N,K)=0$ for all $j=1, \ldots, \ell-1$. Assume further $\gr_R\big(\Tor_i^R(M,N),K\big)>\ell-i$ for all $i=1, \ldots, \ell$. Then the following hold:
\begin{enumerate}[\rm(i)]
\item $\Ext^i_R\big(M,\Hom_R(N,K)\big)\cong\Ext^i_R\big(M\otimes_RN,K\big)$ for all $i=1, \ldots, \ell-1$. 
\item There is an injection $\Ext^\ell_R\big(M,\Hom_R(N,K)\big)\hookrightarrow\Ext^\ell_R\big(M\otimes_RN,K\big)$.
\end{enumerate}	
\end{lem}

\begin{proof} There are two spectral sequences converging to the same point~\cite[10.62 and 10.64]{Roit}:
$$E_2^{pq}=\Ext^p_R(\Tor_q^R(M,N),K) \Longrightarrow H^{p+q} \text{ and }
F_2^{pq}=\Ext^p_R(M,\Ext^q_R(N,K)) \Longrightarrow H^{p+q}.$$
We have that $E_2^{pq}=0$ if $1\leq q\leq \ell$ and $p\leq \ell-q$ since $\gr_R\big(\Tor_i^R(M,N),K\big)>\ell-i$ for all $i=1, \ldots, \ell$. Therefore, we obtain an isomorphism $E_2^{i0}\cong H^i$ for all integers $i\le \ell$. Also, $F_2^{pq}=0$ if $1\le q\le \ell-1$ since $\Ext^i_R(N,K)=0$ for all $1\le i\le \ell-1$. So, there is an isomorphism $F_2^{i0}\cong H^i$ for all integers $i\le \ell-1$, and an injection $F_2^{\ell0}\hookrightarrow H^\ell$. This yields $\Ext^i_R\big(M\otimes_R N,K\big)=E_2^{i0}\cong H^i\cong F_2^{i0}=\Ext^i_R\big(M,\Hom_R(N,K)\big)$ for all integers $i\le \ell-1$, and an injection $\Ext^\ell_R\big(M,\Hom_R(N,K)\big)\cong F_2^{\ell0}\hookrightarrow H^\ell\cong\Ext^\ell_R\big(M\otimes_RN,K\big)$.
\end{proof}

Theorem~\ref{mainofsection} is a special case of what we prove next. Recall that an $R$-module $M$ satisfies Serre's condition $(S_{v})$ for some $v\geq 0$ if $\depth_{R_{\fp}}(M_\fp)\geq\min\big\{v,\Ht_R(\fp)\big\}$ for all $\fp\in\Supp_R(M)$; see~\ref{SL}.

\begin{prop}\label{ll} Let $R$ be a Gorenstein local ring, $M$ be an $R$-module, and let $n\geq 1$ be an integer. Assume $M$ satisfies $(S_{n+1})$ and $M_{\fp}$ is free for all $\fp \in \Spec(R)$ with $\Ht_R(\fp)\leq n-1$. Then the following hold:
\begin{enumerate}[\rm(i)]
\item $\Ext^i_R(M,M)\cong\Ext^i_R(M\otimes_RM^*,R)$ for all $i=1, \ldots, n-1$.
\item There is an injection $\Ext^n_{R}(M,M)\hookrightarrow\Ext^n_R(M\otimes_RM^*,R)$.
\item If $M$ is Tor-rigid and $\hh^{d-n}_{\fm}(M\otimes_RM^*)=0$, then $\pd_R(M)<n$.
\end{enumerate}
\end{prop}

\begin{proof} It follows from~\cite[4.25]{AuBr} that $\Ext^i_R(\Tr_R M,R)=0$ for all $i=1, \ldots, n+1$. Consequently, $M$ is reflexive and $\Ext^i_R(M^*,R)=0$ for all $i=1, \ldots, n-1$ since $M^{\ast} \cong \Omega^2_R \Tr_R M$; see~\ref{a1} and~(\ref{lan}.2). 

Fix $i\geq 1$ and let $\fp \in \Supp_R\big(\Tor_i^R(M,M^*)\big)$. Then $\Ht_R(\fp)\geq n$ by our assumption. This implies that $\gr_R\big(\Tor_i^R(M,M^*),R\big)=\gr_R\big(\Tor_i^R(M,M^*)\big)\geq n$ for all $i\geq 1$; see (\ref{chunk:bildikk}.2). Now we use Lemma~\ref{duality} with $N=M^{\ast}$, $K=R$, and $n=\ell$, and obtain the claims of parts~(i) and~(ii).

Next assume $\hh^{d-n}_{\fm}(M\otimes_RM^*)=0$. Then it follows from~\ref{Ld} that $\Ext^n_R(M\otimes_RM^{\ast},R)=0$. Thus, by part (ii), we have that  $\Ext^n_R(M,M)=0$. Therefore, if $M$ is Tor-rigid, then $\pd_R(M)<n$; see~\ref{TorRigid}(i)
\end{proof}

Next, we give a proof of Theorem~\ref{mainofsection}:

\begin{proof} [Proof of Theorem \ref{mainofsection}] As $R$ is Gorenstein and $M$ is reflexive, we know that $M$ satisfies $(S_2)$; see \ref{SL}(ii). Moreover, $M_{\fp}$ is free over $R_{\fp}$ for all minimal prime ideals $\fp$ of $R$ because $d \geq 1$ and $M$ is locally free on the punctured spectrum of $R$. Therefore, since $M$ is Tor-rigid, we make use of Proposition~\ref{ll}(iii) with $n=1$ and conclude that $M$ is free. 
\end{proof}

%%%%%%%%%%%%%%%%%%%%%%%%%%%%%%%%%%%%%%%%%%%%%%%%%%%%%%%%%%%%%%%%%
%%%%%%%%%%%%%%%%%%% SECTION  %%%%%%%%%%%%%%%%%%%%%%%%%%%%%%%%%%%%%%%%
%%%%%%%%%%%%%%%%%%% SECTION  %%%%%%%%%%%%%%%%%%%%%%%%%%%%%%%%%%%%%%%%
%%%%%%%%%%%%%%%%%%% SECTION  %%%%%%%%%%%%%%%%%%%%%%%%%%%%%%%%%%%%%%%%
%%%%%%%%%%%%%%%%%%%%%%%%%%%%%%%%%%%%%%%%%%%%%%%%%%%%%%%%%%%%%%%%%

\subsection*{On modules of finite complete intersection dimension} In this subsection we provide applications of Theorem~\ref{t12} and Corollary~\ref{c1} for modules of finite complete intersection dimension over local rings. We recall the definition and give some remarks:

\begin{chunk} [{Avramov-Gasharov-Peeva \cite{AGP}}] \label{CI} A \emph{quasi-deformation} of a local ring $R$ is a diagram $R\rightarrow A\twoheadleftarrow Q$ of local ring homomorphisms such that $R\rightarrow A$ is flat and $A\twoheadleftarrow Q$ is surjective with its kernel is generated by a regular sequence. The \emph{complete intersection dimension} of an $R$-module $M$ is defined as: $$\CI_R(M)=\inf\{\pd_Q(M\otimes_RA)-\pd_Q(A)\mid R\rightarrow A\twoheadleftarrow Q \text{ is a quasi-deformation}\}.$$
Therefore, $M$ has finite complete intersection dimension if there exists a quasi-deformation $R \rightarrow A \twoheadleftarrow Q$ for which $\pd_Q(M\otimes_RA)$ is finite. 
To facilitate the discussion, we note a few key properties; see~\cite[1.3 and 1.4]{AGP}.
\begin{enumerate}[\rm(i)] 
\item Given an $R$-module $M$, it follows that $\Gdim_R(M) \leq \CI_R(M) \leq \pd_R(M)$. 
\item $R$ is a complete intersection if and only if $\CI_R(N)<\infty$ for each $R$-module $N$.
\item If $\CI_R(M)<\infty$, then $\Gdim_R(M)=\CI_R(M)=\depth(R)-\depth_{R}(M)$.
\end{enumerate}
\end{chunk}	

\begin{chunk} [{Celikbas-Dao \cite[3.9]{CeD}}] \label{CelDa} Let $R$ be a $d$-dimensional Gorenstein local ring and let $M$ be a maximal Cohen-Macaulay $R$-module such that $\CI_R(M)<\infty$ and $M$ is locally free on the punctured spectrum of $R$. Assume $d\geq 2$ is even. If $\depth_R(M\otimes_RM^{\ast})\geq 1$, or equivalently, $\hh^0_{\fm}(M\otimes_RM^{\ast})=0$, then $M$ is free.
\end{chunk}

We know by Corollary~\ref{c1} that, if $R$ is Cohen-Macaulay and $M$ is an $R$-module such that $\depth_R(M)\geq 1$ and $M$ is locally free on the punctured spectrum of $R$, then the vanishing of $\hh^1_{\fm}(M\otimes_RM^{\ast})$ forces $M$ to be free. Our aim is to study the vanishing of $\hh^i_{\fm}(M\otimes_RM^{\ast})$ for $i\geq 2$ for modules $M$ with $\CI_R(M)<\infty$, and prove an extension of the fact stated in~\ref{CelDa}. We establish this in Proposition~\ref{gor}. In the proof of the proposition, we use several of our previous results, including Proposition~\ref{l2} and Theorem~\ref{tomax}. We also need several facts, some of which have been previously recorded for complete intersections. These facts actually hold over arbitrary local rings for modules of finite complete intersection dimension, so we present them again for the convenience of the reader.

\begin{chunk} \label{CIfacts} Let $R$ be a local ring and let $M$ and $N$ be $R$-modules. Assume $\CI_R(M)<\infty$.
\begin{enumerate}[\rm(i)]
\item If $\Ext^{2i}_R(M,M)=0$ for some $i\geq 0$, then $\pd_R(M)<\infty$; see~\ref{JorJot} and~\cite[4.2]{AvBu}. 
\item If $\mz$ is a specialization-closed subset of $\Spec(R)$ and $\Supp_R\big(\Tor_{i}^R(M,N)\big)\subseteq\mz$ for all $i\gg 0$, then $ \Supp_R\Big(\widehat{\Tor}_j^R(M,N)\Big) \subseteq\mz$ for all $j\in\ZZ$; see~\ref{fp}(i) and~\cite[4.9]{AvBu}.
\item If $\cx_R(M)\leq 1$, then $\Omega_R^t(M) \cong \Omega_R^{t+2}(M)$, where $t= \CI_R(M)$; see~\cite[7.3]{AGP}.
\end{enumerate}
\end{chunk}

\begin{prop}\label{gor} Let $R$ be a $d$-dimensional Cohen-Macaulay local ring, with $d\geq 2$, and let $M$ be a maximal Cohen-Macaulay $R$-module such that $M$ is locally free on the punctured spectrum of $R$ and $\CI_R(M)<\infty$. Assume $\hh^j_{\fm}(M\otimes_RM^{\ast})=0$ for some integer $j$.
Then $M$ is free if at least one of the following holds:
\begin{enumerate}[\rm(i)]
\item $1\leq j\leq d$ and $j$ is odd.
\item $R$ is Gorenstein, $d$ is even, and $0\leq j\leq d-1$.
\end{enumerate}
\end{prop}

\begin{proof} For part (i), we may assume $j\geq 3$ due Corollary~\ref{c1}. In that case, Lemma~\ref{remzii} implies that $\Ext^{j-1}_R(M,M)=0$. As $j-1$ is even and $\CI_R(M)<\infty$, it follows that $\pd_R(M)<\infty$, that is, $M$ is free; see~\ref{CIfacts}(i).

Next we proceed to establish part (ii). In view of~\ref{CelDa} and part (i), we may assume $2\leq j\leq d-1$ and $j$ is even. Note that $M$ satisfies $(S_2)$, that is, $M$ is reflexive, since $\depth_R(M)\geq 2$ and $M$ is locally free on the punctured spectrum of $R$; see~\ref{SL}. Hence, we conclude from Theorem~\ref{tomax} (by setting $N=M^{\ast}$) that $\hh^{d-j+1}_{\fm}(M\otimes_RM^{\ast})=0$. This yields
$$0=\hh^{d-j+1}_{\fm}(M^{\ast}\otimes_RM)\cong \Ext^{d-j+2}_R(\Tr_RM^{\ast},M)\cong \Ext^{d-j}_R(\Omega^2_R \Tr_RM^{\ast},M) \cong \Ext^{d-j}_R(M,M).$$
Here, we obtain the first isomorphism from Proposition~\ref{l2}(i) (by setting $n=d$ since $0\leq d-j+1\leq d-1$) and the third isomorphism from the fact that $M \cong M^{\ast\ast} \cong \Omega^2\Tr_R M^{\ast}$; see~\ref{a1}. As $d-j$ is even, it follows that $M$ is free as in part~(i).
\end{proof}

Our next result is mainly an application of Theorem~\ref{t12} for modules of bounded Betti numbers and of finite complete intersection dimension.

\begin{prop}\label{cx1} Let $R$ be a local ring and let $\mz$ be a specialization-closed subset of $\Spec(R)$. Assume that $M$ and $N$ are $R$-modules, and that $n\geq 0$ is an integer. Assume further:\begin{enumerate}[\rm(i)]
\item $\CI_R(M)<\infty$ and $\cx_R(M)\leq 1$ 
\item $\Supp_R\big(\Tor_i^R(M,N)\big)\subseteq\mz$ for all $i\geq 1$.  
\item $\depth_{R_\fp}(M_\fp)+\depth_{R_\fp}(N_\fp)-\depth(R_{\fp})\geq n$ for all $\fp\in\mz$.
\end{enumerate}
Then the following hold:
\begin{enumerate}[\rm(a)]
\item If $n\geq3$, then $\hh^i_{\mz}(M\otimes_RN)\cong\hh^{i+2}_{\mz}(M\otimes_RN)$ for all $i=0, \ldots, n-3$. 
\item If $\hh^n_{\mz}(M\otimes_RN)=0$, then $\widehat{\Tor}_{2i-n}^R(M,N)=0$ for all $i\in\ZZ$. Therefore, if $n$ is even (respectively, $n$ is odd), then $\widehat{\Tor}_{2j}^R(M,N)=0$ for all $j\in\ZZ$ (respectively, $\widehat{\Tor}_{2j+1}^R(M,N)=0$ for all $j\in\ZZ$).  
\end{enumerate}	
\end{prop}

\begin{proof} It follows from our assumption (ii) that $\Supp_R\big(\widehat{\Tor}_i^R(M,N)\big)\subseteq\mz$ for all $i\in\ZZ$; see~\ref{CIfacts}(ii). Moreover, we set $t=\CI_R(M)$ so that $\Omega_R^t(M) \cong \Omega^{t+2}_R(M)$; see~\ref{CIfacts}(iii). Thus, given $i\in \ZZ$, we have
\begin{equation}\label{cx1eq} \tag{\ref*{cx1}.1}
\widehat{\Tor}_i^R(M,N)\cong \widehat{\Tor}_{i-t}^R(\Omega_R^tM,N)\cong \widehat{\Tor}_{i-t}^R(\Omega_R^{t+2}M,N) \cong \widehat{\Tor}_{i+2}^R(M,N). 
\end{equation}
Assume $n\geq 3$ and $i$ is an integer such that $0\leq i \leq n-3$. Then the claim in part (a) holds due to
\begin{equation}\label{cx2eq} \tag{\ref*{cx1}.2}
\hh^i_{\mz}(M\otimes_RN)\cong \widehat{\Tor}_{-i}^R(M,N) \cong \widehat{\Tor}_{-i-2}^R(M,N) \cong \hh^{i+2}_{\mz}(M\otimes_RN).
\end{equation}
Here the first and the third isomorphisms follow from Theorem~\ref{t12}, and the second one is due to~\eqref{cx1eq}.

It follows from Theorem~\ref{t12} that there is an injection $\widehat{\Tor}_{-n}^R(M,N)\hookrightarrow\hh^n_{\mz}(M\otimes_RN)$.
Therefore, if $\hh^n_{\mz}(M\otimes_RN)=0$, then $\widehat{\Tor}_{-n}^R(M,N)=0$ so that \eqref{cx1eq} implies that $\widehat{\Tor}_{2i-n}^R(M,N)=0$ for all $i\in\ZZ$. This proves the claim of part (b).
\end{proof}

\begin{cor}\label{cq} Let $R$ be a $d$-dimensional local complete intersection ring, with $d\geq 1$, and let $M$ and $N$ be $R$-modules. Set $b=\depth_{R}(M)+\depth_{R}(N)-\depth(R)$, and assume:
\begin{enumerate}[\rm(i)]
\item $\cx_R(M)\leq 1$.
\item $M$ or $N$ is locally free on the punctured spectrum of $R$. 
\item $\hh^r_{\fm}(M\otimes_RN)=\hh^s_{\fm}(M\otimes_RN)=0$ for some integers $r$ and $s$ with $0\leq r, s \leq b$, $r$ is odd, and $s$ is even.
\end{enumerate}	
Then the following hold:
\begin{enumerate}[\rm(a)]
\item $\Tor_i^R(M,N)=0$ for all $i\gg0$.
\item If $R$ is a hypersurface, or $M=N$, then $\pd(M)<\infty$ or $\pd(N)<\infty$.
\item If $M$ and $N$ are maximal Cohen-Macaulay, so is $M\otimes_RN$.
\end{enumerate}	
\end{cor}

\begin{proof} As $\hh^r_{\fm}(M\otimes_RN)=\hh^s_{\fm}(M\otimes_RN)=0$, by setting $\mz=\V(\fm)$, we conclude from Proposition~\ref{cx1}(b) that $\widehat{\Tor}_{i}^R(M,N)=0$ for all $i\in\ZZ$. Thus, $\Tor_i^R(M,N)=0$ for all $i\gg0$; see~\ref{fp}(i). Hence, if $R$ is a hypersurface, or $M=N$, then $\pd(M)<\infty$ or $\pd(N)<\infty$; see~\ref{TVFPD}.

If $M$ is maximal Cohen-Macaulay, then the vanishing of $\Tor_i^R(M,N)=0$ for all $i\gg0$ forces the vanishing of $\Tor_i^R(M,N)$ for all $i\geq 1$; see~\cite[2.1]{Jo1}. Thus, if $M$ and $N$ are maximal Cohen-Macaulay, then the depth formula holds for the pair $(M,N)$, and this implies that $M\otimes_RN$ is maximal Cohen-Macaulay; see~\ref{df}.
\end{proof}

The following example shows that the key hypotheses of Corollary~\ref{cq} are necessary:

\begin{eg} $\phantom{}$
\begin{enumerate}[\rm(i)] 
\item Let $R=k[\![x,y,z,w]\!]/(xy)$. Then $\depth(R)=3$, and $\cx_R(T)\leq 1$ for all $R$-modules $T$ since $R$ is a hypersurface ring. Let $M=N=R/xR$. Since $\depth_R(M\otimes_RN)=\depth_R(M)=\depth(R)=3$, it follows that $\hh^1_{\fm}(M\otimes_RN)=0$. Note that $M$ is not locally free on the punctured spectrum of $R$, and $\widehat{\Tor}_{1}^R(M,N)\cong \Tor_1^R(M,N)\neq 0$. \smallskip
		
\item Let $R=k[\![x,y,z,w,t]\!]/(xy)$. Then $\depth(R)=4$, and $\cx_R(T)\leq 1$ for all $R$-modules $T$ since $R$ is a hypersurface ring. Let $M=N=\fm$. Then $\pd_R(M)=\infty$, $M$ is locally free on the punctured spectrum of $R$, and $b=\depth_{R}(M)+\depth_{R}(N)-\depth(R)<0$. Note that $ \hh^{2}_{\fm}(M\otimes_RN)=\hh^{3}_{\fm}(M\otimes_RN) =0$; see~\cite[5.5(i)]{MA}.
\end{enumerate}
\end{eg}

\begin{rmk} Note that Proposition~\ref{cx1} is broad and applies to rings that are not necessarily complete intersections. For example, if $R$ is a local ring such that $\cod(R)\leq 3$, or if $R$ is a Gorenstein local ring such that $\cod(R)\leq 4$, then each $R$-module that has finite complexity also has finite complete intersection dimension; see~\cite[1.5]{CAP}. Thus, modules of bounded Betti numbers over such rings enjoy the properties stated in part~(a) of  Proposition~\ref{cx1}. For a specific example, let us note that, if $I$ is an Ulrich ideal of a Gorenstein local ring $R$ such that $\cod(R)\leq 4$, then $\cx_R(R/I)\leq 1$ and $\CI_R(R/I)<\infty$; see~\cite[7.4]{GotoUlrich}.
\end{rmk}

Recall that if $R$ is a local hypersurface with an isolated singularity, and $M$ and $N$ are nonfree maximal Cohen-Macaulay $R$-modules, then the depth of $M\otimes_RN$ cannot exceed one; see~\ref{SRT}. We close this section by proving a related result: 

\begin{cor}\label{hhyper} Let $R$ be a $d$-dimensional local hypersurface, with $d\geq 2$, and let $M$ and $N$ be nonfree maximal Cohen-Macaulay $R$-modules. Assume $R$ has an isolated singularity. Then $\hh^d_{\fm}(M\otimes_RN)\neq 0$. Moreover, one of the following conditions holds:
\begin{enumerate}[\rm(i)]
\item $\hh^i_{\fm}(M\otimes_RN)\neq 0$ for all $i=0, \ldots, d-1$. 
\item $\hh^r_{\fm}(M\otimes_RN)\neq 0=\hh^s_{\fm}(M\otimes_RN)$ for all even integers $r$ and all odd integers $s$ with $0\leq r, s\leq d-1$.
\item $\hh^r_{\fm}(M\otimes_RN)\neq 0=\hh^s_{\fm}(M\otimes_RN)$ for all odd integers $r$ and all even integers $s$ with $0\leq r, s\leq d-1$.
\end{enumerate}	
\end{cor}

\begin{proof} Note that, since $d\geq 2$ and $R$ has an isolated singularity, it follows that $R$ is a domain; see, for example,~\cite[23.8]{Mat}. Moreover, $\Supp_R(M)=\Supp_R(N)=\Spec(R)$ as $M$ and $N$ are nonzero torsion-free $R$-modules. Hence, $d=\dim_R(M\otimes_RN)$ and so $\hh^d_{\fm}(M\otimes_RN)\neq 0$; see~\ref{LCo}(iii).

If $\hh^i_{\fm}(M\otimes_RN)=0$ for all $i=0, \ldots, d-1$, then $M\otimes_RN$ must be maximal Cohen-Macaulay, and this yields a contradiction due to~\ref{SRT}; see also Lemma~\ref{grz}. Therefore, $\hh^j_{\fm}(M\otimes_RN)\neq 0$ for some $j$ with $0\leq j \leq d-1$.

Next, we use Proposition~\ref{cx1} with $n=d$ and $\mz=\V(\fm)$, and conclude from part~(a) of the proposition that $\hh^0_{\fm}(M\otimes_RN) \cong \hh^2_{\fm}(M\otimes_RN)$, \ldots, $\hh^{d-3}_{\fm}(M\otimes_RN) \cong \hh^{d-1}_{\fm}(M\otimes_RN)$, and $\hh^{d-4}_{\fm}(M\otimes_RN) \cong \hh^{d-2}_{\fm}(M\otimes_RN)$. As we have already observed that $\hh^i_{\fm}(M\otimes_RN)$ does not vanish for all $i=0, \ldots, d-1$, the isomorphisms establish the claims.
\end{proof}

%%%%%%%%%%%%%%%%%%%%%%%%%%%%%%%%%%%%%%%%%%%%%%%%%%%%%%%%%%%%%%%%%
%%%%%%%%%%%%%%%%%%% SECTION  %%%%%%%%%%%%%%%%%%%%%%%%%%%%%%%%%%%%%%%%
%%%%%%%%%%%%%%%%%%% SECTION  %%%%%%%%%%%%%%%%%%%%%%%%%%%%%%%%%%%%%%%%
%%%%%%%%%%%%%%%%%%% SECTION  %%%%%%%%%%%%%%%%%%%%%%%%%%%%%%%%%%%%%%%%
%%%%%%%%%%%%%%%%%%%%%%%%%%%%%%%%%%%%%%%%%%%%%%%%%%%%%%%%%%%%%%%%%

\section{Proofs of the auxiliary results}\label{PAR}
In this section, we prove the auxiliary results used in the proof of Theorem~\ref{t12}. Specifically, we prove Lemmas~\ref{lem1},~\ref{pr},~\ref{grz}, as well as Proposition~\ref{l2}. We can directly provide proofs for the lemmas. However, to prove Proposition~\ref{l2}, we require an additional auxiliary result, namely Lemma~\ref{t1}, which appears to be of independent interest.

Recall that, unless otherwise stated, $R$ denotes a commutative Noetherian ring (not necessarily local). Moreover, unless otherwise stated, each $R$-module is assumed to be finitely generated. 

\begin{proof}[Proof of Lemma \ref{lem1}] 
Set $\mw=\Supp_R(M)$ and $\mz'=\mz\cap\mw$.
Then $\mw$ and $\mz'$ are specialization-closed subsets
of $\Spec(R)$. Set $T=M\otimes_RN$ and choose a minimal injective resolution
\[
0\longrightarrow T\longrightarrow I^0\longrightarrow I^1
\longrightarrow\cdots.
\]
Since the resolution is minimal, we have
\[
\Supp_R(I^j)\subseteq\Supp_R(T)
=\mw\cap\Supp_R(N)\subseteq\mw
\quad\text{for all }j\geq0.
\]
Hence $\Gamma_{\mw}(I^j)=I^j$ for all $j\geq0$;
see \ref{LC}(iv). Therefore, by~\ref{LC}(vii),
\[
\Gamma_{\mz'}(I^j)
=\Gamma_{\mz\cap\mw}(I^j)
=\Gamma_{\mz}\bigl(\Gamma_{\mw}(I^j)\bigr)
=\Gamma_{\mz}(I^j)
\quad\text{for all }j\geq0.
\]
These equalities identify
$\Gamma_{\mz'}(I^\bullet)$ and $\Gamma_{\mz}(I^\bullet)$
as complexes. Consequently,
\[
\hh^i_{\mz}(T)
\cong H^i\bigl(\Gamma_{\mz}(I^\bullet)\bigr)
=H^i\bigl(\Gamma_{\mz'}(I^\bullet)\bigr)
\cong\hh^i_{\mz'}(T)
\quad\text{for all }i\geq0.
\]
This completes the proof.
\end{proof}

\begin{proof}[Proof of Lemma~\ref{pr}] It follows that
\begin{equation}\label{pr1} \tag{\ref*{pr}.1}
\widehat{\Tor}_j^{R_\fp}(M_\fp,N_\fp)\cong\Tor_{j}^{R_\fp}(M_\fp,N_\fp) \text{ for all } \fp\in\X \text{ and all } j\geq 1,
\end{equation}
since $M_{\fp}$ is totally reflexive over $R_{\fp}$; see~\ref{Tate hom}(i). Moreover, given $\fp\in\X$ and $i\leq 0$, the following hold:
\begin{equation}\label{pr2} \tag{\ref*{pr}.2}
\widehat{\Tor}_i^{R_\fp}(M_\fp,N_\fp)\cong\widehat{\Ext}^{-i-1}_{R_\fp}(M^*_\fp,N_\fp)\cong\widehat{\Ext}^{-i+1}_{R_\fp}(\Tr_{R_\fp} M_\fp,N_\fp)\cong\Ext^{-i+1}_{R_\fp}(\Tr_{R_\fp} M_\fp,N_\fp).
\end{equation}
Here, in~\eqref{pr2}, the first isomorphism is due to~\ref{Tate hom}(ii), the second one follows as $M^{\ast} \cong \Omega_R^2\Tr_R M$, and the third isomorphism holds because $\Tr_{R_\fp} M_\fp$ is totally reflexive over $R_{\fp}$; see~\ref{HD}(i) and~\ref{Tate hom}(i). Now the assertion follows easily from~(\ref{pr}.1) and~(\ref{pr}.2).
\end{proof}

\begin{rmk} \label{g-m} Let $\fp \in \Spec(R)$ and let $\mz$ be a specialization-closed subset of $\Spec(R)$. Set $E=\IH_R(R/\fp)$.

Note that $\fq \in \Supp_R(E)$ if and only if $\fq \in \V(\fp)$. Assume $\fp \in \mz$. If $\fq \in \Supp_R(E)$, then $\fq \in \V(\fp)$ since $\Supp_R(E)\subseteq \V(\fp)$  so that $\mz \owns \fp \subseteq \fq \in \Spec(R)$, which implies $\fq \in \mz$. Thus, $\Supp_R(E)\subseteq \mz$ and hence $\Gamma_{\mz}(E)=E$; see~\ref{LC}(iv).

Assume $\Gamma_{\mz}(E)\neq 0$. Then, since $\Gamma_{\mz}(E)=\{x\in E : \Supp_R(Rx)\subseteq\mz\}$, there is $0\neq x\in E$ such that $\Supp_R(Rx)\subseteq\mz$. On the other hand, $\Ass_R(Rx)\subseteq\Ass_R(E)=\{\fp\}$. This shows that, if $\fp \notin \mz$, then $\Gamma_{\mz}(E)=0$. Consequently, we have:
\begin{equation}\notag{} 
\Gamma_{\mz}(E)=\left\lbrace
\begin{array}{cl}
E\ \ & \text{ if}\ \ \fp\in\mz\\
0\ \ & \text{ if }\ \ \fp\notin\mz.\\
\end{array}\right.
\end{equation}
\end{rmk}

Given an $R$-module $M$ and $\fp \in \Spec(R)$, we set $k(\fp)= R_{\fp} / \fp R_{\fp}$ and $\mu_i(\fp,M)=\rank_{k(\fp)}\Ext^i_{R_{\fp}}(k(\fp), M_{\fp})$ for all $i\geq 0$; see~\cite{BassU}.

\begin{proof}[Proof of Lemma \ref{grz}]
Recall that, given an $R$-module $M$ and a specialization-closed subset $\mz$ of $\Spec(R)$, we have by definition that
$\gr_R(\mz,M)=\inf\{i\geq0:\hh^i_{\mz}(M)\neq0\}$. We want to prove that
\[
\gr_R(\mz,M)
=\inf\{\depth_{R_{\fp}}(M_{\fp}):\fp\in\mz\}.
\]

Consider a minimal injective resolution of $M$, namely an exact
complex of the form
\[
0\longrightarrow M
\xrightarrow{\partial^{-1}} E^0
\xrightarrow{\partial^0} E^1
\xrightarrow{\partial^1} E^2
\xrightarrow{\partial^2}\cdots,
\]
where
\[
E^i=\bigoplus_{\fp\in\Spec(R)}
\IH_R(R/\fp)^{\mu_i(\fp,M)}.
\]
It follows from Remark~\ref{g-m} that
\[
\tag{\ref{grz}.1}
\Gamma_{\mz}(E^i)
=\bigoplus_{\fp\in\Spec(R)}
\Gamma_{\mz}\bigl(\IH_R(R/\fp)\bigr)^{\mu_i(\fp,M)}
\cong
\bigoplus_{\fp\in\mz}
\IH_R(R/\fp)^{\mu_i(\fp,M)}
\quad\text{for all }i\geq0.
\]

Set $n=\inf\{\depth_{R_{\fp}}(M_{\fp}):\fp\in\mz\}$.
By Lemma~\ref{lem1}, applied with $N=R$, and the convention
$\depth(0)=\infty$, we may replace $\mz$ by
$\mz\cap\Supp_R(M)$. So, we may assume that
$\mz\subseteq\Supp_R(M)$.

If $\mz=\varnothing$, then $n=\infty$ and
$\hh^i_{\mz}(M)=0$ for all $i\geq0$, so the desired equality holds.
Hence, we assume that $\mz\neq\varnothing$.
Since $M_{\fp}$ is a nonzero finitely generated
$R_{\fp}$-module for every $\fp\in\mz$, it follows that $0\leq n<\infty$.

For all $\fp\in\mz$, we have
$\depth_{R_{\fp}}(M_{\fp})\geq n$, and hence
$\Ext^i_{R_{\fp}}(k(\fp),M_{\fp})=0$ for all $i=0, \ldots, n-1$.  
This implies that $\mu_i(\fp,M)=0$ for all $\fp\in\mz$
and for all $i=0, \ldots, n-1$. 
By the definition of $n$, there exists $\fq\in\mz$ such that
$\depth_{R_{\fq}}(M_{\fq})=n$.
Thus $\Ext^n_{R_{\fq}}(k(\fq),M_{\fq})$, and hence
$\mu_n(\fq,M)$, are nonzero.
We deduce from~(\ref{grz}.1) that
\[
\tag{\ref{grz}.2}
\Gamma_{\mz}(E^i)=0
\quad\text{for all } i=0, \ldots, n-1
\quad \text{ and } \quad 
0\neq\IH_R(R/\fq)^{\mu_n(\fq,M)}
\hookrightarrow\Gamma_{\mz}(E^n).
\]

As $E^n$ occurs in a minimal injective resolution of $M$,
it is an essential extension of
$\im(\partial^{n-1})$; when $n=0$, this is the image of
the displayed injection $\partial^{-1}:M\to E^0$.
Thus
$\im(\partial^{n-1})\cap L\neq0$
for every nonzero $R$-submodule $L$ of $E^n$.
Since $\Gamma_{\mz}(E^n)\neq0$ by~(\ref{grz}.2),
we conclude that
\[
\tag{\ref{grz}.3}
\ker(\partial^n)\cap\Gamma_{\mz}(E^n)
=\im(\partial^{n-1})\cap\Gamma_{\mz}(E^n)
\neq0.
\]

Note that
$\hh^0_{\mz}(M)\cong\Gamma_{\mz}(M)
\hookrightarrow\Gamma_{\mz}(E^0)$
and
$\hh^i_{\mz}(M)
\cong
\ker\bigl(\Gamma_{\mz}(\partial^i)\bigr)/
\im\bigl(\Gamma_{\mz}(\partial^{i-1})\bigr)\;
\text{ for all }i\geq1$. Since $\ker\bigl(\Gamma_{\mz}(\partial^i)\bigr)
\subseteq\Gamma_{\mz}(E^i)$,
we see from~(\ref{grz}.2) that
$\hh^i_{\mz}(M)=0$ for all $i=0, \ldots, n-1$. Consequently,
\[
\gr_R(\mz,M)=\inf\{i\geq0:\hh^i_{\mz}(M)\neq0\}\geq n.
\]

It remains to prove that $\hh^n_{\mz}(M)\neq0$.
If $n\geq1$, then $\Gamma_{\mz}(E^{n-1})=0$
by~(\ref{grz}.2), so
$\im\bigl(\Gamma_{\mz}(\partial^{n-1})\bigr)=0$.
If $n=0$, left exactness gives $\hh^0_{\mz}(M)
\cong\Gamma_{\mz}(M)
\cong\ker\bigl(\Gamma_{\mz}(\partial^0)\bigr)$. Thus, in either case, left exactness and~(\ref{grz}.3) yield
\[
\begin{aligned}
\hh^n_{\mz}(M)\cong\ker\bigl(\Gamma_{\mz}(\partial^n)\bigr)\cong\Gamma_{\mz}\bigl(\ker(\partial^n)\bigr)=\ker(\partial^n)\cap\Gamma_{\mz}(E^n)
\neq0.
\end{aligned}
\]
This proves the desired equality.
\end{proof}

Next we prepare a lemma which is necessary for our proof of Proposition~\ref{l2}. 

\begin{lem}\label{t1} Let $\mz$ be a specialization-closed subset of $\Spec(R)$, $M$ be an $R$-module, and let $N$ be an $R$-module (not necessarily finitely generated). Set $g=\gr_R(\mz,N)$ and assume there is an integer $t$ such that $2\leq t\leq g$ and $\Supp_R\big(\Ext^i_R(M,N)\big)\subseteq\mz$ for all $i=1, \ldots, t-1$. Then the following hold:
\begin{enumerate}[\rm(i)]
\item $\Ext^i_R(M,N)\cong\hh_{\mz}^{i+1}\left(\Hom_R(M,N)\right)$ for all $i=1, \ldots, t-2$.
\item There is an injection $\Ext^{t-1}_R(M,N)\hookrightarrow\hh_{\mz}^{t}\big(\Hom_R(M,N)\big)$.
\end{enumerate}
\end{lem}

\begin{proof}
Consider a free resolution of $M$, namely an exact sequence
of the form
\[
\tag{\ref{t1}.1}
\cdots\longrightarrow F_1
\overset{\partial_1}{\longrightarrow}F_0
\overset{\partial_0}{\longrightarrow}M
\longrightarrow0,
\]
where $F_i$ is a finitely generated free $R$-module for all
$i\geq0$. Applying the functor
$(-)^{\triangledown}=\Hom_R(-,N)$
to~(\ref{t1}.1), we obtain the following complex:
\[
\tag{\ref{t1}.2}
0\longrightarrow M^{\triangledown}
\overset{\partial_0^{\triangledown}}{\longrightarrow}
(F_0)^{\triangledown}
\overset{\partial_1^{\triangledown}}{\longrightarrow}
(F_1)^{\triangledown}
\overset{\partial_2^{\triangledown}}{\longrightarrow}
(F_2)^{\triangledown}
\longrightarrow\cdots.
\]
Set
$X_i=\im(\partial_{i+1}^{\triangledown})$
and
$T_i=\coker(\partial_i^{\triangledown})$
for all $i\geq0$. The complex~(\ref{t1}.2) yields the
following short exact sequences:
\[
\tag{\ref{t1}.3}
\begin{gathered}
0\longrightarrow\Ext^i_R(M,N)
\longrightarrow T_i
\longrightarrow X_i
\longrightarrow0
\quad\text{for all }i\geq1,\\
0\longrightarrow X_j
\longrightarrow(F_{j+1})^{\triangledown}
\longrightarrow T_{j+1}
\longrightarrow0
\quad\text{for all }j\geq0.
\end{gathered}
\]

We assume that
$\Supp_R(\Ext^i_R(M,N))\subseteq\mz$
for all $i=1,\ldots,t-1$.
Hence, by~\ref{LC}(iv, v), we have
\[
\tag{\ref{t1}.4}
\begin{gathered}
\Ext^i_R(M,N)
=\Gamma_{\mz}\bigl(\Ext^i_R(M,N)\bigr),\\
\hh^j_{\mz}\bigl(\Ext^i_R(M,N)\bigr)=0
\quad\text{for all }i=1,\ldots,t-1
\text{ and all }j\geq1.
\end{gathered}
\]

Since each $(F_i)^{\triangledown}$ is a finite direct sum
of copies of $N$ and $\gr_R(\mz,N)\geq t$, we have
\[
\hh^j_{\mz}\bigl((F_i)^{\triangledown}\bigr)=0
\quad\text{for all }i\geq0
\text{ and all }j=0,\ldots,t-1.
\]
As $X_i\hookrightarrow(F_{i+1})^{\triangledown}$,
left exactness of $\Gamma_{\mz}(-)$ gives
$\Gamma_{\mz}(X_i)=0$ for all $i\geq0$.
Applying $\Gamma_{\mz}(-)$ to the first short exact sequence
in~(\ref{t1}.3), we obtain
\[
\Gamma_{\mz}\bigl(\Ext^i_R(M,N)\bigr)
\cong\Gamma_{\mz}(T_i)
\quad\text{for all }i\geq1.
\]
Consequently, by~(\ref{t1}.4), we see that
\[
\tag{\ref{t1}.5}
\Ext^i_R(M,N)
=\Gamma_{\mz}\bigl(\Ext^i_R(M,N)\bigr)
\cong\Gamma_{\mz}(T_i)
\quad\text{for all }i=1,\ldots,t-1.
\]

Similarly, the vanishing in~(\ref{t1}.4), together with
the first short exact sequence in~(\ref{t1}.3), yields
\[
\tag{\ref{t1}.6}
\hh_{\mz}^{j+1}(T_i)
\cong\hh_{\mz}^{j+1}(X_i)
\quad\text{for all }i=1,\ldots,t-1
\text{ and all }j\geq0.
\]

Since $(F_{i+1})^{\triangledown}$ is a finite direct sum
of copies of $N$ and $\gr_R(\mz,N)\geq t$, we have
$\hh^j_{\mz}\big((F_{i+1})^{\triangledown}\big)=0$
for all $i\geq0$ and all $j=0,\ldots,t-1$. Therefore, the second short exact sequence
in~(\ref{t1}.3) yields
\[
\tag{\ref{t1}.7}
\hh_{\mz}^j(T_{i+1})
\cong\hh_{\mz}^{j+1}(X_i)
\quad\text{for all }i\geq0
\text{ and all }j=0,\ldots,t-2.
\]
Using~(\ref{t1}.6) and~(\ref{t1}.7) repeatedly, we obtain
\[
\tag{\ref{t1}.8}
\Gamma_{\mz}(T_i)
\cong\hh_{\mz}^{i-1}(T_1)
\quad\text{for all }i=1,\ldots,t-1.
\]

Next, we consider the following short exact sequences
induced from~(\ref{t1}.2):
\[
\tag{\ref{t1}.9}
\begin{gathered}
0\longrightarrow M^{\triangledown}
\longrightarrow(F_0)^{\triangledown}
\longrightarrow X_0
\longrightarrow0,\\
0\longrightarrow X_0
\longrightarrow(F_1)^{\triangledown}
\longrightarrow T_1
\longrightarrow0.
\end{gathered}
\]

Since $F_0$ and $F_1$ are finitely generated free,
their duals $(F_0)^{\triangledown}$ and
$(F_1)^{\triangledown}$ are finite direct sums of copies
of $N$. Hence the assumption $\gr_R(\mz,N)\geq t$ gives
\[
\hh^j_{\mz}\bigl((F_0)^{\triangledown}\bigr)
=\hh^j_{\mz}\bigl((F_1)^{\triangledown}\bigr)=0
\quad\text{for all }j=0,\ldots,t-1.
\]
So, the long exact local cohomology sequences associated
to~(\ref{t1}.9) give the following injection and isomorphisms:
\[
\tag{\ref{t1}.10}
\begin{gathered}
\hh_{\mz}^{t-1}(X_0)
\hookrightarrow
\hh_{\mz}^{t}\bigl(M^{\triangledown}\bigr),\\
\hh_{\mz}^{i-1}(T_1)
\cong
\hh_{\mz}^{i+1}\bigl(M^{\triangledown}\bigr)
\quad\text{for all }i=1,\ldots,t-2.
\end{gathered}
\]
Consequently, we have the following isomorphisms,
which establish the first part of the lemma:
\[
\tag{\ref{t1}.11}
\Ext^i_R(M,N)
\cong\Gamma_{\mz}(T_i)
\cong\hh_{\mz}^{i-1}(T_1)
\cong\hh_{\mz}^{i+1}\bigl(M^{\triangledown}\bigr)
\quad\text{for all } i=1,\ldots,t-2.
\]
Here, in (\ref{t1}.11), the isomorphisms follow from~(\ref{t1}.5),~(\ref{t1}.8), and~(\ref{t1}.10), respectively.

Similarly, the second part of the lemma follows from
\[
\tag{\ref{t1}.12}
\begin{aligned}
\Ext^{t-1}_R(M,N) \cong\Gamma_{\mz}(T_{t-1})
\cong\hh_{\mz}^{t-2}(T_1)
\cong\hh_{\mz}^{t-1}(X_0)
\hookrightarrow
\hh_{\mz}^{t}\bigl(M^{\triangledown}\bigr).
\end{aligned}
\]
Here, in~(\ref{t1}.12), the first and second isomorphisms follow
from~(\ref{t1}.5) and~(\ref{t1}.8), respectively,
by setting $i=t-1$; the third isomorphism follows from~(\ref{t1}.7)
by setting $i=0$ and $j=t-2$; and the injection follows from~(\ref{t1}.10).
\end{proof}

We note the following result of Bruns-Vetter~\cite{BV}, which is similar to the first part of Lemma~\ref{t1}:

\begin{chunk} Let $R$ be a local ring, and let $M$ and $N$ be (finitely generated) $R$-modules. If $M$ is locally free on the punctured spectrum of $R$ and $\depth_R(N)=t\geq 3$, then $\Ext^i_R(M,N)\cong\hh_{\fm}^{i+1}(M^{\ast}\otimes_RN)$ for all $i=1, \ldots, t-2$; see~\cite[4.1]{BV}.
\end{chunk}
\iffalse
The conclusion of Lemma~\ref{t1} may fail if the grade of $N$ is strictly smaller than the integer $t\geq 2$:

\begin{eg} Let $R=k[\![x,y, z]\!]/(xy)$, $M=N=R/(z)$, and $\mz=\V(\fm)$. Then $\depth_R(N)=1$ and so $g=\gr_R(\mz,N)=1$. Note that $\hh_{\fm}^{2}\left(\Hom_R(M,N)\right)=\hh_{\fm}^{2}(M)=0$ since $\dim_R(M)=1$. Note also that $\Supp_R\big(\Ext^1_R(M,N)\big)=\Supp_R(M)=\V(z)\subseteq \mz$, and $\Ext^i_R(M,N)=0$ for all $i\geq 2$ since $\pd_R(M)=1$.

As $\Ext^1_R(M,N) \cong M \ncong 0=\hh_{\fm}^{2}\left(\Hom_R(M,N)\right)$, the conclusion of the first part of Lemma~\ref{t1} fails in case $t=3$. Also, as $\Ext^{1}_R(M,N)$ does not embed in $\hh_{\fm}^{2}\left(\Hom_R(M,N)\right)$, the conclusion of the second part of Lemma~\ref{t1} fails in case $t=2$.
\end{eg}
\fi
The injection established in the second part of Lemma~\ref{t1} does not need to be an isomorphism, in general:

\begin{eg} Let $R=k[\![x,y, z, u]\!]/(xy, zu)$ and $M=N=R/(y,u)$. Then $\dim_R(N)=\depth_R(N)=2$ and $\Hom_R(M,N)\cong M$. Set $\mz=\V(\fm)$ so that $\gr_R(\mz,N)=2$. Letting $g=t=2$ as in Lemma~\ref{t1}, it follows that $\Ext^{t-1}_R(M,N)=\Ext^{1}_R(M,N)=0$ and $\hh_{\fm}^{2}\left(\Hom_R(M,N)\right)\neq 0$. 
\end{eg}

We present a consequence of Lemma~\ref{t1}:

\begin{cor} Let $R$ be a $d$-dimensional local complete intersection ring with an isolated singularity, with $d\geq 2$, and let $M$ be a maximal Cohen-Macaulay $R$-module. If $\hh_{\fm}^{t}\left(\Hom_R(M,M)\right)=0$ for some odd integer $n$ such that $2\leq t \leq d$, then $M$ is free.
\end{cor} 

\begin{proof}
We may assume that $M\neq0$.
Setting $\mz=\V(\fm)$, we have
$\gr_R(\mz,M)=\depth_R(M)=d$.
Since $R$ has an isolated singularity, $M$ is locally free
on the punctured spectrum of $R$. Hence, 
$\Supp_R\bigl(\Ext^i_R(M,M)\bigr)\subseteq\mz$
for all $i\geq1$.
Therefore, Lemma~\ref{t1}(ii) gives an injection
\[
\Ext^{t-1}_R(M,M)
\hookrightarrow
\hh_{\fm}^{t}\left(\Hom_R(M,M)\right)=0.
\]
As $t-1$ is a positive even integer, it follows
from~\ref{CIfacts}(i) that $\pd_R(M)<\infty$.
Since $M$ is maximal Cohen-Macaulay, the
Auslander-Buchsbaum formula implies that $M$ is free.
\end{proof}

We conclude this section with a proof of Proposition~\ref{l2}.

\begin{proof}[Proof of Proposition \ref{l2}] We assume there is an integer $n$ such that $0\leq n\leq\gr_R(\mz,N)$ and $\Supp_R\left(\Ext^i_R(\Tr_R M,N)\right)\subseteq\mz$ for all $i=1, \ldots, n+1$. Consider the four-term exact sequence~\ref{lan}:
\[
\tag{\ref{l2}.1} 0\rightarrow\Ext^1_R(\Tr_R M,N)\rightarrow M\otimes_RN\rightarrow\Hom_R(M^*,N)\rightarrow\Ext^2_R(\Tr_R M,N)\rightarrow0.
\]
Note that, since $\Supp_R\left(\Ext^i_R(\Tr_R M,N)\right)\subseteq\mz$, it follows, for all $i=1, \ldots, n+1$, from \ref{LC}(iv, v) that:
\[
\tag{\ref{l2}.2} \Gamma_{\mz}\left(\Ext^i_R(\Tr_R M,N)\right)=\Ext^i_R(\Tr_R M,N) \;\;\text{ and }\;\;  \hh^j_{\mz}\left(\Ext^i_R(\Tr_R M,N)\right)=0 \text{ for all } j\geq 1.
\]
Therefore, since $\Gamma_{\mz}\left(\Ext^1_R(\Tr_R M,N)\right)=\Ext^1_R(\Tr M,N)$, by applying the functor $\Gamma_{\mz}(-)$ to~(\ref{l2}.1), we obtain an injection $\Ext^{1}_R(\Tr M,N)\hookrightarrow\Gamma_{\mz}(M\otimes_RN)$. This proves the second part of the proposition (note that the first part of the proposition is vacuous if $n=0$). Hence, we assume $n\geq 1$.

The four-term exact sequence~(\ref{l2}.1) yields two short exact sequences
\[
\tag{\ref{l2}.3} 0\rightarrow\Ext^1_R\left(\Tr_RM,N\right)\rightarrow M\otimes_RN\rightarrow X\rightarrow0,
\]
and
\[
\tag{\ref{l2}.4} 0\rightarrow X\rightarrow\Hom_R(M^*,N)\rightarrow\Ext^2_R\left(\Tr_RM,N\right)\rightarrow0.
\]

We have, by Lemma \ref{gr}, that: 
\[
\tag{\ref{l2}.5} \gr_R\left(\mz,\Hom_R(M^{\ast},N)\right) \geq \min \left\{2,\gr_R(\mz,N)\right\}.
\]
As $\gr_R(\mz,N)\geq n\geq 1$, we see from (\ref{l2}.5) that $\gr_R\left(\mz,\Hom_R(M^{\ast},N)\right) \geq 1$. Thus, $\gr_R(\mz,X)\geq 1$ since 
$X\hookrightarrow \Hom_R(M^*,N)$. Therefore, $\Gamma_{\mz}(X)=0=\Gamma_{\mz}\left(\Hom_R(M^{\ast},N)\right)$; see~\ref{gsc}. We  apply the functor $\Gamma_{\mz}(-)$ to~(\ref{l2}.3) and obtain by~(\ref{l2}.2) that:
\[
\tag{\ref{l2}.6} \Ext^1_R(\Tr_R M,N)=\Gamma_{\mz}\left(\Ext^1_R(\Tr M,N)\right) \cong \Gamma_{\mz}(M\otimes_RN) \;\;\text{ and }\;\; \hh^i_{\mz}(M\otimes_RN) \cong \hh^i_{\mz}(X) \text{ for all } i\geq 1.
\]
Similarly, applying the functor $\Gamma_{\mz}(-)$ to~(\ref{l2}.4), we obtain by~(\ref{l2}.2) the exact sequence:
\[
\tag{\ref{l2}.7} 0 \rightarrow \Ext^2_R(\Tr_R M,N)=\Gamma_{\mz}\left(\Ext^2_R(\Tr M,N)\right)\hookrightarrow\hh^1_{\mz}(X) \rightarrow \hh^1_{\mz}\left(\Hom_R(M^*,N)\right),   
\]
and the isomorphisms:
\[
\tag{\ref{l2}.8} \hh^{i}_{\mz}(X) \cong \hh^i_{\mz}\left(\Hom_R(M^*,N)\right)  \text{ for all } i\geq 2.
\]
Note that, we now have $\Gamma_{\mz}(M\otimes_RN) \cong \Ext^1_R(\Tr_R M,N)$ and $ \Ext^2_R(\Tr_R M,N) \hookrightarrow \hh^1_{\mz}(X) \cong \hh^1_{\mz}(M\otimes_RN)$; see~(\ref{l2}.6) and~(\ref{l2}.7). So, both of the claims of the proposition follow for the case where $n=1$.

Next assume $n\geq 2$. Recall that $M^{\ast} \cong \Omega^2\Tr_RM$; see \ref{a1}. So, by the hypothesis of the proposition, we have that:
\[
\tag{\ref{l2}.9} \Supp_R\left(\Ext^i_R(M^{\ast},N)\right)=\Supp_R\left(\Ext^{i+2}_R(\Tr_R M,N)\right)\subseteq\mz \text{ for all } i=1, \ldots, n-1. 
\]
Therefore, as $\gr_R(\mz,N)\geq n\geq 2$, we may use Lemma \ref{t1} for the pair $(M^{\ast}, N)$ with $t=n$, and conclude:
\[
\tag{\ref{l2}.10} \Ext^{n-1}_R(M^{\ast},N)\hookrightarrow\hh^n_{\mz}\left(\Hom_R(M^{\ast},N)\right).
\]
Consequently, we obtain:
\[
\tag{\ref{l2}.11} \Ext^{n+1}_R(\Tr M,N)\cong \Ext^{n-1}_R(M^{\ast},N)\hookrightarrow\hh^n_{\mz}\left(\Hom_R(M^{\ast},N)\right) \cong \hh^{n}_{\mz}(X) \cong \hh^n_{\mz}(M\otimes_RN). 
\]
Here, in~(\ref{l2}.11), the second and third isomorphisms follow from~(\ref{l2}.8) and~(\ref{l2}.6), respectively (recall that $n\geq 2$). This establishes the second part of the proposition.

We proceed to prove the first part of the proposition when $n\geq 2$. We first consider the case where $n=2$. Note, since $\gr_R(\mz,N)\geq n$, it follows from~(\ref{l2}.5) that $\gr_R\left(\mz,\Hom_R(M^{\ast},N)\right)\geq2$. Thus, $\hh^1_{\mz}\left(\Hom_R(M^{\ast},N)\right)=0$. Hence,~(\ref{l2}.6) and~(\ref{l2}.7) yield the following isomorphisms, which complete the proof of the proposition for $n=2$:
\[
\Ext^2_R(\Tr_R M,N) \cong \hh^1_{\mz}(X) \cong \hh^1_{\mz}(M\otimes_RN).
\]

Finally, we assume $n\geq 3$. In that case, due to (\ref{l2}.9), we may use the first part of Lemma~\ref{t1} for the pair $(M^{\ast}, N)$ with $t=n$, and conclude
\[
\tag{\ref{l2}.12} \hh_{\mz}^{i+1}\left(\Hom_R(M^{\ast},N)\right) \cong \Ext^i_R(M^{\ast},N)  \; \text{ for all } i=1, \ldots, n-2.
\]
So, if $j$ is an integer such that $2\leq j \leq n-1$, we have
\[
\tag{\ref{l2}.13}  \hh_{\mz}^{j}\left(M\otimes_RN\right) \cong \hh_{\mz}^{j}\left(\Hom_R(M^{\ast},N)\right) \cong \Ext^{j-1}_R(M^{\ast},N) \cong \Ext^{j+1}_R\left(\Tr_R M,N\right).
\]
In~(\ref{l2}.13), the first isomorphism holds due to~(\ref{l2}.6) and~(\ref{l2}.8) since $j\geq 2$, the second isomorphism follows from~(\ref{l2}.12), and the third isomorphism follows by the fact that $M^{\ast} \cong \Omega^2\Tr_RM$. Consequently,~(\ref{l2}.13) confirms the first assertion of the proposition.
\end{proof}

%%%%%%%%%%%%%%%%%%%%

\section*{Acknowledgements}
The authors thank to K\k {e}stutis \v{C}esnavi\v{c}ius, Brian Laverty, Hiroki Matsui, and Takeshi Yoshizawa for their valuable suggestions and comments on earlier versions of the manuscript. The authors also thank the anonymous referee for checking the manuscript.

Part of this work was carried out during Celikbas's visit to the Max Planck Institute for Mathematics in Bonn, Germany, from January to June 2025. The author acknowledges the financial support and the hospitality provided by the institute.

\medskip
\noindent\textbf{AI assistance.} The results and original proofs of this paper were obtained by the authors without AI assistance. After the paper was accepted for publication, the authors used ChatGPT (OpenAI, GPT-6.1) to assist with language editing, typographical and bibliographic corrections, and checking the manuscript for mathematical errors and inconsistencies. All changes incorporated into the final version were reviewed and verified by the authors, who take full responsibility for the mathematical correctness and entire content of the paper.
 
%%%%%%%%%%%%%%%%%%


\begin{thebibliography}{10}

\bibitem{Tome1}
\emph{Anneaux de {G}orenstein, et torsion en alg\`ebre commutative},
  S\'eminaire d'Alg\`ebre Commutative dirig\'e par Pierre Samuel, 1966/67.
  Texte r\'edig\'e, d'apr\`es des expos\'es de Maurice Auslander, Marguerite
  Mangeney, Christian Peskine et Lucien Szpiro. \'Ecole Normale Sup\'erieure de
  Jeunes Filles, Secr\'etariat math\'ematique, Paris, 1967.

\bibitem{MA}
Mohsen Asgharzadeh, \emph{Finite support of tensor products}, Preprint; posted
  at arXiv:1902.10509 (2019).

\bibitem{AOM}
Mohsen Asgharzadeh, Olgur Celikbas, and Arash Sadeghi, \emph{A study of the
  cohomological rigidity property}, Preprint; posted at arXiv:2009.06481
  (2020).

\bibitem{Au}
Maurice Auslander, \emph{Modules over unramified regular local rings}, Illinois
  J. Math. \textbf{5} (1961), 631--647.

\bibitem{AuBr}
Maurice Auslander and Mark Bridger, \emph{Stable module theory}, Memoirs of the
  American Mathematical Society, No. 94, American Mathematical Society,
  Providence, R.I., 1969.

\bibitem{Av1}
Luchezar~L. Avramov, \emph{Modules of finite virtual projective dimension},
  Invent. Math. \textbf{96} (1989), no.~1, 71--101.

\bibitem{Av2}
\bysame, \emph{Infinite free resolutions}, six lectures on commutative algebra
  ({B}ellaterra, 1996), Progr. Math., vol. 166, Birkh\"auser, Basel, 1998,
  pp.~1--118.

\bibitem{AvBu}
Luchezar~L. Avramov and Ragnar-Olaf Buchweitz, \emph{Support varieties and
  cohomology over complete intersections}, Invent. Math. \textbf{142} (2000),
  no.~2, 285--318.

\bibitem{AGP}
Luchezar~L. Avramov, Vesselin~N. Gasharov, and Irena~V. Peeva, \emph{Complete
  intersection dimension}, Inst. Hautes \'Etudes Sci. Publ. Math. (1997),
  no.~86, 67--114 (1998).

\bibitem{AM}
Luchezar~L. Avramov and Alex Martsinkovsky, \emph{Absolute, relative, and
  {T}ate cohomology of modules of finite {G}orenstein dimension}, Proc. London
  Math. Soc. (3) \textbf{85} (2002), no.~2, 393--440.

\bibitem{BassU}
Hyman Bass, \emph{On the ubiquity of {Gorenstein} rings}, Math. Z. \textbf{82}
  (1963), 8--28 (English).

\bibitem{Bsh}
Markus~P. Brodmann and Rodney~Y. Sharp, \emph{Local cohomology. {An} algebraic
  introduction with geometric applications}, 2nd ed., Camb. Stud. Adv.
  Math., vol. 136, Cambridge: Cambridge University Press, 2012 (English).

\bibitem{BH}
Winfried Bruns and J\"{u}rgen Herzog, \emph{Cohen-{M}acaulay rings}, Cambridge
  Studies in Advanced Mathematics, vol.~39, 2nd ed. Cambridge University Press,
  Cambridge, 1998.

\bibitem{BV}
Winfried Bruns and Udo Vetter, \emph{Length formulas for the local cohomology
  of exterior powers}, Math. Z. \textbf{191} (1986), 145--158 (English).

\bibitem{Ce}
Olgur Celikbas, \emph{Vanishing of {T}or over complete intersections}, J.
  Commut. Algebra \textbf{3} (2011), no.~2, 169--206.

\bibitem{CeD}
Olgur Celikbas and Hailong Dao, \emph{Asymptotic behavior of {E}xt functors for
  modules of finite complete intersection dimension}, Math. Z. \textbf{269}
  (2011), no.~3, 1005--1020.

\bibitem{CSTT}
Olgur Celikbas, Shiro Goto, Ryo Takahashi, and Naoki Taniguchi, \emph{On the
  ideal case of a conjecture of {H}uneke and {W}iegand}, Proc. Edinb. Math.
  Soc. (2) \textbf{62} (2019), no.~3, 847--859.

\bibitem{CIPW}
Olgur Celikbas, Srikanth Iyengar, Greg Piepmeyer, and Roger Wiegand,
  \emph{Criteria for vanishing of {T}or over complete intersections}, Pacific
  J. Math. \textbf{276} (2015), no.~1, 93--115.

\bibitem{CeTos}
Olgur Celikbas and Toshinori Kobayashi, \emph{On a class of {B}urch ideals and
  a conjecture of {H}uneke and {W}iegand}, Collect. Math. \textbf{73} (2022),
  no.~2, 221--236.

\bibitem{Fukano}
Olgur Celikbas, Uyen Le, Hiroki Matsui, and Arash Sadeghi, \emph{Remarks on a
  conjecture of {Huneke} and {Wiegand} and the vanishing of (co)homology}, J.
  Math. Soc. Japan \textbf{76} (2024), no.~3, 951--974 (English).

\bibitem{AOMCM}
Olgur Celikbas and Arash Sadeghi, \emph{Maximal {C}ohen-{M}acaulay tensor
  products}, Ann. Mat. Pura Appl. (4) \textbf{200} (2021), no.~3, 923 -- 944.

\bibitem{CST}
Olgur Celikbas, Arash Sadeghi, and Ryo Takahashi, \emph{Bounds on depth of
  tensor products of modules}, J. Pure Appl. Algebra \textbf{219} (2015),
  no.~5, 1670--1684.

\bibitem{CeRo}
Olgur Celikbas and Roger Wiegand, \emph{Vanishing of {T}or, and why we care
  about it}, J. Pure Appl. Algebra \textbf{219} (2015), no.~3, 429--448.

\bibitem{KC}
K{\c{e}}stutis {\v{C}}esnavi{\v{c}}ius, \emph{Grothendieck-{Lefschetz} for
  vector bundles}, Algebr. Geom. \textbf{7} (2020), no.~4, 503--511 (English).

\bibitem{Gdimbook}
Lars~Winther Christensen, \emph{Gorenstein dimensions}, Lecture Notes in
  Mathematics, vol. 1747, Springer-Verlag, Berlin, 2000.

\bibitem{cfh}
Lars~Winther Christensen, Anders Frankild, and Henrik Holm, \emph{On
  {G}orenstein projective, injective and flat dimensions -- a functorial
  description with applications}, J. Algebra \textbf{302} (2006), no.~1,
  231--279.

\bibitem{CJ}
Lars~Winther Christensen and David~A. Jorgensen, \emph{Vanishing of {T}ate
  homology and depth formulas over local rings}, J. Pure Appl. Algebra
  \textbf{219} (2015), no.~3, 464--481.

\bibitem{CHKV}
Alberto Corso, Craig Huneke, Daniel Katz, and Wolmer~V. Vasconcelos,
  \emph{Integral closure of ideals and annihilators of homology}, Commutative
  algebra, Lect. Notes Pure Appl. Math., vol. 244, Chapman \& Hall/CRC, Boca
  Raton, FL, 2006, pp.~33--48.

\bibitem{Da2}
Hailong Dao, \emph{Asymptotic behaviour of {T}or over complete intersections
  and applications}, Preprint; posted at arXiv:0710.5818 (2007).

\bibitem{Da3}
\bysame, \emph{Some observations on local and projective hypersurfaces}, Math.
  Res. Lett. \textbf{15} (2008), no.~2, 207--219.

\bibitem{Da1}
\bysame, \emph{Decent intersection and {T}or-rigidity for modules over local
  hypersurfaces}, Trans. Amer. Math. Soc. \textbf{365} (2013), no.~6,
  2803--2821.

\bibitem{DaoSurvey}
\bysame, \emph{Some homological properties of modules over a complete
  intersection, with applications}, Commutative algebra. Expository papers
  dedicated to David Eisenbud on the occasion of his 65th birthday, New York,
  NY: Springer, 2013, pp.~335--371 (English).

\bibitem{Ei}
David Eisenbud, \emph{Homological algebra on a complete intersection, with an
  application to group representations}, Trans. Amer. Math. Soc. \textbf{260}
  (1980), no.~1, 35--64.

\bibitem{EG}
E.~Graham Evans and Phillip Griffith, \emph{Syzygies}, London Mathematical
  Society Lecture Note Series, vol. 106, Cambridge University Press, Cambridge,
  1985.

\bibitem{GotoUlrich}
Shiro Goto, Kazuho Ozeki, Ryo Takahashi, Kei-Ichi Watanabe, and Ken-Ichi
  Yoshida, \emph{Ulrich ideals and modules}, Math. Proc. Cambridge Philos. Soc.
  \textbf{156} (2014), no.~1, 137--166.

\bibitem{sga2}
Alexander Grothendieck, \emph{S{\'e}minaire de g{\'e}om{\'e}trie alg{\'e}brique
  du {Bois} {Marie} 1962. {Cohomologie} locale des faisceaux coh{\'e}rents et
  th{\'e}or{\`e}mes de {Lefschetz} locaux et globaux ({SGA} 2)}, revised and
  updated edition of the 1968 original published by {North} {Holland} ed., Doc.
  Math. (SMF), vol.~4, Paris: Soci{\'e}t{\'e} Math{\'e}matique de France, 2005
  (French).

\bibitem{CAP}
Melvin Hochster, Craig Huneke, and Judy~D. Sally (eds.), \emph{Commutative
  {A}lgebra, {P}roceedings of a {M}icroprogram}, Mathematical Sciences Research
  Institute Publications, no.~1, Springer New York, NY, 1989.

\bibitem{HJW}
Craig Huneke, David~A. Jorgensen, and Roger Wiegand, \emph{Vanishing theorems
  for complete intersections}, J. Algebra \textbf{238} (2001), no.~2, 684--702.

\bibitem{HW1}
Craig Huneke and Roger Wiegand, \emph{Tensor products of modules and the
  rigidity of {${\textrm Tor}$}}, Math. Ann. \textbf{299} (1994), no.~3,
  449--476.

\bibitem{HW2}
\bysame, \emph{Tensor products of modules, rigidity and local cohomology},
  Math. Scand. \textbf{81} (1997), no.~2, 161--183.

\bibitem{Jo2}
David~A. Jorgensen, \emph{Tor and torsion on a complete intersection}, J.
  Algebra \textbf{195} (1997), no.~2, 526--537.

\bibitem{Jo1}
\bysame, \emph{Complexity and {T}or on a complete intersection}, J. Algebra
  \textbf{211} (1999), no.~2, 578--598.

\bibitem{J2}
\bysame, \emph{A generalization of the {A}uslander-{B}uchsbaum formula}, J.
  Pure Appl. Algebra \textbf{144} (1999), no.~2, 145--155.

\bibitem{JorExt}
\bysame, \emph{Finite projective dimension and the vanishing of {Ext}{{\(_R(M,
  M)\)}}}, Commun. Algebra \textbf{36} (2008), no.~12, 4461--4471 (English).

\bibitem{Jot}
P.~Jothilingam, \emph{A note on grade}, Nagoya Math. J. \textbf{59} (1975),
  149--152 (English).

\bibitem{Kem}
George~R. Kempf, \emph{A criterion for the splitting of a vector bundle}, Forum
  Math. \textbf{2} (1990), no.~5, 477--480.

\bibitem{Li}
Stephen Lichtenbaum, \emph{On the vanishing of {T}or in regular local rings},
  Illinois J. Math. \textbf{10} (1966), 220--226.

\bibitem{Kum}
N. Mohan Kumar, Prabhakar Rao, and Ravindra Girivaru,
\emph{Arithmetically Cohen-Macaulay bundles on hypersurfaces},
Commentarii Mathematici Helvetici \textbf{82} (2007),
no. 4, 829--843.

\bibitem{Luk}
Hing~Sun Luk and Stephen S.-T. Yau, \emph{Cohomology and splitting criterion
  for holomorphic vector bundles on $\textbf{CP}^n$}, Math. Nachr. \textbf{161}
  (1993), 233--238.

\bibitem{Mat}
Hideyuki Matsumura, \emph{Commutative ring theory}, second ed., Cambridge
  Studies in Advanced Mathematics, vol.~8, Cambridge University Press,
  Cambridge, 1989, Translated from the Japanese by M. Reid.

\bibitem{Mi}
Claudia Miller, \emph{Complexity of tensor products of modules and a theorem of
  {H}uneke-{W}iegand}, Proc. Amer. Math. Soc. \textbf{126} (1998), no.~1,
  53--60.

\bibitem{Roit}
Joseph~J. Rotman, \emph{An {I}ntroduction to {H}omological {A}lgebra}, second
  ed., Universitext, Springer, New York, 2009.

\bibitem{ST2}
Arash Sadeghi and Ryo Takahashi, \emph{Resolving subcategories closed under
  certain operations and a conjecture of {D}ao and {T}akahashi}, Michigan Math.
  J. \textbf{70} (2021), no.~2, 341--367.

\bibitem{TTY}
Ryo Takahashi, Yuji Yoshino, and Takeshi Yoshizawa, \emph{Local cohomology
  based on a nonclosed support defined by a pair of ideals}, J. Pure Appl.
  Algebra \textbf{213} (2009), no.~4, 582--600 (English).

\bibitem{Yo}
Yuji Yoshino, \emph{Cohen-{M}acaulay modules over {C}ohen-{M}acaulay rings},
  London Mathematical Society Lecture Note Series, vol. 146, Cambridge
  University Press, Cambridge, 1990.

\end{thebibliography}
\end{document}